\pdfoutput=1
\RequirePackage{ifpdf}
\ifpdf 
\documentclass[pdftex]{sigma}
\else
\documentclass{sigma}
\fi

\usepackage{mathrsfs}

\newtheorem{Theorem}{Theorem}[section]
\newtheorem{Corollary}[Theorem]{Corollary}
\newtheorem{Lemma}[Theorem]{Lemma}
\newtheorem{Proposition}[Theorem]{Proposition}
 { \theoremstyle{definition}
\newtheorem{Definition}[Theorem]{Definition}
\newtheorem{Remark}[Theorem]{Remark} }

\numberwithin{equation}{section}

\newcommand{\SF}{\ensuremath{ {S\!\!\!\;F}}}

\begin{document}

\allowdisplaybreaks

\renewcommand{\PaperNumber}{049}

\FirstPageHeading

\ShortArticleName{Hermite and Laguerre Symmetric Functions}

\ArticleName{Hermite and Laguerre Symmetric Functions\\ Associated with Operators\\ of Calogero--Moser--Sutherland Type}

\Author{Patrick DESROSIERS~$^\dag$ and Martin HALLN\"AS~$^\ddag$}

\AuthorNameForHeading{P.~Desrosiers and M.~Halln\"as}

\Address{$^\dag$~Instituto Matem\'atica y F\'isica, Universidad de Talca, 2 Norte 685, Talca, Chile}
\EmailD{\href{mailto:Patrick.Desrosiers@inst-mat.utalca.cl}{Patrick.Desrosiers@inst-mat.utalca.cl}}

\Address{$^\ddag$~Department of Mathematical Sciences, Loughborough University,\\
\hphantom{$^\ddag$}~Leicestershire, LE11 3TU, UK}
\EmailD{\href{mailto:M.A.Hallnas@lboro.ac.uk}{M.A.Hallnas@lboro.ac.uk}}

\ArticleDates{Received March 22, 2012, in f\/inal form July 25, 2012; Published online August 03, 2012}

\Abstract{We introduce and study natural generalisations of the Hermite and Laguerre polynomials in the ring of symmetric functions as eigenfunctions of inf\/inite-dimensional analogues of partial dif\/ferential operators of Calogero--Moser--Sutherland (CMS) type. In particular, we obtain generating functions, duality relations, limit transitions from Jacobi symmetric functions, and Pieri formulae, as well as the integrability of the corresponding operators. We also determine all ideals in the ring of symmetric functions that are spanned by either Hermite or Laguerre symmetric functions, and by restriction of the corresponding inf\/inite-dimensional CMS operators onto quotient rings given by such ideals we obtain so-called deformed CMS operators. As a consequence of this restriction procedure, we deduce, in particular, inf\/inite sets of polynomial eigenfunctions, which we shall refer to as super Hermite and super Laguerre polynomials, as well as the integrability, of these deformed CMS operators. We also introduce and study series of a generalised hypergeometric type, in the context of both symmetric functions and `super' polynomials.}


\Classification{05E05; 13J05; 81R12}

\tableofcontents

\section{Introduction}
The main purpose of this paper is to introduce and study two (non-homogenous) bases for the graded ring of symmetric functions $\Lambda$, which we shall refer to as Hermite and Laguerre symmetric functions due to their relations with the corresponding classical orthogonal polynomials.

\looseness=1
One of the main themes in the theory of symmetric functions is the description of va\-rious homogeneous bases for $\Lambda$. Important examples include the Schur symmetric functions and, more generally, Jack's and Macdonald's symmetric functions, which are one- and two-parameter deformations thereof; see, e.g., Macdonald~\cite{Mac95}. These symmetric functions arise as (inverse) limits with respect to the number of variables $n$ of corresponding symmetric polynomials.

In the early 1990s, Lassalle \cite{Las91c, Las91a,Las91b} and Macdonald \cite{Mac} independently introduced natural $n$-variable generalisations of the classical orthogonal polynomials (of Hermite, Laguerre and Jacobi), depending on one `extra' parameter~$\alpha$.\footnote{The generalised Jacobi polynomials were  previously introduced by Debiard \cite{Deb87}, and independently by Heckman and Opdam~\cite{HO87} (see also~\cite{Op89}).   However, the authors in~\cite{Deb87,HO87}  do not explicitly make use of Jack polynomial theory, which contrasts with the approach of Lassalle and Macdonald~\cite{Las91a,Mac}.} For the special value $\alpha=2$ they are naturally realised as functions on the $n\times n$ real symmetric matrices, and as such had been introduced already in work of James \cite{Jam75} in the Hermite, Herz \cite{Her55} and Constantine \cite{Con66} in the Laguerre, and James and Constantine~\cite{JC74} in the Jacobi case; see also Muirhead~\cite{Mui82}. In contrast to the Schur polynomials, these $n$-variable generalisations of the classical orthogonal polynomials do not posses limits as $n$ goes to inf\/inity.

We shall circumvent this problem by introducing an additional (formal) parameter $p_0$, which can be viewed as a zeroth-order power sum symmetric function, and consider symmetric functions over the f\/ield $\mathbb{Q}(p_0)$. The Hermite and Laguerre symmetric functions are then the unique elements in $\Lambda$ such that if we set $p_0=n\in\mathbb{N}$, and restrict to $\Lambda_n$, the ring of symmetric polynomials in $n$ variables, then we recover the corresponding $n$-variable polynomials of Lassalle and Macdonald.

\subsection{Description of main results}
1.~The Hermite and Laguerre symmetric functions will be indexed by partitions, and def\/ined as eigenfunctions of certain inf\/inite-dimensional analogues of partial dif\/ferential operators of Calogero--Moser--Sutherland (CMS) type; see Def\/initions \ref{HermiteDef} and \ref{LaguerreDef}. When restricted to a~f\/inite number of variables, they reduce to precisely such partial dif\/ferential operators, a fact that directly yields the relation to the multivariable Hermite and Laguerre polynomials introduced by Lassalle and Macdonald.

2.~In Propositions \ref{HermiteDualityProp} and \ref{LaguerreDualityProp} we establish particular duality relations, which are not present at the level of symmetric polynomials. For Jack's symmetric functions the corresponding duality is well-known, and amounts to taking the conjugate of the labelling partition and sending $\alpha\to1/\alpha$; see Section~VI.10 in Macdonald~\cite{Mac95}.

 3.~We show that a number of well-known and important results on, and properties of, the $n$-variable Hermite and Laguerre polynomials generalise in a straightforward manner to the symmetric functions setting. These include generating functions, limit transitions from Jacobi symmetric functions, recently introduced by Sergeev and Veselov \cite{SV09}, as well as higher-order eigenoperators.

\looseness=1
4.~We also establish Pieri formulae: the expansion of products of algebraic generators and linear basis elements of $\Lambda$ in the same basis elements. They are exhibited in Propositions~\ref{HermitePieriFormulaeProp} and~\ref{LaguerrePieriFormulaeProp}. In this case, the algebraic generators are the elementary symmetric functions. In the $n$-variable case such formulae were obtained by van Diejen \cite{vD97}. However, his formulae do not directly lift to the level of symmetric functions, since they depend in a non-trivial manner on~$n$.

5.~These Pieri formulae allow us to completely describe the set of ideals in $\Lambda$ that are inva\-riant under the full algebra of eigenoperators of either the Hermite or the Laguerre symmetric functions. In the Hermite case we show that such an ideal exists only if we set $p_0=n-\alpha m$ for some $n,m\in\mathbb{N}_0\equiv\mathbb{N}\cup\lbrace 0\rbrace$. In that case the ideal is unique, and we give a basis in terms of Hermite symmetric functions; see Theorem~\ref{IdealThm}. In the Laguerre case such an ideal exists also for $p_0=n+1-\alpha(m+a+1)$, which again is unique, and we provide a~basis in terms of Laguerre symmetric functions; see Theorem~\ref{IdealThmLaguerre}.

6.~The restriction of a dif\/ferential operator on $\Lambda$ onto a given quotient ring $\Lambda/I$ is possible if and only if the ideal $I$ is invariant under the operator in question. For $p_0=n-\alpha m$, we give an explicit realisation of the restrictions of the eigenoperators of both the Hermite and the Laguerre symmetric functions onto the corresponding quotient ring. More precisely, in Sections~\ref{sHermitePols} and~\ref{sLaguerrePols} we show that they are given by particular partial dif\/ferential operators of so-called deformed CMS type. The operators in question were previously considered by Feigin~\cite{Fei08}, Halln\"as and Langman~\cite{HL10}.  They also appeared, in a disguised form, in the work of Guhr and Kohler~\cite{KG05}, and more recently in the context of Random Matrix Theory~\cite{DL11}.    For $n=0$ or $m=0$ these operators reduce to ordinary CMS operators.

7.~Through this restriction procedure our previous results in the paper immediately yield corresponding results for these deformed CMS operators. In particular, the Hermite and Laguerre symmetric functions restrict to corresponding eigenfunctions that are polyno\-mials, which, following previous results in the literature, we shall refer to as super Hermite and super Laguerre polynomials, respectively.

8.~In order to establish a wider context for our results, we also discuss the notion of dif\/ferential operators on the ring of symmetric functions in Appendix~\ref{AppDiffOps}. In particular, we establish an explicit description of such operators in terms of their action on power sum symmetric functions. This will make it clear that the inf\/inite-dimensional CMS operators mentioned above are indeed dif\/ferential operators on the symmetric functions. In Appendix~\ref{AppCMSOps} we isolate certain results that hold true not only in the Hermite and Laguerre cases. This includes the fact that a generic inf\/inite-dimensional CMS operator of second order has a~complete set of eigenfunctions in the ring of symmetric functions.

One of our main motivations for this paper, and perhaps its most important consequence, is that lifting multivariable Hermite and Laguerre polynomials to the level of symmetric functions unif\/ies the CMS operators in question with their corresponding deformed analogues, as mentioned under item~(6). In particular, this provides a conceptual understanding of the latter operators, and many of their key properties are thus inherited from the undeformed case. A~further appraisal of this point of view can be found in a paper by Sergeev and Veselov \cite{SV10}.

\subsection{Notes}

The trick of lifting a family of (non-stable) symmetric polynomials to symmetric functions by introducing an additional formal parameter, which represents the dependence on the number of variables~$n$, has previously been used by Rains \cite{Rai05} and Sergeev and Veselov \cite{SV09}. Rains deals with the Koornwinder polynomials, whereas Sergeev and Veselov consider Jacobi symmetric polynomials.

In addition, this trick was recently used by Olshanski \cite{Ols11a} (see also \cite{Ols11b}) to introduce Laguerre and Meixner symmetric functions. His notion of Laguerre symmetric functions is a special case of ours, corresponding to $\alpha=1$. There is certain overlap between Olshanski's paper and the present one, but they are to a large extent complementary. Indeed, he addresses a number of problems that are not considered here. For example, orthogonality of the Laguerre symmetric functions, as well as a corresponding inf\/inite-dimensional dif\/fusion process. On the other hand, Olshanski does not discuss the Hermite case, the classif\/ication of invariant ideals in the ring of symmetric functions, and `super' polynomials.

The generalised Hermite and Laguerre polynomials in $n$-variables are, up to an overall (groundstate) factor, eigenfunctions of Schr\"odinger operators that def\/ine integrable quantum $n$-body systems. For the Hermite case the corresponding system is essentially the one originally considered by Calogero~\cite{Cal71}, and the system corresponding to the Laguerre polynomials appear in Section~11 of Olshanetsky and Perelomov~\cite{OP83}, as a generalisation of Calogero's system corresponding to the root system~$B_n$. A detailed discussion of this relationship between multivariable Hermite and Laguerre polynomials on the one hand and integrable quantum many-body systems on the other can be found in Baker and Forrester~\cite{BF97} and van Diejen~\cite{vD97}. These papers also contain a number of important results on such polynomials, as well as references to further related results in the literature.

\subsection{Notation and conventions}

We conclude this introduction with a few remarks on notation. In particular, on the parameter~$\alpha$: in the context of integrable systems it is typically replaced by its inverse~$1/\alpha$, denoted by a~number of dif\/ferent letters, e.g., $-k$ in \cite{SV09, SV04} or~$\theta$ in~\cite{KOO98,OO97,SV05}; and in literature related to Random Matrix Theory $\beta=2/\alpha$ is often used (as, e.g., in Baker and Forrester~\cite{BF97}). In an attempt to minimise confusion we shall throughout this paper only make use of the para\-me\-ter~$\alpha$. Regarding the natural numbers, we shall require both the set including and the set excluding the element zero. For that reason, we make use of the conventions $\mathbb{N}_0\equiv \lbrace 0,1,2,\ldots\rbrace$ and $\mathbb{N}\equiv\lbrace 1,2,\ldots\rbrace$.

\section{Symmetric functions}\label{symFuncsSection}
This section is largely a brief review of def\/initions and results from the theory of symmetric functions that we shall make use of. This review is intended to serve two purposes: f\/irstly, to f\/ix our notation, and secondly we hope that it will make the paper accessible to a somewhat wider audience. Throughout this section we shall in most cases adhere to the notation in Macdonald's book~\cite{Mac95}, to which the reader is referred for further details.

\subsection{Partitions}

A partition $\lambda = (\lambda_1,\lambda_2,\ldots,\lambda_i,\ldots)$ is a sequence of non-negative integers $\lambda_i$ such that
\begin{gather*}
    \lambda_1\geq\lambda_2\geq\cdots\geq\lambda_i\geq\cdots
\end{gather*}
and only a f\/inite number of the terms $\lambda_i$ are non-zero. The number of non-zero terms is referred to as the length of $\lambda$, and is denoted $\ell(\lambda)$. We shall not distinguish between two partitions that dif\/fer only by a string of zeros. The weight of a partition $\lambda$ is the sum
\begin{gather*}
    |\lambda|:= \lambda_1+\lambda_2+\cdots
\end{gather*}
of its parts, and its diagram is the set of points $(i,j)\in\mathbb{N}^2$ such that $1\leq j\leq\lambda_i$. Ref\/lection in the diagonal produces the conjugate partition
$\lambda^\prime=(\lambda_1',\lambda_2',\ldots)$. We use the notation $e_i$, $i\in\mathbb{N}$, for the sequence def\/ined by $(e_i)_j = \delta_{ij}$, where $\delta_{ij}$ is the Kronecker delta. In addition, we shall make use of the notation
\begin{gather*}
    \lambda^{(i)} = \lambda + e_i,\qquad \lambda_{(i)} = \lambda - e_i.
\end{gather*}

The set of all partitions of a given weight are partially ordered by the dominance order: $\lambda\leq \mu $ if and only if $\sum\limits_{i=1}^k\lambda_i\leq \sum\limits_{i=1}^k \mu_i$ for all $k\in\mathbb{N}$. One easily verif\/ies that $\lambda\leq\mu$ if and only if $\mu'\leq\lambda'$. We shall also require the inclusion order on the set of all partitions, def\/ined by $\mu\subseteq\lambda$ if and only if $\lambda_i\leq\mu_i$ for all $i$, or equivalently, if and only if the diagram of $\mu$ is contained in that of~$\lambda$.

To a partition $\lambda$ is associated the following product of $\alpha$-deformed hook lengths:
\begin{gather} \label{defhook}
    h_\lambda=\prod_{(i,j)\in\lambda}\left(1+a_\lambda(i,j)+\frac{1}{\alpha}l_\lambda(i,j)\right),
\end{gather}
involving the arm-lengths and leg-lengths
\begin{gather}\label{lengths}
    a_\lambda(i,j) = \lambda_i-j,\qquad l_\lambda(i,j) = \lambda^\prime_j-i.
\end{gather}
Closely related is the following $\alpha$-deformation of the Pochhammer symbol:
\begin{gather}\label{defpochhammer}
    [x]_\lambda = \prod_{1\leq i\leq \ell(\lambda)}\left(x-\frac{i-1}{\alpha}\right)_{\lambda_i} = \prod_{(i,j)\in\lambda}\left(x+a^\prime_\lambda(i,j)-\frac{1}{\alpha}l^\prime_\lambda(i,j)\right)
\end{gather}
with $(x)_n\equiv x(x+1)\cdots(x+n-1)$ the ordinary Pochhammer symbol, to which $\lbrack x\rbrack_\lambda$ clearly reduces for $\ell(\lambda)=1$, and where the second expression for $\lbrack x\rbrack_\lambda$ involves the co-arm-lengths and co-leg-lengths
\begin{gather}\label{colengths}
    a^\prime_\lambda(i,j) = j-1,\qquad l^\prime_\lambda(i,j) = i-1.
\end{gather}

\subsection{Symmetric functions}
The ring of symmetric polynomials in $n$ indeterminants $x=(x_1,\ldots,x_n)$ with integer coef\/f\/icients,
\begin{gather*}
	\Lambda_n = \mathbb{Z}\lbrack x_1,\ldots,x_n\rbrack^{S_n},
\end{gather*}
has a natural grading given by the degree of the polynomials:
\begin{gather*}
    \Lambda_n = \bigoplus_{k\geq 0}\Lambda_n^k,
\end{gather*}
where $\Lambda_n^k$ is the submodule consisting of all homogeneous symmetric polynomials of degree~$k$.

For a given $k$ and $n\in\mathbb{N}$, consider the homomorphism $\rho^k_{n,n-1}:\Lambda_n^k\to\Lambda_{n-1}^k$ def\/ined by
\begin{gather*}
	(\rho^k_{n,n-1}f)(x_1,\ldots,x_{n-1}) = f(x_1,\ldots,x_{n-1},0).
\end{gather*}
Let $\Lambda^k$ denote the module consisting of all sequences $(f_1,f_2,\ldots,f_n,\ldots)$ such that $f_n\in\Lambda_n^k$ and $\rho^k_{n,n-1}f_n=f_{n-1}$, and with the module structure given by term wise operations. The ring of symmetric functions can then be def\/ined as the graded ring
\begin{gather*}
    \Lambda = \bigoplus_{k\geq 0}\Lambda^k.
\end{gather*}
We note the restriction homomorphisms $\rho_n^k:\Lambda^k\to\Lambda_n^k$, which sends $f\in\Lambda^k$ to $f_n$, and
$\rho_n\equiv\oplus_{k\geq 1}\rho^k_n:\Lambda\to\Lambda_n$.

Given a (commutative) ring $A$ we will use the notation $\Lambda_A$ for the tensor product $A\otimes_{\mathbb{Z}}\Lambda$, and similarly for $\Lambda_n$. In this paper we shall mainly be concerned with either the f\/ields
\begin{gather}\label{fields}
    \mathbb{F} = \mathbb{Q}(\alpha),\qquad \mathbb{K} = \mathbb{Q}(a,\alpha),
\end{gather}
or the corresponding extensions generated by the indeterminate $p_0$,
\begin{gather}\label{fieldExtensions}
    \mathbf{F} = \mathbb{F}(p_0),\qquad \mathbf{K} = \mathbb{K}(p_0).
\end{gather}

There are several important and useful generators of $\Lambda$. We shall make use of the elementary- and power sum symmetric functions, given by
\begin{gather*}
    e_r := \lim_\leftarrow e_r(x),\qquad e_r(x) = \sum_{1\leq i_1<\cdots<i_r\leq n}x_{i_1}\cdots x_{i_r},\\
    p_r := \lim_\leftarrow p_r(x),\qquad p_r(x) = \sum_{i=1}^n x_i^r,
\end{gather*}
respectively. Here, $\lim\limits_\leftarrow g_r$, where $g_r(x_1,\dots,x_n)\in\Lambda^r_n$ for all $n\in \mathbb N$,
denotes the sequence $(g_r(x_1),g_r(x_1,x_2),g_r(x_1,x_2,x_3),\ldots)\in \Lambda^r$, and $r$ is allowed to be any non-negative integer. To be precise, the $e_r$ generate $\Lambda$, and are algebraically independent over $\mathbb{Z}$; while, on the other hand, the $p_r$ generate $\Lambda_{\mathbb{Q}}$, but not $\Lambda$, and are algebraically independent over~$\mathbb{Q}$. We recall the standard notation
\begin{gather*}
    e_\lambda = e_{\lambda_1}e_{\lambda_2}\cdots,\qquad p_\lambda = p_{\lambda_1}p_{\lambda_2}\cdots,
\end{gather*}
where $\lambda$ is any partition. In addition, we shall make use of the monomial symmetric functions
\begin{gather*}
    m_\lambda := \lim_\leftarrow m_\lambda(x),\qquad m_\lambda(x) = \sum_P x_1^{\lambda_{P(1)}}\cdots x_n^{\lambda_{P(n)}},
\end{gather*}
where the sum extends over all distinct permutations $P$ of $\lambda$. As $\lambda$ runs through all partitions, the $e_\lambda$ and $m_\lambda$ form a basis for $\Lambda$ and the $p_\lambda$ for $\Lambda_{\mathbb{Q}}$.

\subsection{Jack's symmetric functions}\label{JackFuncsSection}
Jack's symmetric functions form a further important, albeit more intricate, basis for $\Lambda_{\mathbb{F}}$, which will be a key ingredient in many constructions in this paper. In order to recall their def\/inition we start from the CMS operator
\begin{gather}\label{stableOp}
    D_n = \sum_{i=1}^nx_i^2\frac{\partial^2}{\partial x_i^2} + \frac{2}{\alpha}\sum_{i\neq j}\frac{x_ix_j}{x_i-x_j}\frac{\partial}{\partial x_i}.
\end{gather}
It is important to note that this operator preserves $\Lambda_{\mathbb{F},n}$. This follows from invariance under permutations of $x$ and the observation that, for $p\in\Lambda_{\mathbb{F},n}$, the polynomial $(\partial/\partial x_i-\partial/\partial x_j)p$ is antisymmetric under the interchange of $x_i$ and $x_j$, and hence divisible by $x_i-x_j$. Moreover, it is stable with respect to $\rho_{n,n-1}$:
\begin{gather*}
    \rho_{n,n-1}\circ D_n = D_{n-1}\circ\rho_{n,n-1}.
\end{gather*}
Hence, there is a unique operator $D$ on $\Lambda_{\mathbb{F}}$ such that $\rho_n\circ D=D_n\circ\rho_n$, given by the sequence $(D_1,D_2,\ldots,D_n,\ldots)$. Now, for a partition $\lambda$, the (monic) Jack's symmetric function $P_\lambda$ is the unique eigenfunction of this operator of the form
\begin{gather}\label{triangular}
    P_\lambda = m_\lambda + \sum_{\mu<\lambda}c_{\lambda\mu}m_\mu,\quad c_{\lambda\mu}\in\mathbb{F}.
\end{gather}
Note that $P_\lambda$ is homogeneous of degree $|\lambda|$. For further details see, e.g., Chapter~VI of Macdo\-nald~\cite{Mac95}.

It is a remarkable fact that to each so-called shifted symmetric function corresponds an eigenoperator of Jack's symmetric functions. We recall that the algebra of shifted symmetric polynomials $\Lambda_{\alpha,n}\equiv \Lambda_{\mathbb{F},\alpha,n}$ consists of all polynomials $p(x_1,\ldots,x_n)$ (over the f\/ield $\mathbb{F}\equiv\mathbb{Q}(\alpha)$) that are symmetric in the shifted variables $x_i-i/\alpha$; see, e.g., Okounkov and Olshanski~\cite{OO97}. It has a natural f\/iltration, given by the degree of the polynomials:
\begin{gather*}
    \Lambda_{\alpha,n}^0\subset\Lambda_{\alpha,n}^1\subset\cdots\subset\Lambda_{\alpha,n}^k\subset\cdots,
\end{gather*}
where $\Lambda_{\alpha,n}^k$ is the subspace of polynomials of degree at most $k$. As in the construction of the ring of symmetric functions, we can introduce linear spaces $\Lambda^k_\alpha$, and the algebra of shifted symmetric functions is then given as the f\/iltered algebra
\begin{gather*}
    \Lambda_\alpha=\bigcup_{k\geq 0}\Lambda_\alpha^k.
\end{gather*}
In particular, this algebra is freely generated by the shifted power sums
\begin{gather*}
    \pi_{r,\alpha}:= \lim_\leftarrow \pi_{r,\alpha}(x),\qquad \pi_{r,\alpha}(x) = \sum_{i=1}^n\left(\left(x_i-\frac{i}{\alpha}\right)^r-\left(-\frac{i}{\alpha}\right)^r\right).
\end{gather*}
Now, for any $f\in\Lambda_\alpha$, there is a unique operator $\mathcal{L}_f$ on $\Lambda_{\mathbb{F}}$ such that
\begin{gather}\label{JackEigenOps}
	\mathcal{L}_fP_\lambda = f(\lambda)P_\lambda
\end{gather}
for all partitions $\lambda$. In particular, the operator $D$ can be obtained as a linear combination of~$\mathcal{L}_{\pi_{1,\alpha}}$ and~$\mathcal{L}_{\pi_{2,\alpha}}$. A construction of these operators using Cherednik--Dunkl operators~\cite{BGHP93,Che94,Op00} can be found in Section~4 of Sergeev and Veselov~\cite{SV05}.

We recall the natural analogue of the specialisation $(x_1,\ldots,x_n)=(1,\ldots,1)$: for any $X\in\mathbb{F}$ def\/ine a homomorphism $\epsilon_X:\Lambda_{\mathbb{F}}\to\mathbb{F}$ by setting
\begin{gather*}
    \epsilon_X(p_r) = X,\qquad r\in\mathbb{N}.
\end{gather*}
Stanley~\cite{Sta89} (see also Section~VI.10 in Macdonald~\cite{Mac95}) has shown that the corresponding specialisation of Jack's symmetric functions is given by
\begin{gather}\label{epsilonXJackEq}
    \epsilon_X\big(P_\lambda\big) = \prod_{(i,j)\in\lambda}\frac{X+\alpha a_\lambda^\prime(i,j)-l_\lambda^\prime(i,j)}{\alpha a_\lambda(i,j)+l_\lambda(i,j)+1},
\end{gather}
c.f.~\eqref{lengths} and~\eqref{colengths}.

\subsection{CMS operators on the symmetric functions}
The $n$-variable Hermite and Laguerre polynomials introduced by Lassalle and Macdonald can be def\/ined as eigenfunctions of CMS operators of the form
\begin{gather}\label{defopLn}
	\mathcal{L}_n = \sum_{k=0}^2 a_kD_n^k + \sum_{\ell=0}^1 b_\ell E_n^\ell
\end{gather}
for some choice of coef\/f\/icients $a_k$ and $b_\ell$, and where the `building blocks'
\begin{gather}\label{EellDef}
	E_n^\ell = \sum_{i=1}^n x_i^\ell\frac{\partial}{\partial x_i}
\end{gather}
and
\begin{gather}\label{DkDef}
	D_n^k = \sum_{i=1}^n x_i^k\frac{\partial^2}{\partial x_i^2} + \frac{2}{\alpha}\sum_{i\neq j}\frac{x_i^k}{x_i-x_j}\frac{\partial}{\partial x_i}.
\end{gather}
More specif\/ically, the Hermite and Laguerre cases correspond to the following choice of coef\/f\/icients:
\begin{gather*}
	(a_2,a_1,a_0) = (0,0,1),\qquad (b_1,b_0) = (-2,0),\qquad \text{Hermite}\\
	(a_2,a_1,a_0) = (0,1,0),\qquad (b_1,b_0) = (-1,a+1),\qquad \text{Laguerre}.
\end{gather*}
A generic operator $\mathcal{L}_n$ is not stable under restrictions of the number of variables, i.e., $\rho_{n,n-1}\circ\mathcal{L}_n\neq \mathcal{L}_{n-1}\circ\rho_{n,n-1}$. Consequently, it does not directly lift to an operator on $\Lambda_{\mathbb{F}}$. In fact, only among the eigenoperators for the Jack polynomials can such stable CMS operators be found.

Nevertheless, by introducing a new indeterminate $p_0$, which ef\/fectively encodes the dependence on the number of variables $n$, to each CMS operator $\mathcal{L}_n$ we can in a natural manner assign an operator $\mathcal{L}$ on $\Lambda_{\mathbb{F}}$. In the Jacobi case ($(a_2,a_1,a_0)=(1,2,0)$ and $(b_1,b_0)=(-p-2q+1,-2p-2q+1)$) this fact was demonstrated by Sergeev and Veselov \cite{SV09}. Closely related is an earlier paper by Rains \cite{Rai05}, which concerns a symmetric function analogue of the Koornwinder polynomials.

As a f\/irst step towards making these remarks precise, we shall rewrite the CMS opera\-tors~\eqref{EellDef}, and~\eqref{DkDef} in a more convenient form. Fix $n\in\mathbb{N}$ and let $r=1,\ldots,n$. Then, we can def\/ine a~dif\/ferential operator $\partial^{(n)}(p_r)$ on $\Lambda_{\mathbb{F},n}$ by requiring that $\partial^{(n)}(p_r)1=0$ and
\begin{gather*}
    \partial^{(n)}(p_r)p_s =
    \begin{cases}
        1, & r=s,\\
        0, & r\neq s
    \end{cases}
\end{gather*}
for $s=1,\ldots,n$.

\begin{Lemma}\label{psLemma}
Set $p_0 = n$. Then, the differential operators $E^\ell_n$ and $D^k_n$ are given by
\begin{subequations}
\begin{gather}\label{Eps}
    E^\ell_n = \sum_{r=1}^n r p_{r+\ell-1}\partial^{(n)}(p_r)
\end{gather}
and
\begin{gather}
    D^k_n = \sum_{r,q=1}^n rq p_{r+q+k-2}\partial^{(n)}(p_r)\partial^{(n)}(p_q) + \sum_{r=2}^n r(r - 1)p_{r+k-2}\partial^{(n)}(p_r)\nonumber\\
     \hphantom{D^k_n =}{}
     + \frac{1}{\alpha}\sum_{r=1}^n r\sum_{m=0}^{r+k-2}(p_{r+k-2-m}p_m - p_{r+k-2})\partial^{(n)}(p_r),\label{Dps}
\end{gather}
\end{subequations}
respectively.
\end{Lemma}

\begin{proof}
We recall that $E^\ell_n$ are f\/irst-order dif\/ferential operators, and that $\Lambda_{\mathbb{F},n}$ is generated by the power sums~$p_r(x)$ with $r=1,\ldots,n$. Hence, it is suf\/f\/icient to compute their action on said power sums. This yields~\eqref{Eps}.

We thus turn to the dif\/ferential operators~$D^k_n$, and observe that their f\/irst-order terms act on the power sums as follows:
\begin{gather}
    2\sum_{i\neq j}\frac{x_i^k}{x_i-x_j}\frac{\partial}{\partial x_i}p_r(x)
    = \sum_{i\neq j}\frac{1}{x_i-x_j}\left(x_i^k\frac{\partial}{\partial x_i} - x_j^k\frac{\partial}{\partial x_j}\right)p_r(x)
      = \sum_{i\neq j}r \sum_{m=0}^{r+k-2}x_i^{r+k{-m}-2}x_j^m\nonumber\\
\hphantom{2\sum_{i\neq j}\frac{x_i^k}{x_i-x_j}\frac{\partial}{\partial x_i}p_r(x) }{}
= r\sum_{m=0}^{r+k-2}(p_{r+k-2-m}(x)p_m(x) - p_{r+k-2}(x)).\label{DAction1}
\end{gather}
For the second-order terms, it suf\/f\/icient to know the action on the power sums $p_r$ and the products of two power sums, i.e., on the terms $p_rp_q$ with $r,q = 1,\ldots,n$. If we allow $r,q=0$ and set $p_r(x)\equiv 0$ for $r<0$, then these cases are all included in the formula
\begin{gather}
    \sum_{i=1}^n x_i^k\frac{\partial^2}{\partial x_i^2}p_r(x)p_q(x)  = r(r - 1)p_{r+k-2}(x)p_q(x) + 2rqp_{r+q+k-2}(x)\nonumber\\
    \hphantom{\sum_{i=1}^n x_i^k\frac{\partial^2}{\partial x_i^2}p_r(x)p_q(x)  =}{}
    + q(q - 1)p_r(x)p_{q+k-2}(x).\label{DAction2}
\end{gather}
Combining these facts we obtain \eqref{Dps}.
\end{proof}

\begin{Remark}
These expressions for the dif\/ferential operators $E^\ell_n$ and $D^k_n$ involve power sums $p_r(x_1,\ldots,x_n)$ with $r>n$.
In principle, such terms can be rewritten in terms of power sums with $r\leq n$.
However, we have refrained from doing so since this would lead to rather complicated expressions. In addition, we are ultimately interested in operators on the algebra of symmetric functions, and there are no non-trivial relations between the power sum symmetric functions~$p_r$.
\end{Remark}

We now let $p_0$ be an indeterminate, and consider the f\/ield $\mathbf{F}\equiv\mathbb{F}(p_0)$.
It is clear that we can not specialise all $f\in\mathbf{F}$ to $p_0=n$.
Indeed, this is possible if and only if $f\in\mathbf{F}_{(p_0-n)}$: the (local) algebra of rational functions $g/h$ in $p_0$ over $\mathbb{F}$ such that $h(n)\neq 0$.
For simplicity of exposition, we shall make use of the short-hand notation $\mathbf{F}^{(n)}\equiv \mathbf{F}_{(p_0-n)}$.
We can now introduce, for each $n\in\mathbb{N}$, the specialisation map $\phi_n: \mathbf{F}^{(n)}\to\mathbb{F}$ by setting
\begin{gather*}
    \phi_n(f) = f\arrowvert_{p_0=n},
\end{gather*}
and thereby the  homomorphism $\varphi_n: \Lambda_{\mathbf{F}^{(n)}}\to\Lambda_{\mathbb{F},n}$ by{\samepage
\begin{gather*}
    \varphi_n(f\otimes p) = \phi_n(f)\otimes\rho_n(p).
\end{gather*}
We note that $\varphi_n$ is surjective for all $n\in\mathbb{N}$.}

On $\Lambda_{\mathbf{F}}$ we have obvious analogues $\partial(p_r)$ of the dif\/ferential operators $\partial^{(n)}(p_r)$; see Appendix \ref{AppDiffOps}. Moreover, in a natural sense, the former dif\/ferential operators are of degree $-r$; see the discussion preceding Lemma \ref{degreeLemma}. Lemma \ref{psLemma} thus suggests the following def\/inition of dif\/ferential opera\-tors~$E^\ell$ and~$D^k$ on $\Lambda_{\mathbf{F}}$:

\begin{Definition}\label{EDDef}
Let $\ell,k\in\mathbb{N}_0$. We then def\/ine dif\/ferential operators $E^\ell$ and $D^k$ on $\Lambda_{\mathbf{F}}$ by
\begin{gather*}
    E^\ell = \sum_{r=1}^\infty rp_{r+\ell-1}\partial(p_r)
\end{gather*}
and
\begin{gather*}
    D^k = \sum_{r,q=1}^\infty rq p_{r+q+k-2}\partial(p_r)\partial(p_q) + \sum_{r=2}^\infty r(r - 1)p_{r+k-2}\partial(p_r)\\
    \hphantom{D^k =}{}
    + \frac{1}{\alpha}\sum_{r=1}^\infty r\sum_{m=0}^{r+k-2}(p_{r+k-2-m}p_m - p_{r+k-2})\partial(p_r),
\end{gather*}
respectively.
\end{Definition}

That this is a natural def\/inition is conf\/irmed by the following lemma:

\begin{Lemma}\label{EDLemma}
Fix $k,\ell\in\mathbb{N}_0$. Then, $E^\ell$ and $D^k$ are homogeneous differential operators on $\Lambda_{\mathbf{F}}$ of degree $\ell-1$ and $k-2$, respectively. Moreover, they are the unique operators on $\Lambda_{\mathbf{F}}$ such that the diagrams
\begin{subequations}
\begin{gather}\label{EProj}
    \begin{CD}
        \Lambda_{\mathbf{F}^{(n)}} @>E^\ell>>\Lambda_{\mathbf{F}^{(n)}}\\
        @V\varphi_nVV @VV\varphi_nV\\
        \Lambda_{\mathbb{F},n} @>E^\ell_n>> \Lambda_{\mathbb{F},n}
    \end{CD}
\end{gather}
and
\begin{gather}\label{DProj}
    \begin{CD}
        \Lambda_{\mathbf{F}^{(n)}} @>D^k>>\Lambda_{\mathbf{F}^{(n)}}\\
        @V\varphi_nVV @VV\varphi_nV\\
        \Lambda_{\mathbb{F},n} @>D^k_n>> \Lambda_{\mathbb{F},n}
    \end{CD}
\end{gather}
\end{subequations}
are commutative for all $n\in\mathbb{N}$.
\end{Lemma}

\begin{proof}
The fact that both $E^\ell$ and $D^k$ are dif\/ferential operators on $\Lambda_{\mathbf{F}}$ is a direct consequence of Proposition \ref{diffOpsProp}. The stated homogeneity and degrees of $E^\ell$ and $D^k$ follows immediately from the observation that $\partial(p_r)$ and $p_r$ are homogeneous of degree $-r$ and $r$, respectively.

It follows from Def\/inition \ref{EDDef} and \eqref{EellDef} that
\begin{gather*}
    \varphi_n(E^\ell p_r) = rp_{r+\ell-1}(x_1,\ldots,x_n)=E^\ell_n(\varphi_np_r),\qquad r\in\mathbb{N},
\end{gather*}
where $p_0(x_1,\ldots,x_n)\equiv n$. We note that $E^\ell p_0=0$. Since $E^\ell$ and $E^\ell_n$ are f\/irst-order dif\/ferential operators, and $\varphi_n$ a $\mathbb{F}$-algebra homomorphism, this implies \eqref{EProj}. We observe that~\eqref{DAction1} and~\eqref{DAction2} hold true for any $r,q\in\mathbb{N}$. Comparing these formulae with Def\/inition \ref{EDDef} we f\/ind that $\varphi_n(D^kp_rp_q)=D^k_n\varphi_n(p_rp_q)$, $r,q\in\mathbb{N}$. Commutativity of the diagram \eqref{DProj} thus follows from the fact that $D^k$ and $D^k_n$ are dif\/ferential operators of order two, and that $D^kp_0=0$.

There remains only to prove uniqueness. Suppose that $D,D^\prime\in\mathscr{D}(\Lambda_{\mathbf{F}})$ are such that $\varphi_n\circ(D - D^\prime) = 0$ for all $n\in\mathbb{N}$. For any non-zero $p\in\Lambda_{\mathbf{F}}$ there exists $n\in\mathbb{N}$ such that $p\in\Lambda_{\mathbf{F}^{(n)}}$ and $\varphi_n(p)\neq 0$. Hence, $D = D^\prime$ and the statement follows.
\end{proof}

From Lemma \ref{EDLemma} we can immediately infer the following:

\begin{Proposition}\label{CMSOpsProp}
Let
\begin{gather}\label{GenCMSOps}
    \mathcal{L} = \sum_{k=0}^\infty a_k D^k + \sum_{\ell=0}^\infty { b_\ell} E^\ell
\end{gather}
for some coefficients $a_k,{b_\ell}\in\mathbf{F}$ such that only finitely many of them are non-zero. Moreover, let~$\mathcal{L}_n$ stand for the operator defined in \eqref{defopLn}.  Then, $\mathcal{L}$ is a differential operator on $\Lambda_{\mathbf{F}}$. Moreover, it is the unique operator on $\Lambda_{\mathbf{F}}$ such that the diagram
\begin{gather*}
    \begin{CD}
        \Lambda_{\mathbf{F}^{(n)}} @>\mathcal{L}>>\Lambda_{\mathbf{F}^{(n)}}\\
        @V\varphi_nVV @VV\varphi_nV\\
        \Lambda_{\mathbb{F},n} @>\mathcal{L}_n>> \Lambda_{\mathbb{F},n}
    \end{CD}
\end{gather*}
is commutative for all $n\in\mathbb{N}$.
\end{Proposition}

\section{Generalised hypergeometric series}\label{HypergSection}
In this section we def\/ine and study a natural analogue of hypergeometric series in the context of symmetric functions,
given as formal series in Jack's symmetric functions. When restricted to a f\/inite number of variables, these formal series coincide with (generalised) hypergeometric series studied, in particular, by Kor\'anyi \cite{Kor91}, Yan \cite{Yan92}, Kaneko \cite{Kan93},  and Macdonald \cite{Mac}.

We shall f\/irst introduce an analogue of Macdonald's hypergeometric series in two sets of variables. For that, we require the graded algebra
\begin{gather*}
    \Lambda_{\mathbf{F}}\otimes\Lambda_{\mathbf{F}} = \bigoplus_{k\geq 0}(\Lambda_{\mathbf{F}}\otimes\Lambda_{\mathbf{F}})^k,
\end{gather*}
where
\begin{gather*}
    (\Lambda_{\mathbf{F}}\otimes\Lambda_{\mathbf{F}})^k\equiv \big\lbrace p_1\otimes p_2: p_i\in\Lambda_{\mathbb{F}}^{k_i}~\text{with}~k_1+k_2=k\big\rbrace;
\end{gather*}
c.f., \eqref{fieldExtensions}. We consider the ideal
\begin{gather*}
    U = \bigoplus_{k\geq 1}(\Lambda_{\mathbf{F}}\otimes\Lambda_{\mathbf{F}})^k\subset \Lambda_{\mathbf{F}}\otimes\Lambda_{\mathbf{F}},
\end{gather*}
and equip $\Lambda_{\mathbf{F}}\otimes\Lambda_{\mathbf{F}}$ with the structure of a topological ring by requiring that the sequence of ideals~$U^n$, $n\in\mathbb{N}_0$, form a base of neighbourhoods of $0\in\Lambda_{\mathbf{F}}\otimes\Lambda_{\mathbf{F}}$. The corresponding completion, hereafter denoted by $\Lambda_{\mathbf{F}}\hat{\otimes}\Lambda_{\mathbf{F}}$, can be identif\/ied with the algebra of formal power series
\begin{gather*}
    \hat{p}=\sum_{\lambda,\mu}a_{\lambda\mu}\ p_\mu\otimes p_\lambda,\qquad a_{\mu\lambda}\in\mathbf{F}.
\end{gather*}

We are now ready to give the  precise def\/inition of the hypergeometric series in question.

\begin{Definition}\label{FpqDef}
Fix $p,q\in\mathbb{N}_0$ and let $(a_1,\ldots,a_p)\in\mathbf{F}^p$ and $(b_1,\ldots,b_q)\in\mathbf{F}^q$ be such that $(i-1)/\alpha-b_j\notin\mathbb{N}_0$ for all $i\in\mathbb{N}_0$.
We then def\/ine ${}_p\mathscr{F}_q(a_1,\ldots,a_p;b_1,\ldots,b_q;\alpha,p_0)\in\Lambda_{\mathbf{F}}\hat{\otimes}\Lambda_{\mathbf{F}}$ by
\begin{gather*}
    {}_p\mathscr{F}_q(a_1,\ldots,a_p;b_1,\ldots,b_q;\alpha,p_0) = \sum_\lambda \frac{1}{h_\lambda}\frac{\lbrack a_1\rbrack^{}_\lambda\cdots\lbrack a_p\rbrack^{}_\lambda}{\lbrack b_1\rbrack^{}_\lambda
    \cdots\lbrack b_q\rbrack^{}_\lambda}\frac{P_\lambda^{}\otimes P_\lambda^{}}{\epsilon_{p_0}(P_\lambda^{})},
\end{gather*}
where $h_\lambda$ and $\lbrack u \rbrack_\lambda$ are given by \eqref{defhook} and \eqref{defpochhammer}, respectively.
\end{Definition}

As in the f\/inite variable case, ${}_2\mathscr{F}_1$ satisf\/ies a simple dif\/ferential equation of second order. In order to make this remark precise, we f\/irst note that we can equip also $\Lambda_{\mathbf{F}}$ with the structure of a topological ring by starting from the ideal $U=\oplus_{k\geq 1}\Lambda_{\mathbf{F}}^k$. Then, any two continuous dif\/ferential operators $D_1$ and $D_2$ on $\Lambda_{\mathbf{F}}$ yield a continuous dif\/ferential operator $D_1\hat{\otimes}D_2$ on $\Lambda_{\mathbf{F}}\hat{\otimes}\Lambda_{\mathbf{F}}$ by
\begin{gather*}
    D_1\hat{\otimes}D_2\left(\sum_{\lambda_1,\lambda_2}a_{\lambda_1,\lambda_2}p_{\lambda_1}\otimes p_{\lambda_2}\right)\equiv \sum_{\lambda_1,\lambda_2}a_{\lambda_1,\lambda_2}(D_1p_{\lambda_1})\otimes (D_2p_{\lambda_2}).
\end{gather*}
It is important to note that dif\/ferential operators that we consider -- $E^\ell$ and $D^k$ for $\ell,k\in\mathbb{N}_0$ -- are all continuous. For a simple way to see this fact see Lemma \ref{degreeLemma} in Appendix \ref{AppDiffOps}. With this fact in mind, we proceed to state and prove the following:

\begin{Proposition}\label{prop2F1twosets}
Let $a,b,c\in\mathbf{F}$ be such that $(i-1)/\alpha-c\notin\mathbb{N}_0$ for all $i\in\mathbb{N}_0$. Then, ${}_2\mathscr{F}_1(a,b;c;\alpha,p_0)$ is the unique solution of the differential equation
\begin{gather}
    \big(D^1\hat{\otimes} 1\big) F + \left(c - \frac{p_0 - 1}{\alpha}\right)\big(E^0\hat{\otimes} 1\big) F - \big(1\hat{\otimes} D^3\big) F\nonumber\\
    \qquad{}- \left(a + b + 1 - \frac{2(p_0 - 1)}{\alpha}\right)\big(1\hat{\otimes} E^2\big) F = ab(1\hat{\otimes}p_1) F\label{2F1Eq}
\end{gather}
that is of the form
\begin{gather}\label{Fansatz}
    F = \sum_\lambda A_\lambda\frac{P_\lambda\otimes P_\lambda}{h_\lambda\epsilon_{p_0}(P_\lambda)},\qquad A_\lambda\in\mathbf{F},\qquad A_0=1.
\end{gather}
\end{Proposition}

\begin{proof} The proof follows closely that of Proposition~A.1 in Baker and Forrester~\cite{BF97}. Firstly, we observe that setting $k=2$ in \eqref{raisek} yields
\begin{gather*}
    D^3 = \frac{1}{2}\lbrack D^2,E^2\rbrack + \left(\frac{p_0-1}{\alpha}-1\right)E^2.
\end{gather*}
If we now take \eqref{Fansatz} as an ansatz for the solution $F$, then a straightforward, albeit somewhat lengthy, computation using Lemma~\ref{actionLemma} shows that the dif\/ferential equation~\eqref{2F1Eq} is satisf\/ied if and only if the coef\/f\/icients $A_\lambda$ solve the recurrence relation
\begin{gather*}
    \left(c+\lambda_i-\frac{(i-1)}{\alpha}\right)A_{\lambda^{(i)}}=\left(a+\lambda_i-\frac{(i-1)}{\alpha}\right)
    \left(b+\lambda_i-\frac{(i-1)}{\alpha}\right)A_{\lambda}.
\end{gather*}
Since we have f\/ixed $A_0=1$ and assumed that $(i-1)/\alpha-c\notin\mathbb{N}_0$, it is clear that this recurrence relation has a unique solution. Moreover, it follows immediately from the relation
\begin{gather*}
    \lbrack x\rbrack_{\lambda^{(i)}} = \lbrack x\rbrack_\lambda\left(x+\lambda_i-\frac{i-1}{\alpha}\right)
\end{gather*}
that this solution is given by
\begin{gather*}
    A_\lambda=\frac{[a]_\lambda  [b]_\lambda}{[c]_\lambda},
\end{gather*}
which clearly implies that the series $F$ is equal to ${}_2\mathscr{F}_1(a,b;c;\alpha,p_0)$.
\end{proof}

The hypergeometric series ${}_1\mathscr{F}_1$, ${}_0\mathscr{F}_1$ and ${}_0\mathscr{F}_0$ can be shown to satisfy analogous dif\/ferential equations. Since we shall make use of this fact in later parts of the paper, we proceed to deduce these dif\/ferential equations by exploiting suitable limit transitions from ${}_2\mathscr{F}_1$. To consider such limits, requires a topology of term-wise convergence of formal power series. For reasons that will become evident below, we shall work with symmetric functions over the real numbers, i.e., with $\Lambda_{\mathbb{R}}\hat{\otimes}\Lambda_{\mathbb{R}}$. Consequently, whenever they occur, we assume that $\alpha,p_0\in\mathbb{R}_+$. The restriction to positive numbers is made in order to avoid potential singularities of Jack's symmetric functions and ${}_p\mathscr{F}_q$. However, it is important to note that, since both Jack's symmetric functions as well as all coef\/f\/icients in ${}_p\mathscr{F}_q$ are rational functions of $\alpha$ and $p_0$, and the dif\/ferential operators that are involved are all of f\/inite degree, the dif\/ferential equations we deduce will hold true also in $\Lambda_{\mathbf{F}}\hat{\otimes}\Lambda_{\mathbf{F}}$.

In order to simplify the exposition somewhat, we shall write $\bar{\lambda}$ to indicate that $\bar{\lambda} = (\lambda^{(1)},\lambda^{(2)})$ for some partitions $\lambda^{(1)}$ and $\lambda^{(2)}$. It will also be convenient to use the corresponding short-hand notation $p_{\bar{\lambda}} = p_{\lambda^{(1)}}\otimes p_{\lambda^{(2)}}$. To each such `double-partition' $\bar{\lambda}$ we associate a function $C_{\bar{\lambda}}: \Lambda_{\mathbb{R}}\hat{\otimes}\Lambda_{\mathbb{R}}\rightarrow\mathbb{R}$ by the expansion
\begin{gather*}
    f = \sum_{\bar{\lambda}}C_{\bar{\lambda}}(f)p_{\bar{\lambda}},\qquad f\in\Lambda_{\mathbb{R}}\hat{\otimes}\Lambda_{\mathbb{R}}.
\end{gather*}
We note that any such function $C_{\bar{\lambda}}$ def\/ines a semi-norm $|\cdot|_{\bar{\lambda}}$ on $\Lambda_{\mathbb{R}}\hat{\otimes}\Lambda_{\mathbb{R}}$ by
\begin{gather*}
    |f|_{\bar{\lambda}} = \left|C_{\bar{\lambda}}(f)\right|,\qquad f\in\Lambda_{\mathbb{R}}\hat{\otimes}\Lambda_{\mathbb{R}},
\end{gather*}
where $|\cdot|$ in the right hand side denotes the standard (absolute value) norm on $\mathbb{R}$. The topology of term-wise convergence on $\Lambda_{\mathbb{R}}\hat{\otimes}\Lambda_{\mathbb{R}}$ is now the corresponding natural topology, def\/ined as the weakest topology in which all of these semi-norms, along with addition, are continuous. We note that, equipped with this topology, $\Lambda_{\mathbb{R}}\hat{\otimes}\Lambda_{\mathbb{R}}$ becomes a complete and metrisable locally convex vector space~-- a so-called Fr\'echet space. It is important to note that this topology of term-wise convergence does not depend on our specif\/ic choice of basis~-- in the discussion above $p_{\bar{\lambda}}$ with~$\bar{\lambda}$ running through all pairs of partitions $(\lambda^{(1)},\lambda^{(2)})$. These latter facts are all easy to infer from the general theory of locally convex vector spaces; see, e.g., Sections V.1-2 in Reed and Simon~\cite{RS80}.

We proceed to brief\/ly consider the relation to the $\hat{U}$-adic topology introduced at the beginning of this section. In particular, we observe that, for a sequence $\lbrace p_n\rbrace$ of elements $p_n\in\Lambda_{\mathbb{R}}\hat{\otimes}\Lambda_{\mathbb{R}}$, convergence in the $\hat{U}$-adic topology implies term-wise convergence. Moreover, we have the following lemma:

\begin{Lemma}\label{secondContinuityLemma}
If a differential operator $D$ on $\Lambda_{\mathbb{R}}\hat{\otimes}\Lambda_{\mathbb{R}}$ is continuous in the $\hat{U}$-adic topology, then it is continuous in the topology of term-wise convergence.
\end{Lemma}

\begin{proof}
Let $\lbrace q_n\rbrace$ be a sequence of elements $q_n\in\Lambda_{\mathbb{R}}\hat{\otimes}\Lambda_{\mathbb{R}}$ such that $q_n\rightarrow 0$ term-wise. Fix a~`double'-partition $\bar{\mu}$. By assumption, $D$ is continuous in the $\hat{U}$-adic topology. It follows that there exists $m\in\mathbb{N}_0$ such that
\begin{gather*}
    D\hat{U}_m\subset\hat{U}_{|\bar{\mu}|+1}.
\end{gather*}
We can thus deduce that
\begin{gather*}
    |Dq_n|_{\bar{\mu}}  = \left|\sum_{\bar{\lambda}}C_{q_n}(\bar{\lambda})Dp_{\bar{\lambda}}\right|_{\bar{\mu}} = \left|\sum_{|\bar{\lambda}|<m}C_{q_n}(\bar{\lambda})\sum_{\bar{\lambda}^\prime}
    C_{Dp_{\bar{\lambda}}}(\bar{\lambda}^\prime)p_{\bar{\lambda}^\prime}\right|_{\bar{\mu}}
     \leq \sum_{|\bar{\lambda}|<m}|C_{q_n}(\bar{\lambda})||C_{Dp_{\bar{\lambda}}}(\bar{\mu})|.
\end{gather*}
Hence, the fact that the latter sum is f\/inite implies that $|Dq_n|_{\bar{\mu}}\rightarrow 0$.
\end{proof}

We continue by considering limit transitions from the hypergeometric series ${}_2\mathscr{F}_1$. For $\gamma\in\mathbb{R}$, let $\sigma_\gamma:\Lambda_{\mathbb{R}}\rightarrow\Lambda_{\mathbb{R}}$ be the automorphism given by
\begin{gather}\label{defsigma}
    \sigma_\gamma(p_r)=\gamma^r p_r, \qquad r\in\mathbb{N}.
\end{gather}
Since $1\otimes\sigma_\gamma$ is degree preserving, it is continuous, and extends uniquely to a homomor\-phism~$1\hat{\otimes}\sigma_\gamma$ on $\Lambda_{\mathbb{R}}\hat{\otimes}\Lambda_{\mathbb{R}}$. In particular, we have that
\begin{gather*}
    (1\hat{\otimes}\sigma_{1/b}) {}_2\mathscr{F}_1(a,b;c;\alpha;p_0) = \sum_\lambda\frac{\lbrack a\rbrack_\lambda\lbrack b\rbrack_\lambda}{b^{|\lambda|}\lbrack c\rbrack_\lambda}\frac{P_\lambda\otimes P_\lambda}{\epsilon_{p_0}(P_\lambda)h_\lambda}.
\end{gather*}
In the sense of term-wise convergence, this implies the limit
\begin{gather*}
    \lim_{b\rightarrow\infty} (1\hat{\otimes}\sigma_{1/b}) {}_2\mathscr{F}_1(a,b;c;\alpha,p_0) = {}_1\mathscr{F}_1(a;c;\alpha,p_0).
\end{gather*}

{\sloppy Consider now the dif\/ferential equation \eqref{2F1Eq} for $F = {}_2\mathscr{F}_1$, and apply the homomor\-phism~$1\hat{\otimes}\sigma_{1/b}$. For any homogeneous dif\/ferential operator $D$ of f\/inite degree $\mathrm{deg}(D)$, we have that
\begin{gather}\label{sigmaD}
    \sigma_\gamma\circ D = \gamma^{\mathrm{deg}(D)}D\circ\sigma_\gamma.
\end{gather}
It follows from Lemma \ref{degreeLemma} and Lemma \ref{secondContinuityLemma} that such a dif\/ferential operator $D$ is continuous with respect to the topology of term-wise convergence, and thereby that it commutes with the limit in question. Using this fact, a direct computation yields the dif\/ferential equation satisf\/ied by ${}_1\mathscr{F}_1$. After computing similar limits in the parame\-ters~$a$ and~$c$, we arrive at the following proposition:

}

\begin{Proposition} \label{PropEq0F1}
Let $a,c\in\mathbf{F}$ be such that $(i-1)/\alpha-c\notin\mathbb{N}_0$ for all $i\in\mathbb{N}_0$. Then, ${}_1\mathscr{F}_1(a,c;\alpha,p_0)$ is a solution of
\begin{gather}
    (D^1\hat{\otimes} 1) F + \left(c - \frac{p_0 - 1}{\alpha}\right)(E^0\hat{\otimes} 1) F - (1\hat{\otimes} E^2) F = a(1\hat{\otimes}p_1) F,
\end{gather}
${}_0\mathscr{F}_1(c;\alpha,p_0)$ is a solution of
\begin{gather}\label{0F1DiffEq}
    (D^1\hat{\otimes} 1) F + \left(c - \frac{p_0 - 1}{\alpha}\right)(E^0\hat{\otimes} 1) F = (1\hat{\otimes}p_1) F,
\end{gather}
and ${}_0\mathscr{F}_0(\alpha,p_0)$ is a solution of
\begin{gather}\label{0F0DiffEq}
    (E^0\hat{\otimes} 1) F = (1\hat{\otimes}p_1) F.
\end{gather}
\end{Proposition}

We conclude this section by brief\/ly considering the hypergeometric series
\begin{gather}\label{eqpFqoneset}
    {}_pF_q(a_1,\ldots,a_p;b_1,\ldots,b_q;\alpha,p_0) = \sum_\lambda \frac{1}{h_\lambda}\frac{\lbrack a_1\rbrack_\lambda\cdots\lbrack a_p\rbrack_\lambda}{\lbrack b_1\rbrack_\lambda\cdots\lbrack b_q\rbrack^{}_\lambda} P_\lambda,
\end{gather}
which can be obtained by applying the homomorphism $1\otimes\epsilon_{p_0}$ to each term in ${}_p\mathscr{F}_q$. In this equation, it is assumed that the indeterminates $(b_1,\ldots,b_q)$ comply with the conditions stated in Def\/inition~\ref{FpqDef}. The next Proposition generalises a result of Yan~\cite{Yan92} and Kaneko~\cite{Kan93} on the solution of a multivariable generalisation of Euler's hypergeometric equation. The proof is omitted since it closely parallels that of  Proposition~\ref{prop2F1twosets}.

\begin{Proposition}\label{propeq2F1}
Let $a,b,c\in\mathbf{F}$ be such that $(i-1)/\alpha-c\notin\mathbb{N}_0$ for all $i\in\mathbb{N}_0$. Then, $_{2}F_{1}(a,b;c;\alpha,p_0)$ is the unique solution  of the differential equation
\begin{gather*}
    D^1F-D^2F +\left(c-\frac{p_0-1}{\alpha}\right)E^0F-\left(a+b+1-\frac{p_0-1}{\alpha}\right)E^1F=abp_0 F
\end{gather*}
that is of the form
\begin{gather*}
    F = \sum_\lambda \frac{A_\lambda}{h_\lambda}P_\lambda,\qquad A_\lambda\in\mathbf{F},\qquad A_0=1.
\end{gather*}
\end{Proposition}

\section{Hermite symmetric functions}\label{HermiteSec}
In this section we introduce and study Hermite symmetric functions as eigenfunctions of the dif\/ferential operator
\begin{gather}\label{HermiteOP}
    \mathcal{L}^H = D^0 - 2\nu^2 E^1
\end{gather}
with the parameter $\nu\in\mathbf{F}$. As in the f\/inite variable case, one can essentially remove the dependence on the parameter $\nu$. More precisely, since $D^0$ and $E^1$ are of degree $-2$ and $0$, respectively, we have that
\begin{gather}\label{sigmaInter}
    \sigma_{1/\nu}\circ(D^0 - 2\nu^2 E^1) = \nu^2(D^0 - 2E^1)\circ\sigma_{1/\nu};
\end{gather}
c.f., \eqref{sigmaD}. Using this fact, we can reduce most of the statements below to that for a f\/ixed value of~$\nu$. However, $\nu$ will play an important role in our discussion of a particular duality of the Hermite symmetric functions. For the moment we therefore refrain from specifying a f\/ixed value for~$\nu$.

It is readily inferred from Lemma \ref{actionLemma} that
\begin{gather*}
    \mathcal{L}^HP_\lambda = -2\nu^2|\lambda|P_\lambda + \sum_{\mu\subset\lambda}c_{\lambda_\mu}P_\mu
\end{gather*}
for some coef\/f\/icients $c_{\lambda\mu}\in\mathbf{F}$. By Theorem \ref{TheoExistenceEigenfunction}, it is thus clear that we can make the following def\/inition:

\begin{Definition}\label{HermiteDef}\sloppy
Let $\lambda$ be a partition. We then def\/ine the Hermite symmetric function $H_\lambda(\alpha,p_0,\nu^2)$ as the unique symmetric function such that
\begin{enumerate}\itemsep=0pt
\item $H_\lambda = P_\lambda + \sum\limits_{\mu\subset\lambda}u_{\lambda\mu}P_\mu$ for some $u_{\lambda\mu}\in\mathbf{F}$,
\item $\mathcal{L}^HH_\lambda=-2\nu^2|\lambda|H_\lambda$.
\end{enumerate}
\end{Definition}

\begin{Remark}
The generalised Hermite polynomials are recovered by setting $p_0=n$ and restric\-ting to $n$ indeterminates $x=(x_1,\ldots,x_n)$. Indeed, using Proposition~\ref{CMSOpsProp} it is readily verif\/ied that the resulting symmetric polynomials satisfy def\/initions given by Lassalle~\cite{Las91b} and Macdonald~\cite{Mac}.
\end{Remark}

Before proceeding to further investigate the properties of the Hermite symmetric functions, we detail a constructive def\/inition in terms of the Jack symmetric functions; c.f., (3.21) in Baker and Forrester~\cite{BF97} for the corresponding result in the f\/inite variable case. This requires the following notation: given a dif\/ferential operator $D$ on $\Lambda_{\mathbf{F}}$ and $L\in\mathbb{N}$, we let
\begin{gather}\label{expL}
    \exp_L(D)=1+\sum_{k=1}^L\frac{1}{k!}(D)^k.
\end{gather}
Clearly, $\exp_L(D)$ is a dif\/ferential operator on $\Lambda_{\mathbf{F}}$. Furthermore, if $D$ has f\/inite degree (see the paragraph preceding Lemma \ref{degreeLemma}), then so has $\exp_L(D)$, which, by Lemma~\ref{degreeLemma}, implies continuity. We stress the importance of truncating the series in the right-hand side of \eqref{expL} at some positive integer~$L$. Indeed, if this is not done, then we do not obtain a dif\/ferential operator on $\Lambda_{\mathbf{F}}$, c.f.~the paragraph containing~\eqref{infOrder}.

\begin{Proposition}\label{PropLassalle}
For any $L\geq \lfloor|\lambda|/2\rfloor$, we have that
\begin{gather}\label{HermiteExpRep}
    H_\lambda=\exp_L\left(-\frac{1}{4\nu^2}D^0\right)(P_\lambda).
\end{gather}
\end{Proposition}

\begin{proof}
For simplicity of exposition, we let $\Delta=-\frac{1}{4\nu^2}D^0$. It follows immediately from \eqref{D0Action} in Lemma \ref{actionLemma} that $\exp_L(\Delta)(P_\lambda)$ satisf\/ies property~(1) in Def\/inition~\ref{HermiteDef}. Since $E^1P_\lambda=|\lambda|P_\lambda$ and $[E^1,\Delta]=-2\Delta$, we have
\begin{gather*}
    E^1 \big(\Delta^k(P_\lambda)\big)=(|\lambda|-2k)\Delta^k(P_\lambda).
\end{gather*}
Consequently,
\begin{gather*}
    \mathcal{L}^H \big(\exp_L(\Delta)(P_\lambda)\big)  =-2\nu^2\big(E^1+2\Delta\big)\left( P_\lambda+\Delta (P_\lambda)+\frac{1}{2!}\Delta^2 (P_\lambda)+\cdots +\frac{1}{L!}\Delta^L(P_\lambda)\right)\\
\hphantom{\mathcal{L}^H \big(\exp_L(\Delta)(P_\lambda)\big)}{}
     = -2\nu^2\Bigg(|\lambda|P_\lambda+ (|\lambda|-2)\Delta (P_\lambda) +\frac{(|\lambda|-4)}{2!}\Delta^2 (P_\lambda)+\cdots\\
\hphantom{\mathcal{L}^H \big(\exp_L(\Delta)(P_\lambda)\big)=}{}
     +\frac{|\lambda|-2L}{L!}\Delta^L(P_\lambda)+2\Delta (P_\lambda)+2\Delta^2(P_\lambda)+\cdots+\frac{2}{(L-1)!}\Delta^L(P_\lambda)\Bigg)\\
\hphantom{\mathcal{L}^H \big(\exp_L(\Delta)(P_\lambda)\big)}{}
 =-2\nu^2|\lambda|\exp_L(\Delta)(P_\lambda),
\end{gather*}
i.e., also property (2) is satisf\/ied by $\exp_L(\Delta)(P_\lambda)$.
\end{proof}

\subsection{A duality relation}
We proceed to establish a particular duality relation for the Hermite symmetric functions that is not present at the level of the corresponding symmetric polynomials. To this end, we recall the standard automorphism $\omega_\gamma$, $\gamma\in\mathbb{F}$, of $\Lambda_{\mathbb{F}}$, given by
\begin{gather}\label{omegaDef}
    \omega_\gamma(p_r) = {(-1)^{r-1}}\gamma p_r,\qquad r\in\mathbb{N}.
\end{gather}
It is well known that, for a given value of the parameter $\alpha$, Jack's symmetric functions corresponding to the inverse parameter value $1/\alpha$ can be obtained by the following duality relation:
\begin{gather}\label{JackDuality}
    \omega_\alpha\big(P_\lambda(\alpha)\big) = Q_{\lambda^\prime}(1/\alpha),
\end{gather}
where
\begin{gather}\label{QDef}
    Q_\lambda = b_\lambda P_\lambda,\qquad b_\lambda = \prod_{(i,j)\in\lambda}\frac{l_\lambda(i,j)+1+\alpha  a_\lambda (i,j)}{l_\lambda(i,j)+\alpha+\alpha  a_\lambda(i,j)};
\end{gather}
see, e.g., Section VI.10 in Macdonald~\cite{Mac95}. This duality relation can be inferred from the identity
\begin{gather*}
    -\alpha\big(\omega_\alpha\circ D(\alpha)\big) = D(1/\alpha),\qquad D=D^2-\frac{2}{\alpha}(p_0-1)E^1,
\end{gather*}
and the fact that Jack's symmetric functions can be def\/ined as the unique eigenfunctions of $D$ that are of the form \eqref{triangular}; c.f.,~\eqref{stableOp}; and see Lemma~\ref{lemmadualityDE} below and note that~$D$ is independent of~$p_0$.

In order to deduce an analogous duality relation for the Hermite symmetric functions, we must consider also the parameter~$p_0$. The reason being that these symmetric functions have no eigenoperators that are independent of $p_0$. We therefore extend the automorphism~$\omega_\gamma$ to $\Lambda_{\mathbf{F}}$ by setting
\begin{gather*}
    \omega_\gamma(p_0) = -\gamma p_0,
\end{gather*}
or, equivalently, by replacing $\mathbb{N}$ by $\mathbb{N}_0$ in~\eqref{omegaDef}. With this extension in force, it is straightforward to determine the ef\/fect of $\omega_\alpha$ on the CMS operators $E^\ell$ and $D^k\equiv D^k(\alpha,p_0)$.

\begin{Lemma}\label{lemmadualityDE}
We have that
\begin{gather*}
    {\omega_\alpha \circ E^\ell=(-1)^{\ell-1}E^\ell\circ\omega_\alpha}
\end{gather*}
and
\begin{gather*}
    (-1)^{k-1}\alpha\big(\omega_\alpha\circ D^k(\alpha,p_0)\big) = D^{k}(1/\alpha,-\alpha p_0)\circ\omega_\alpha - (\alpha+1)k\big(E^{k-1}\circ\omega_\alpha\big).
\end{gather*}
\end{Lemma}

\begin{proof}
It follows immediately from \eqref{omegaDef} that
\begin{subequations}
\begin{gather}\label{promegaRels}
    \omega_\alpha\circ p_r  = (-1)^{r-1}\alpha (p_r\circ\omega_\alpha),\qquad r\in\mathbb{N}_0,\\
    \omega_\alpha\circ\partial(p_r)  = \frac{(-1)^{r-1}}{\alpha}\big(\partial(p_r)\circ\omega_\alpha\big),\qquad r\in\mathbb{N}.
\end{gather}
\end{subequations}
Using these relations, a direct computation yields the statement for $E^\ell$. We continue by obser\-ving that
\begin{gather*}
    D^k + \frac{k}{\alpha}E^{k-1} = \sum_{r,q=1}^\infty rqp_{r+q+k-2}\partial(p_r)\partial(p_q)\\
    \hphantom{D^k + \frac{k}{\alpha}E^{k-1} =}{}
    + \sum_{r=2}^\infty r(r-1)\left(1-\frac{1}{\alpha}\right)p_{r+k-2}\partial(p_r)+\frac{1}{\alpha}\sum_{r=1}^\infty r\sum_{m=0}^{r+k-2} p_{r+k-2-m}p_m\partial(p_r).
\end{gather*}
Using again \eqref{promegaRels}, it is readily seen that
\begin{gather*}
    \omega_\alpha\circ\left(D^k(\alpha,p_0) + \frac{k}{\alpha}E^{k-1}\right) = \frac{(-1)^{k-1}}{\alpha}\big(D^k(1/\alpha,-\alpha p_0) - k\alpha E^{k-1}\big)\circ\omega_\alpha,
\end{gather*}
which clearly implies the statement for $D^k$.
\end{proof}

There are now (at least) two dif\/ferent methods by which we can establish a duality relation for the Hermite symmetric functions. Firstly, we can follow the method sketched above for Jack's symmetric functions; and, secondly, we can make use of the representation \eqref{HermiteExpRep}. Here, we shall employ the latter method, since it yields a somewhat shorter proof.

\begin{Proposition}\label{HermiteDualityProp}
We have the duality relation
\begin{gather}\label{HermiteDuality}
    \omega_\alpha\big(H_\lambda(\alpha,p_0,\nu^2)\big) = b_{\lambda^\prime}(1/\alpha) H_{\lambda^{\prime}}\big(1/\alpha,-\alpha p_0,-\alpha \nu^2\big).
\end{gather}
\end{Proposition}

\begin{proof}
Starting from \eqref{HermiteExpRep}, we infer from Lemma \ref{lemmadualityDE} that
\begin{gather*}
    \omega_\alpha\big(H_\lambda(\alpha,p_0,\nu^2)\big) = \exp_L\left(\frac{1}{{ 4}\alpha\nu^2}D^0(1/\alpha,-\alpha p_0)\right)\big(Q_{\lambda^\prime}(1/\alpha)\big).
\end{gather*}
The statement is now a direct consequence of \eqref{QDef}.
\end{proof}

We stress that the duality relation \eqref{HermiteDualityProp} has no direct analogue in the f\/inite variable case. Indeed, the `restriction' homomorphism $\varphi_n$, which maps $H_\lambda$ to $H_\lambda(x_1,\ldots,x_n)$, f\/ixes $p_0=n$, whereas $\omega_\alpha$ maps $p_0$ to $-\alpha p_0$, and thus can not be restricted to $\Lambda_{\mathbb{F},n}=\varphi_n(\Lambda_{\mathbf{F}^{(n)}})$. However, this duality relation does have a natural analogue for the super Hermite polynomials, introduced in Section~\ref{sHermitePols}.

In the remainder of this section the parameter $\nu$ will not play any particular role. From hereon, we shall therefore assume that $\nu=1$. If needed, then this parameter can be reintroduced by applying the automorphism $\sigma_{\nu}$; c.f.~the paragraph containing~\eqref{sigmaInter}.

\subsection{A generating function}\label{SectionGenFunc}
We proceed to establish a generating function for the Hermite symmetric functions. As a f\/irst example of its usefulness, we shall then use this generating function to construct higher-order eigenoperators for the Hermite symmetric functions. These results will be obtained as rather direct generalisations of corresponding results due to Baker and Forrester \cite{BF97} on the generalised Hermite polynomials -- in turn based on an unpublished manuscript by Lassalle.

\begin{Proposition}\label{HermiteGenFuncProp}
We have that
\begin{gather}\label{HermiteGenFunc}
    \sum_\lambda \frac{1}{h_\lambda\epsilon_{p_0}(P_\lambda)} H_\lambda\otimes P_\lambda = {}_0\mathscr{F}_0 e^{-\frac{1}{4}(1\otimes p_2)}
\end{gather}
with
\begin{gather*}
    e^{-\frac{1}{4}(1\otimes p_2)}:= \sum_{n=0}^\infty \frac{1\otimes (-p_2/4)^n}{n!}.
\end{gather*}
\end{Proposition}

\begin{proof}
Let $x = (x_1,\ldots,x_n)$ and $y = (y_1,\ldots,y_n)$ be two sequences of indeterminates. By wri\-ting $\varphi_{n,x}$ and $\varphi_{n,y}$ we indicate that the homomorphism map $\Lambda_{\mathbf{F}^{(m)}}$ onto the algebra of symmetric polynomials in the indeterminates $x$ and $y$, respectively. We have that
\begin{gather}
    (\varphi_{n,x}\hat{\otimes}\varphi_{n,y}){}_0\mathscr{F}^{}_0 e^{-\frac{1}{4}(1\otimes p_2)} = {}_0\mathscr{F}^{}_0(x,y) e^{-\frac{1}{4}(1\otimes p_2(y))},
\end{gather}
where
\begin{gather*}
    {}_0\mathscr{F}_0(x;y) = \sum_\lambda \frac{1}{h_\lambda P_\lambda(1^n)} P_\lambda(x)\otimes P_\lambda(y)
\end{gather*}
and $e^{-\frac{1}{4}(1\otimes p_2(y))}$ is def\/ined in the obvious way. As shown by Baker and Forrester \cite{BF97} (see their Proposition 3.1),
\begin{gather}\label{expansion1}
    {}_0\mathscr{F}_0(x;y) e^{-\frac{1}{4}(1\otimes p_2(y))} = \sum_\lambda \frac{1}{h_\lambda P_\lambda(1^n)}H_\lambda(x)\otimes P_\lambda(y),
\end{gather}
where the sum is over all partitions $\lambda$ such that $\ell(\lambda)\leq n$. We note that Baker and Forrester use the normalisation $C_\lambda(x) = |\lambda|! P_\lambda(x)/h_\lambda$ for the Jack polynomials, and that their generalised Hermite polynomials are equal to $2^{|\lambda|}H_\lambda(x)/P_\lambda(1^n)$. As a consequence, the generating function expansion \eqref{expansion1} dif\/fers slightly from that stated by Baker and Forrester.

By a direct expansion of the right hand side of \eqref{HermiteGenFunc} in terms of Jack symmetric functions we obtain
\begin{gather}\label{expansion2}
    {}_0\mathscr{F}_0 e^{-\frac{1}{4}(1\otimes p_2)} = \sum_\lambda \frac{1}{h_\lambda\epsilon_{p_0}(P_\lambda)}U_\lambda\otimes P_\lambda
\end{gather}
for some $U_\lambda\in\Lambda_{\mathbf{F}}$. Using \eqref{Etop} for $k=2$ and $\ell=0$, \eqref{E2Action} and Proposition \ref{PieriJackProp}, it is readily verif\/ied that these symmetric functions are of the form
\begin{gather}\label{UExpansion}
    U_\lambda = \sum_{\mu\subset\lambda} v_{\lambda\mu}P_\lambda,\qquad v_{\lambda\mu}\in\mathbf{F}.
\end{gather}
If we compare the two expansions \eqref{expansion1} and \eqref{expansion2}, then we f\/ind that
\begin{gather*}
    \varphi_{n,x}(U_\lambda) = H_\lambda(x)\equiv \varphi_{n,x}(H_\lambda),\qquad \forall\, n\geq\ell(\lambda).
\end{gather*}
Since both $U_\lambda$ and $H_\lambda$ depend rationally on $p_0$, i.e., when expanded in, e.g., Jack's symmetric functions, the coef\/f\/icients are rational functions of $p_0$, it follows that $U_\lambda = H_\lambda$.
\end{proof}

Proposition \ref{HermiteGenFuncProp} can also be established from f\/irst principles by essentially the same method used by Baker and Forrester~\cite{BF97} to prove~\eqref{expansion1}.

The generating function~\eqref{HermiteGenFunc} is an ef\/fective tool for establishing a number of basic properties of the Hermite symmetric functions. For example, the ef\/fect of multiplication by $p_1$ and application of the dif\/ferential operator $E_0$. However, since these results can be obtained in complete analogy with the proofs of Corollaries~3.4 and~3.5 in Baker and Forrester~\cite{BF97}, we leave it to the interested reader to work out the details. Furthermore, the recurrence relation corresponding to multiplication by $p_1$ is the simplest special case of the complete set of recurrence relations we shall obtain in Section~\ref{section4.4}; see Proposition \ref{HermitePieriFormulaeProp}.

We proceed to use Proposition~\ref{HermiteGenFuncProp} to obtain higher-order eigenoperators for the Hermite symmetric functions. Also this result can be deduced in close analogy with Proposition~3.2 in Baker and Forrester~\cite{BF97}. At this point it might be helpful to recall the discussion of eigenoperators for Jack's symmetric functions in Section~\ref{JackFuncsSection}.

Let $D$ be a dif\/ferential operator of order $k$ on $\Lambda_{\mathbf{F}}\hat{\otimes}\Lambda_\mathbf{F}$ (see Appendix \ref{AppDiffOps} for the def\/inition of order). It follows from the Baker--Campbell--Hausdorf\/f formula that, for any $f\in\Lambda_{\mathbf{F}}\hat{\otimes}\Lambda_{\mathbf{F}}$,
\begin{gather}
    D fe^{-\frac{1}{4}(1\otimes p_2)} = e^{-\frac{1}{4}(1\otimes p_2)}\Bigg(D + \frac{1}{4}\big\lbrack (1\hat{\otimes}p_2),D \big\rbrack + \frac{1}{4^2 2!}\big\lbrack(1\hat{\otimes}p_2),\big\lbrack(1\hat{\otimes}p_2),D\big\rbrack\big\rbrack\nonumber\\
\hphantom{D fe^{-\frac{1}{4}(1\otimes p_2)} = e^{-\frac{1}{4}(1\otimes p_2)}\Bigg(}{}
    +\cdots+\frac{1}{4^k k!}\big\lbrack(1\hat{\otimes}p_2),\ldots,\big\lbrack(1\hat{\otimes}p_2),D\big\rbrack\cdots\big\rbrack\Bigg) f,\label{DBCH}
\end{gather}
where $1\hat{\otimes}p_2$ denotes the operator of multiplication by $1\otimes p_2$. From Def\/inition \ref{FpqDef} (for $p=q=0$) we can directly infer that
\begin{gather*}
    (\mathcal{L}_f\hat{\otimes} 1) {}_0\mathscr{F}_0 = (1\hat{\otimes}\mathcal{L}_f) {}_0\mathscr{F}_0,\quad \forall f\in\Lambda_{\mathbb{F},\alpha},
\end{gather*}
where $\mathcal{L}_f$ denotes the eigenoperator for Jack's symmetric functions given by~\eqref{JackEigenOps}. Moreover, since $p_2=\lbrack E^2,p_1\rbrack$, we can infer from Proposition~\ref{PieriJackProp} and Lemma~\ref{actionLemma} that
\begin{gather*}
	(D^0\hat{\otimes} 1){}_0\mathscr{F}_0 = (1\hat{\otimes}p_2){}_0\mathscr{F}_0.
\end{gather*}
Using these facts, as well as \eqref{DBCH} for $D=1\hat{\otimes}\mathcal{L}_f$, it is a matter of straightforward computations to verify that the following Proposition holds true:

\begin{Proposition}\label{PropLfH}
Let $f\in\Lambda_{\mathbb{F},\alpha}$, and let $k$ be the degree of $f$. Then, we have that
\begin{gather*}
    \big(1\hat{\otimes}\mathcal{L}_f\big){}_0\mathscr{F}_0 e^{-\frac{1}{4}(1\otimes p_2)} = \big(\mathcal{L}^H_f\hat{\otimes}1\big){}_0\mathscr{F}_0 e^{-\frac{1}{4}(1\otimes p_2)},
\end{gather*}
where
\begin{gather}\label{HermiteEigenOps}
    \mathcal{L}^H_f = \mathcal{L}_f - \frac{1}{4}\big\lbrack D^0,\mathcal{L}_f \big\rbrack
    + \frac{1}{4^2 2!}\big\lbrack D^0,\big\lbrack D^0,\mathcal{L}_f\big\rbrack\big\rbrack+\cdots+\frac{(-1)^k}{4^k k!}\big\lbrack D^0,\ldots,\big\lbrack D^0,\mathcal{L}_f\big\rbrack\cdots\big\rbrack.
\end{gather}
\end{Proposition}

In particular, the set of dif\/ferential operators $\mathcal{L}^H_f$, $f\in\Lambda_{\mathbb{F},\alpha}$, contains the CMS operator $\mathcal{L}^H$. Indeed, it is readily verif\/ied that
\begin{gather*}
    \mathcal{L}^H = \mathcal{L}^H_{-2\pi_{1,\alpha}};
\end{gather*}
note \eqref{HermiteOP} and \eqref{equalk}, and use the fact that $\mathcal{L}_{-2\pi_{1,\alpha}}=-2E^1$. We also note that if we substitu\-te~$\mathcal{L}_{f,n}$ for $\mathcal{L}_f$ and $D^0_n$ for $D^0$ in~\eqref{HermiteEigenOps}, then we obtain an eigenoperator $\mathcal{L}^H_{f,n}$ for the generalised Hermite polynomials, which satisf\/ies the intertwining relation
\begin{gather*}
    \varphi_n\circ\mathcal{L}_f = \mathcal{L}_{f,n}\circ\varphi_n.
\end{gather*}
Now,  in any dif\/ferential operator in $n$ variables,  let $\text{l.o.}$ denote terms of lower order.  Given that the Cherednik operators~$\xi_i$ satisfy $\xi_i=x_i\frac{\partial}{\partial x_i}+\mathrm{l.o.}$ (see for instance \cite[Section~4]{SV05}), it is clear from~\eqref{JackEigenOps} that
\begin{gather*}
    \mathcal{L}_{f,n} = f\left(x_1\frac{\partial}{\partial x_1},\ldots,x_n\frac{\partial}{\partial x_n}\right) + \text{l.o.}.
\end{gather*}
Given that $D^0_n$ is a dif\/ferential operator of order two, this implies that the order of  $[D^0_n,\mathcal{L}_{f,n}]$ is $k+1$, that of $[D^0_n,[D^0_n,\mathcal{L}_{f,n}]]$ is $k+2$, etc. In particular, this means that the order of $\mathcal{L}^H_{f,n}$, and therefore also of $\mathcal{L}^H_f$, is $2k$. As a consequence of Propositions~\ref{HermiteGenFuncProp} and~\ref{PropLfH}, we thus obtain the following corollary:

\begin{Corollary}\label{HermiteEigenOpsCorollary}
Let $f$ be as in Proposition~{\rm \ref{PropLfH}}. Then, $\mathcal{L}^H_f$ is a differential operator on $\Lambda_{\mathbf{F}}$
of order $2k$. Moreover, it is the unique operator on $\Lambda_{\mathbf{F}}$ such that
\begin{gather*}
    \mathcal{L}^H_f H_\lambda = f(\lambda)H_\lambda
\end{gather*}
for all partitions $\lambda$.
\end{Corollary}

\begin{proof}
There remains only to prove uniqueness, but this is immediate from the fact that the Hermite symmetric functions span $\Lambda_{\mathbf{F}}$; c.f., Corollary~\ref{basisCor}.
\end{proof}

Referring again to the fact that the Hermite symmetric functions form a basis for $\Lambda_{\mathbf{F}}$, we can conclude that the eigenoperators $\mathcal{L}^H_f$ pairwise commute.

\begin{Corollary}
We have that
\begin{gather*}
    \big[\mathcal{L}^H_f, \mathcal{L}^H_g\big]=0
\end{gather*}
for all $f,g\in\Lambda_{\mathbb{F},\alpha}$.
\end{Corollary}

We also note that the set of eigenoperators $\mathcal{L}^H_f$, $f\in\Lambda_{\mathbb{F},\alpha}$, separate the Hermite symmetric functions.

\begin{Lemma}\label{separationLemma}
For any two partitions $\lambda$ and $\mu$ such that $\lambda\neq\mu$, there exists $f\in\Lambda_{\mathbb{F},\alpha}$ such that $f(\lambda)\neq f(\mu)$.
\end{Lemma}

\begin{proof}
Let
\begin{gather*}
    \lambda^\alpha = \big(\lambda_1,\lambda_2-1/\alpha,\ldots,\lambda_i-(i-1)/\alpha,\ldots,\lambda_{\ell(\lambda)}-(\ell(\lambda)-1)/\alpha\big),
\end{gather*}
and similarly for $\mu^\alpha$. For any $f\in\Lambda_{\mathbb{F},\alpha}$, there exists a unique $p_f\in\Lambda_{\mathbb{F}}$ such that $f(\lambda) = p_f(\lambda^\alpha)$, and vice versa. If we expand $p_f(\lambda)$ in powers of $\alpha$, then we obtain
\begin{gather*}
    p_f(\lambda^\alpha) = p_f(\lambda) + \text{l.d.},
\end{gather*}
where $\text{l.d.}$ stands for terms of lower degree in $\alpha$. Since $\alpha$ is an indeterminate, we can conclude that $f(\lambda) = f(\mu)$ if and only if $p_f(\lambda) = p_f(\mu)$. The fact that the symmetric functions separate partitions thus implies the statement.
\end{proof}

\subsection{A limit from the Jacobi symmetric functions}
As indicated in the introduction, Sergeev and Veselov \cite{SV09} introduced and studied Jacobi symmetric functions as eigenfunctions of the dif\/ferential operator
\begin{gather}\label{JacobiDiffOp}
    \mathcal{L}^J = D^2 + 2D^1 - (p + 2q - 1)E^1 - (2p + 2q - 1)E^0.
\end{gather}
To make matters precise, let $\lambda$ be a partition. By Theorem \ref{TheoExistenceEigenfunction}, we can then def\/ine a corresponding Jacobi symmetric function $\mathcal{J}_\lambda(\alpha,p_0,p,q)$ as the unique eigenfunction of the dif\/ferential operator $\mathcal{L}^J$ that is of the form
\begin{gather}\label{JacobiEigenValue}
    \mathcal{J}_\lambda = P_\lambda + \sum_{\mu\subset\lambda}u_{\lambda\mu}P_\mu,\qquad u_{\lambda\mu}\in\Lambda_{\mathbf{F}(p,q)}.
\end{gather}
The associated eigenvalue is given by
\begin{gather*}
    e^J(\lambda) = \sum_i \lambda_i\left(\lambda_i + \frac{2}{\alpha}(p_0-i)\right) - (p+2q-1)|\lambda|.
\end{gather*}

We recall that Sergeev and Veselov used the parameter $k = -1/\alpha$, specif\/ied the form of the Jacobi symmetric functions in terms of the symmetric monomials $m_\lambda$, and f\/ixed the leading coef\/f\/icient to~$2^{|\lambda|}$. However, it is readily inferred from the triangular expansion \eqref{triangular} and the fact that the dominance order is compatible with the order given by inclusion of diagrams (see the discussion succeeding Theorem~\ref{TheoExistenceEigenfunction}) that the def\/inition given above is, up to a dif\/ference in normalisation, equivalent to that given by Sergeev and Veselov.

We note that the one-variable polynomials
\begin{gather*}
    \mathcal{J}_n(x)=\varphi_1(\mathcal{J}_{(n)}),\qquad n\in\mathbb{N}_0,
\end{gather*}
have a somewhat non-standard form. Indeed,
\begin{gather*}
    \varphi_1\circ\mathcal{L}^J = \left(x(x+2)\frac{d^2}{dx^2} - \big((p+2q-1)x+2p+2q-1\big)\frac{d}{dx}\right)\circ\varphi_1,
\end{gather*}
and the $\mathcal{J}_n(x)$ can be seen to form a sequence of orthogonal polynomials on the interval $\lbrack-2,0\rbrack$ with respect to the weight function
\begin{gather*}
    w(x)=x^{-p-q-1/2}(2+x)^{-q-1/2}
\end{gather*}
for appropriate parameter values. In order to obtain a more standard one-variable restriction, we can instead start from the symmetric functions
\begin{gather*}
    J_\lambda:= \frac{1}{(-2)^{|\lambda|}}\sigma_{-2}\big(\mathcal{J}_\lambda\big),
\end{gather*}
which are (monic) eigenfunctions of $\sigma_{-2}\circ\mathcal{L}^J\circ \sigma_{-1/2}$. Using \eqref{sigmaD}, we f\/ind that
\begin{gather*}
    \varphi_1\circ\big(\sigma_{-2}\circ\mathcal{L}^J\circ\sigma_{-1/2}\big)  = \varphi_1\big(D^2-D^1-(p+2q-1)E^1+(p+q-1/2)E^0\big)\\
     \hphantom{\varphi_1\circ\big(\sigma_{-2}\circ\mathcal{L}^J\circ\sigma_{-1/2}\big)}{}
     = \left(x(x-1)\frac{d^2}{dx^2} - \big((p+2q-1)x-p-q+1/2\big)\frac{d}{dx}\right)\circ \varphi_1,
\end{gather*}
and it is readily inferred that the polynomials $J_n(x):=\varphi_1(J_{(n)})$ are orthogonal on the interval $\lbrack 0,1\rbrack$ with respect to the weight function
\begin{gather*}
    \tilde{w}(x) = x^{-p-q-1/2}(x-1)^{-q-1/2}.
\end{gather*}
Moreover, after a suitable reparameterisation and renormalisation, the symmetric polynomials $J_\lambda(x_1,\ldots,x_n):=\varphi_n(J_\lambda)$ coincide with the generalised Jacobi polynomials, as considered by Lassalle~\cite{Las91a}, Macdonald~\cite{Mac} and also by Baker and Forrester~\cite{BF97}.

We shall now use a standard method to obtain the Hermite symmetric functions as a particular limit of the Jacobi symmetric functions. In order to do so, we shall again work with symmetric functions over real numbers, i.e., with $\Lambda_{\mathbb{R}}$, and thus assume that $\alpha,p_0,p,q\in\mathbb{R}$. The starting point is the representation \eqref{FRep}, which in this case yields
\begin{gather}\label{JacobiRep}
    \mathcal{J}_\lambda = \prod_{\mu\subset\lambda}\frac{\mathcal{L}^J-e^J(\mu)}{e^J(\lambda)-e^J(\mu)}P_\lambda.
\end{gather}
If we replace $\mathcal{L}^J$ by $\mathcal{L}^H$ and $e^J(\lambda)$ by $-2|\lambda|$, then we obtain the corresponding representation for the Hermite symmetric functions $H_\lambda$.

 We now   introduce a homomorphism $t_\gamma: \Lambda_{\mathbf{F}}\to\Lambda_{\mathbf{F}}$, $\gamma\in\mathbb{F}$, by setting
\begin{gather}\label{tgammaDef}
    t_\gamma(p_r) = \sum_{m=0}^r\gamma^{r-m}\binom{r}{m}p_m,\qquad r\geq 1.
\end{gather}
For a f\/inite number of indeterminates $x=(x_1,\ldots,x_n)$, this simply yields the translation of each~$x_i$ by~$\gamma$.
It follows from Lemma~\ref{tCommuteLemma} and~\eqref{sigmaD}, that
\begin{gather*}
    \big(\sigma_{(-q)^{-1/2}}\circ t_{-1}\big)\circ\mathcal{L}^J
  =  \big(D^2 + qD^0 - (p+2q-1)E^1 - (-q)^{1/2} pE^0\big)\circ\big(\sigma_{(-q)^{-1/2}}\circ t_{-1}\big).
\end{gather*}
We note that, by Lemma~\ref{EDLemma}, Lemma~\ref{degreeLemma} and (the obvious analogue for $\Lambda_{\mathbb{R}}$ of) Lemma~\ref{secondContinuityLemma}, the dif\/ferential operator $\mathcal{L}^J$ is continuous with respect to the topology of term-wise convergence. Combining the observations above with the binomial formula in Proposition~\ref{genBinomProp}, as well as the fact that Jack's symmetric functions~$P_\lambda$ are homogeneous of degree~$|\lambda|$, we readily deduce the following proposition:

\begin{Proposition}\label{JacobiToHermiteProp}
Let $\lambda$ be a partition. Then, for generic parameter values, we have that
\begin{gather*}
    H_\lambda(\alpha,p_0)
    = \lim_{q\rightarrow\infty}(-q)^{|\lambda|/2}\big(\sigma_{(-q)^{-1/2}}\circ t_{-1}\big)\big(\mathcal{J}_\lambda(\alpha,p_0,p,q)\big)
\end{gather*}
in the sense of term-wise convergence.
\end{Proposition}

\begin{Remark}
By generic we mean on a dense set in parameter space with respect to the Zariski topology. The validity of this part of the statement is a direct consequence of the fact that the Jacobi symmetric functions~$J_\lambda(\alpha,p_0,p,q)$ depend rationally on all parameters.
\end{Remark}

\subsection{Structure of Pieri formulae and invariant ideals}\label{section4.4}

Throughout this section we shall assume $p_0\in\mathbb{F}$ f\/ixed. The main purpose is to obtain the ideals $I\subset\Lambda_{\mathbb{F}}$ that are invariant under the action of all dif\/ferential operators $\mathcal{L}^{H}_f$, $f\in\Lambda_{\mathbb{F},\alpha}$. This is the case if and only if~$I$ has a basis consisting of Hermite symmetric functions. The f\/irst part of this claim is trivial, while the second part is a consequence of Lemma~\ref{separationLemma}. For future reference, we state this fact in the form of a~lemma.

\begin{Lemma}\label{IdealLemma}
Let $I\subset\Lambda_{\mathbb{F}}$ be an ideal such that $\mathcal{L}^H_f I\subset I$ for all $f\in\Lambda_{\mathbb{F},\alpha}$. Then, we have that
\begin{gather*}
    I = \mathbb{F}\big\langle H_\lambda: \lambda\in\text{\rm Par}_I\big\rangle
\end{gather*}
for some set of partitions~$\text{\rm Par}_I$.
\end{Lemma}

We proceed to deduce Pieri type recurrence relations for the Hermite symmetric functions. In the Jacobi case, Sergeev and Veselov \cite{SV09} (see Theorem 4.4) obtained such recurrence relations by generalising corresponding recurrence relations for generalised Jacobi polynomials due to van Diejen \cite{vD99} (see Theorem 6.4). However, in this generalisation part of the explicit nature of van Diejen's formulae were lost. In fact, for our purposes, we require a more explicit version of Sergeev and Veselov's result, stated below in Theorem \ref{JacobiPieriFormulaeThm}. By applying the limit transition in Proposition \ref{JacobiToHermiteProp} we shall then obtain the desired recurrence relations for the Hermite symmetric functions.

An important ingredient is the specialisation formula of the Jacobi symmetric functions at $p_r=0$, $r\in\mathbb{N}$, as deduced by Sergeev and Veselov \cite[Proposition 4.3]{SV09} from  the analogue formula for the f\/inite-dimensional case (i.e., $p_0=n$).  The latter can be obtained from  Corollary~5.2 \cite{Op89} by specialising  to the root system~$BC_n$.

\begin{Proposition}[Sergeev and Veselov~\protect{\cite{SV09}}]\label{JacobiSpecialisationProp}
For any partition $\lambda$, let
\begin{subequations}\label{CDef}
\begin{gather}
    C^+_\lambda(z;\alpha)  = \prod_{(i,j)\in\lambda}\big(\lambda_i+j-(\lambda^\prime_j+i)/\alpha+z\big),\\
\label{Cminus}
    C^-_\lambda(z;\alpha)  = \prod_{(i,j)\in\lambda}\big(\lambda_i-j+(\lambda^\prime_j-i)/\alpha+z\big),\\
\label{Czero}
    C^0_\lambda(z;\alpha)  = \prod_{(i,j)\in\lambda}\big(j-1-(i-1)/\alpha+z\big).
\end{gather}
\end{subequations}
Then, we have that
\begin{gather}\label{JacobiSpecialisation}
    \epsilon_0\big(\mathcal{J}_\lambda(\alpha,p_0,p,q)\big) = 2^{|\lambda|}\frac{C^0_\lambda(p_0/\alpha)C^0_\lambda\big((p_0-1)/\alpha-p-q+1/2\big)}{C^-_\lambda(1/\alpha)C^+_\lambda(2p_0/\alpha-p-2q-1)}.
\end{gather}
\end{Proposition}

For $m\in\mathbb{N}$, we let $I(m)$ denote the set consisting of the $m$ smallest non-negative integers, i.e.,
\begin{gather*}
    I(m) = \lbrace 1,\ldots,m \rbrace\subset\mathbb{N}.
\end{gather*}
Given any subset $J\subseteq\mathbb{N}$, and corresponding sequence $\epsilon(J) = \lbrace\epsilon_j\rbrace_{j\in J}$ of signs $\epsilon_j = \pm 1$, $j\in J$, we let $\lambda+e_{\epsilon(J)}$ denote the sequence def\/ined by
\begin{gather*}
    (\lambda + e_{\epsilon(J)})_i = \lambda_i + \epsilon_i,\qquad i\in\mathbb{N},
\end{gather*}
where we have set $\epsilon_i = 0$ if $i\notin J$. With this notation in mind, we are now ready to state the recurrence relations for the Jacobi symmetric functions in a form that is convenient for our purposes.

\begin{Theorem}\label{JacobiPieriFormulaeThm}
Let $J\subseteq I\subset\mathbb{N}$ be two finite subsets of the set of $($positive$)$ natural numbers $\mathbb{N}$, and fix a sequence $\epsilon(J) = \lbrace\epsilon_j\rbrace_{j\in J}$ of signs $\epsilon_j = \pm 1$, $j\in J$. Introduce the rational function
\begin{gather*}
    R_{\epsilon(J)}(z;m) = \prod_{j\in J}\frac{(\epsilon_j z_j+z_m+1/\alpha)(\epsilon_j z_j+p/2+q+1/\alpha)}{(\epsilon_j z_j-p/2-q)(\epsilon_j z_j-z_m)}.
\end{gather*}
Let, furthermore,
\begin{gather*}
    \hat{v}^J(z) = \frac{z+1/\alpha}{z},\qquad \hat{w}^J(z) = \frac{(z-p/2-q)(z+(1-p)/2)}{z(z+1/2)},
\end{gather*}
and introduce the following two rational functions:
\begin{subequations}\label{V}
\begin{gather}
    \hat{V}^{(+)}_{I,\epsilon(J)}(z)  = \prod_{j\in J}\hat{w}^J(\epsilon_j z_j)\prod_{\substack{j,j^\prime\in J\\ j<j^\prime}}\hat{v}^J(\epsilon_j z_j + \epsilon_{j^\prime}z_{j^\prime})\hat{v}^J(\epsilon_j z_j + \epsilon_{j^\prime}z_{j^\prime} + 1)\nonumber \\
\hphantom{\hat{V}^{(+)}_{I,\epsilon(J)}(z)  =}{}
\times \prod_{\substack{j\in J\\ i\in I\setminus J}}\hat{v}^J(\epsilon_j z_j + z_i)\hat{v}^J(\epsilon_j z_j - z_i),\label{Vplus}\\
    \hat{V}^{(-)}_{I,\epsilon(J)}(z)  = \prod_{j\in J}\hat{w}^J(\epsilon_j z_j)\prod_{\substack{j,j^\prime\in J\\ j<j^\prime}}\hat{v}^J(\epsilon_j z_j + \epsilon_{j^\prime}z_{j^\prime})\hat{v}^J(-\epsilon_j z_j - \epsilon_{j^\prime}z_{j^\prime} - 1)\nonumber \\
\hphantom{\hat{V}^{(-)}_{I,\epsilon(J)}(z)  =}{}
     \times \prod_{\substack{j\in J\\ i\in I\setminus J}}\hat{v}^J(\epsilon_j z_j+z_i)\hat{v}^J(\epsilon_j z_j-z_i).\label{Vminus}
\end{gather}
\end{subequations}
To each $z\in\mathbb{F}$, associate the sequence
\begin{gather*}
    \rho^J(z) = \big\lbrace (z-i)/\alpha-p/2-q\big\rbrace_{i\in\mathbb{N}}.
\end{gather*}
For each $r\in\mathbb{N}$, let {
\begin{gather*}
    E_r = 2^rm_{(1^r)},
\end{gather*}
that is, $ E_r$ is equal to $2^r$ times the $r$th elementary symmetric function.}
Then, the re-normalised Jacobi symmetric functions $\mathcal{J}_\lambda/\epsilon_0(\mathcal{J}_0)$ satisfy, for generic values of $p_0$, the recurrence relations
\begin{gather}
    E_r\frac{\mathcal{J}_\lambda}{\epsilon_0(\mathcal{J}_\lambda)} = \sum_{\epsilon(J),\epsilon(K)}(-1)^{|K|}\hat{V}^{(+)}_{I(\ell(\lambda)+r),\epsilon(J)}\big(\rho^J(p_0)
    +\lambda\big)\hat{V}^{(-)}_{I(\ell(\lambda)+r)\setminus J,\epsilon(K)}\big(\rho^J(p_0)+\lambda\big)\label{JacobiPieri}\\
  \hphantom{E_r\frac{\mathcal{J}_\lambda}{\epsilon_0(\mathcal{J}_\lambda)} =}  {}
  \times R_{\epsilon(J)}\big(\rho^J(p_0)+\lambda;\ell(\lambda)+r+1\big) R_{\epsilon(K)}\big(\rho^J(p_0)+\lambda;\ell(\lambda)+r+1\big)  \frac{\mathcal{J}_{\lambda+e_{\epsilon(J)}}}{\epsilon_0(\mathcal{J}_{\lambda+e_{\epsilon(J)}})},\nonumber
\end{gather}
where the sum extends over all sequences of signs $\epsilon(J)$ and $\epsilon(K)$ with $J,K\subset I(\ell(\lambda)+r)$ such that $J\cap K = \varnothing$, $|J|+|K| = r$, and $\lambda + e_{\epsilon(J)}$ is a partition.
\end{Theorem}

\begin{Remark}
It is clear from the representation \eqref{FRep} that $\mathcal{J}_\lambda$, and thereby also $\epsilon_{p_0}(\mathcal{J}_\lambda)$, is a~rational function of $p_0$; c.f., \eqref{JacobiDiffOp} and~\eqref{JacobiEigenValue}. It follows that $\epsilon_{p_0}(\mathcal{J}_\lambda)\neq 0$ on a dense (open) set in the Zariski topology. It is for these `generic' values of $p_0$ that the recurrence relations~\eqref{JacobiPieri} are valid.
\end{Remark}

\begin{Remark}
As discussed above, Theorem \ref{JacobiPieriFormulaeThm} is the inf\/inite-dimensional generalisation of a~result of van Diejen~\cite[Theorem~6.4]{vD99} and is essentially due to Sergeev and Veselov~\cite{SV09} (see Theorem~4.4)~-- with the dif\/ference that the latter authors did not provide an explicit formula for the coef\/f\/icients in~\eqref{JacobiPieri}. We shall require this explicit information in order to obtain corresponding recurrence relations for the Hermite symmetric functions. For the convenience of the reader, we have included a full proof of Theorem~\ref{JacobiPieriFormulaeThm} in Appendix~\ref{JacobiPieriFormulaeThmProof}, expanding on the proof of Theorem~4.4 in Sergeev and Veselov~\cite{SV09}.
\end{Remark}

It is important to note that we can not just simply apply $\sigma_{(-q)^{-1/2}}\circ t_{-1}$ to~\eqref{JacobiPieri}, and then take the limit $q\rightarrow\infty$, as in Proposition~\ref{JacobiToHermiteProp}. Indeed, for $r>1$, the symmetric function $(-q)^{(\ell(\lambda)+r)/2}(\sigma_{(-q)^{-1/2}}\circ t_{-1})({E}_r\mathcal{J}_\lambda(\alpha,p_0,p,q))$ contains terms which diverge as $q\rightarrow\infty$. However, this problem can be resolved by considering instead appropriate linear combinations of the recurrence relations~\eqref{JacobiPieri}. For example, if we are interested in the case $r = 2$, then we should observe that
\begin{gather*}
    t_{-1}\big(E_2+2(p_0-1)E_1+2p_0(p_0-1)\big) = E_2,
\end{gather*}
and consider the corresponding linear combination of recurrence relations~\eqref{JacobiPieri}. For a detailed discussion of this point, in the context of a f\/inite number of variables, see van Diejen~\cite{vD97}.

Another issue, which is one of convenience rather than necessity, is the choice of normalisation of the Hermite symmetric functions. In order to f\/ind the normalisation for which the corresponding recurrence relations take the simplest possible form, we note the $q\to\infty$ limit of the normalisation factors $\epsilon_0(\mathcal{J}_\lambda)$:
\begin{gather}\label{normFactorLimit}
    \lim_{q\rightarrow\infty} \epsilon_0\big(\mathcal{J}_\lambda(\alpha,p_0,p,q)\big) = \frac{C^0_\lambda(p_0/\alpha)}{C^-_\lambda(1/\alpha)} = \epsilon_{p_0}(P_\lambda),
\end{gather}
where the second equality follows from a direct comparison of \eqref{epsilonXJackEq} and \eqref{Cminus}, \eqref{Czero}. As will become clear below, it will be convenient to extract from this limit the factor $C^0_\lambda(p_0/\alpha)$, which contains all the dependence on the parameter $p_0$, and re-normalise the Hermite symmetric functions by the factor $C^-_\lambda(1/\alpha)$ only.

We shall make use of the following notation: given a subset $J\subset\mathbb{N}$, and a corresponding sequence of signs $\epsilon(J)$, we shall write $J_+$ and $J_-$ for the subsets of $J$ given by
\begin{gather*}
    J_+ = \lbrace j\in J: \epsilon_j = +1\rbrace,\qquad J_- = \lbrace j\in J: \epsilon_j = -1\rbrace.
\end{gather*}

With the above remarks in mind, we continue by stating and proving the analogy of Theo\-rem~\ref{JacobiPieriFormulaeThm} for the Hermite symmetric functions.

\begin{Proposition}\label{HermitePieriFormulaeProp}
The re-normalised Hermite symmetric functions
\begin{gather*}
    \mathcal{H}_\lambda:= C^-_\lambda(1/\alpha)H_\lambda
\end{gather*}
satisfy recurrence relations of the form
\begin{gather}\label{HermitePieri}
    e_r\mathcal{H}_\lambda = \sum_{J_+,J_-}\hat{W}_{I(\ell(\lambda)+r);J_+,J_-}(\lambda)\mathcal{H}_{\lambda+e_{J_+}-e_{J_-}},
\end{gather}
where the sum is over all subsets $J_+,J_-\subset\mathbb{N}$ such that $J_+\cap J_- = \varnothing$, $|J_+|+|J_-|\leq r$, $r-|J_+|-|J_-|$ is even, and $\lambda+e_{J_+}-e_{J_-}$ is a partition.

Moreover, the coefficients $\hat{W}_{I(\ell(\lambda)+r);J_+,J_-}$ are of the form
\begin{gather}\label{PieriCoeffs}
    \hat{W}_{I(\ell(\lambda)+r);J_+,J_-}(\lambda) = \frac{1}{2^{|J_-|}}\prod_{j\in J_-}\left(\frac{p_0-j+1}{\alpha}+\lambda_j-1\right)\hat{U}_{I(\ell(\lambda)+r);J_+,J_-}(\lambda),
\end{gather}
where ${U}_{I(\ell(\lambda)+r);J_+,J_-}$ is a polynomial in $p_0$, and if  $|J_+|+|J_-| = r$, then
\begin{gather}
    \hat{U}_{I(\ell(\lambda)+r);J_+,J_-}(\lambda)  = \!\!\prod_{j\in J_+,j^\prime\in J_-}\!\!\left(1+\frac{1}{j^\prime-j+\alpha(\lambda_j-\lambda_{j^\prime})}\right)
    \left(1+\frac{1}{j^\prime-j+\alpha(\lambda_j-\lambda_{j^\prime}+1)}\right)\nonumber\\
\hphantom{\hat{U}_{I(\ell(\lambda)+r);J_+,J_-}(\lambda)  =}{}
    \times \! \prod_{j\in J_-} \! \big(\ell(\lambda)+r-j+\alpha\lambda_j\big)\!\!\prod_{i\in I(\ell(\lambda)+r)\setminus J}\!\!\left(1+\frac{1}{j-i+\alpha(\lambda_i\!-\lambda_j)}\right)\!\!\!\!\!\label{U}\\
\hphantom{\hat{U}_{I(\ell(\lambda)+r);J_+,J_-}(\lambda)  =}{}
    \times\!\prod_{j\in J_+}\!\frac{1}{\ell(\lambda)+r+1-j+\alpha\lambda_j} \!\prod_{i\in I(\ell(\lambda)+r)\setminus J}\!\!\left(1-\frac{1}{j-i+\alpha(\lambda_i-\lambda_j)}\right).\nonumber
\end{gather}
\end{Proposition}

\begin{proof}
As noted above, for a unique set of coef\/f\/icients $c_r,\ldots,c_0\in\Lambda_{\mathbb{F}}$, we have that
\begin{gather*}
    t_{-1}(c_r E_r + c_{r-1}E_{r-1} +\cdots+ c_0) = E_r.
\end{gather*}
Consider the corresponding linear combination of recurrence relations \eqref{JacobiPieri}. For the left-hand side of the resulting relation, Proposition \ref{JacobiToHermiteProp} and \eqref{normFactorLimit} yield the limit
\begin{gather*}
    \lim_{q\to\infty}(-q)^{(\ell(\lambda)+r)/2}\big(\sigma_{(-q)^{-1/2}}\circ t_{-1}\big)\left((c_r E_r + c_{r-1}E_{r-1} +\cdots+ c_0)\frac{\mathcal{J}_\lambda}{\epsilon_0(\mathcal{J}_\lambda)}\right)  = E_r\frac{H_\lambda}{\epsilon_{p_0}(P_\lambda)}.
\end{gather*}

Furthermore, it is clear that the limit of the right-hand side is of the form
\begin{gather*}
    2^r\sum_{J_+,J_-}\hat{W}_{I(\ell(\lambda)+r);J_+,J_-}(\lambda)\frac{H_{\lambda+e_{J_+}-e_{J_-}}}{\epsilon_{p_0}(P_{\lambda+e_{J_+}-e_{J_-}})}
\end{gather*}
with the coef\/f\/icients $\hat{W}_{I(\ell(\lambda)+r);J_+,J_-}$ given by
\begin{gather}
    \lim_{q\to\infty}\frac{(-q)^{(r-|J_+|+|J_-|)/2}}{2^r}\frac{C^0_\lambda(p_0/\alpha)}
    {C^0_{\lambda+e_{J_+}-e_{J_-}}(p_0/\alpha)}\sum_{K_+,K_-}(-1)^{|K|}\nonumber\\
    \qquad{} \times\hat{V}^{(+)}_{I(\ell(\lambda)+r),\epsilon(J)}\big(\rho^J(p_0)+\lambda\big)
    R_{\epsilon(J)}\big(\rho^J(p_0)+\lambda;\ell(\lambda)+r+1\big)\nonumber\\
    \qquad{} \times\hat{V}^{(-)}_{I(\ell(\lambda)+r)\setminus J,\epsilon(K)}\big(\rho^J(p_0)+\lambda\big) R_{\epsilon(K)}\big(\rho^J(p_0)+\lambda;\ell(\lambda)+r+1\big).\label{coeffLimit}
\end{gather}
As a direct computation shows, we have that
\begin{gather}
    \lim_{q\rightarrow\infty}\!R_{\epsilon(L)}\big(\rho^J(p_0)+\lambda;m\big)  = \!\prod_{j\in L_+}\!\!\frac{p_0-j+1+\alpha\lambda_j}{m-j+\alpha(\lambda_j-\lambda_m)} \! \prod_{j\in L_-}\!\!\frac{m-j-1+\alpha(\lambda_j-\lambda_m)}{p_0-j+\alpha\lambda_j},\label{RLim}
\end{gather}
and that
\begin{gather*}
    \lim_{q\rightarrow\infty}(-q)^{|L_-|}\prod_{j\in L}\hat{w}^J\big(\epsilon_j(\rho^J(p_0)+\lambda)_j\big) = 2^{|L_-|}\prod_{j\in L_-}\left(\frac{p_0-j}{\alpha}+\lambda_j\right)
\end{gather*}
for $L=J,K$. We observe that, for all arguments $z$ appearing in \eqref{JacobiPieri}, $\hat{v}(z)$ is a bounded function of $q$; c.f., \eqref{V}. It follows that a given term in \eqref{coeffLimit} provides a non-zero contribution only if
\begin{gather*}
    \frac{r-|J_+|-|J_-|}{2} - |K_-| = 0,
\end{gather*}
which clearly can hold true only if $r-|J_+|-|J_-|$ is even. This concludes the proof of the f\/irst part of the statement.

In order to establish the stated structure of the coef\/f\/icients $\hat{W}_{I(\ell(\lambda)+r);J_+,J_-}$, we observe that
\begin{gather*}
    \frac{C^0_\lambda(p_0/\alpha)}{C^0_{\lambda+e_{J_+}-e_{J_-}}(p_0/\alpha)} = \frac{\prod\limits_{j\in J_-}\big(p_0-j+1+\alpha(\lambda_j-1)\big)}{\prod\limits_{j\in J_+}\big(p_0-j+1+\alpha\lambda_j\big)}.
\end{gather*}
If we now set $m=\ell(\lambda)+r+1$ (c.f., \eqref{JacobiPieri}), and combine the observations made thus far, we readily deduce  \eqref{PieriCoeffs}. Moreover, in case $|J_+|+|J_-|=r$, we have $K=\emptyset$. It follows that the sum in \eqref{coeffLimit} contains only one term, and a direct computation yields \eqref{U}.
\end{proof}

\begin{Remark}
From the representation \eqref{FRep} we can directly infer that $H_\lambda$, and thereby also~$\mathcal{H}_\lambda$, is a polynomial in~$p_0$; c.f., Def\/inition~\ref{HermiteDef}. In contrast to the Jacobi case, this entails that the recurrence relations~\eqref{HermitePieri} are valid not only for generic but indeed all values of the parameter~$p_0$.
\end{Remark}

If we restrict our attention to $r=1$, then the statement can be simplif\/ied considerably. In particular, all coef\/f\/icients can be specif\/ied explicitly.

\begin{Corollary}\label{FirstHermitePieriCor}
The re-normalised Hermite symmetric functions $\mathcal{H}_\lambda$ satisfy the recurrence relation
\begin{gather}\label{FirstHermitePieri}
    e_1\mathcal{H}_\lambda = \sum_{j=1}^{{ \ell(\lambda)}+1}\big(\hat{W}_j(\lambda)\mathcal{H}_{\lambda+e_j} + \hat{W}_{-j}(\lambda)\mathcal{H}_{\lambda-e_j}\big)
\end{gather}
with the coefficients
\begin{gather*}
    \hat{W}_j(\lambda)  = \frac{1}{\ell(\lambda)+2-j+\alpha\lambda_j}\prod_{\substack{1\leq i\leq{ \ell(\lambda)}+1\\ i\neq j}}\left(1-\frac{1}{j-i+\alpha(\lambda_i-\lambda_j)}\right),\\
    \hat{W}_{-j}(\lambda)  = \frac{1}{2}\left(\frac{p_0-j+1}{\alpha}+\lambda_j-1\right)\big(\ell(\lambda)+1-j+\alpha\lambda_j\big)\\
    \hphantom{\hat{W}_{-j}(\lambda)  =}{}
    \times\prod_{\substack{1\leq i\leq{ \ell(\lambda)}+1\\ i\neq j}}\left(1+\frac{1}{j-i+\alpha(\lambda_i-\lambda_j)}\right).
\end{gather*}
\end{Corollary}

\begin{Remark}
 When restricted to the polynomial case,  Corollary \ref{FirstHermitePieriCor} and Proposition \ref{PieriCoeffs} respectively reduce to Propositions~2.5 and~2.6 of \cite{vD97}.  The f\/inite-dimensional analogue of Corollary~\ref{FirstHermitePieriCor} can also be found in \cite[Proposition~3.5]{BF97}.
\end{Remark}

We proceed to consider how the recurrence relation~\eqref{FirstHermitePieri} is related to the question of exis\-tence of invariant ideals. To this end, let $I\subset\Lambda_{\mathbb{F}}$ be an ideal invariant under the dif\/ferential operators $\mathcal{L}^H_f$, $f\in\Lambda_{\mathbb{F},\alpha}$. By Lemma \ref{IdealLemma}, there exists at least one partition $\lambda$ such that $H_\lambda\in I$. In case $\lambda=(0)$, we have that $H_\lambda=1$, and consequently that $I=\Lambda_{\mathbb{F}}$. Suppose that $\lambda\neq 0$. Then, we can always f\/ind an integer $j=1,\ldots,|\lambda|$ such that $\lambda-e_j$ is a partition. For example, $j=|\lambda|$. It is clear from Corollary~\ref{FirstHermitePieriCor} that
\begin{gather*}
    \hat{W}_{-j}(\lambda) = 0
\end{gather*}
if and only if
\begin{gather*}
    p_0=j-1-\alpha(\lambda_j-1).
\end{gather*}
Hence, if $p_0$ is not of this form, then we can conclude that also $H_{\lambda-e_j}\in I$. Moreover, assuming this to be the case for all partitions $\mu\subset\lambda$, it follows by induction on the weight $|\lambda|$ of $\lambda$ that again $1\in I$. This observation forms one part of the main result of this section, as stated in the following theorem:

\begin{Theorem}\label{IdealThm}
$\Lambda_{\mathbb{F}}$ contains a non-trivial ideal invariant under all differential operators $\mathcal{L}^H_f$, $f\in\Lambda_{\mathbb{F},\alpha}$, for and only for non-zero $p_0$ of the form
\begin{gather}\label{p0ForIdeal}
    p_0 = n-\alpha m,\qquad n,m\in\mathbb{N}_0.
\end{gather}
If that is the case, then there is a unique such ideal, spanned by the Hermite symmetric func\-tions~$H_\lambda$ labelled by the partitions $\lambda$ such that $(n+1,m+1)\in\lambda$.
\end{Theorem}

\begin{proof}
There remains only to prove that, for $p_0 = n-\alpha m$, there exists a unique non-trivial invariant ideal
\begin{gather*}
    I = \mathbb{F}\langle H_\lambda: (n+1,m+1)\in\lambda\rangle\subset\Lambda_{\mathbb{F}}.
\end{gather*}
Suppose that a partition $\lambda$ and subsets $J_+,J_-\subset\mathbb{N}$ are such that $(n+1,m+1)\in\lambda$ but $(n+1,m+1)\notin\lambda+e_{J_+}-e_{J_-}$. Clearly, this is only possible if $\lambda_{n+1} = m+1$, and $n+1\in J_-$. Fix $r\in\mathbb{N}$, and consider the corresponding recurrence relation~\eqref{HermitePieri}. Since $\hat{W}_{I(\ell(\lambda)+r);J_+,J_-}$ is a~polynomial in~$p_0$, it follows from \eqref{PieriCoeffs} that
\begin{gather*}
    \hat{W}_{I(\ell(\lambda)+r);J_+,J_-}(\lambda) = 0.
\end{gather*}
Hence, $I$ is indeed an ideal, which, by construction, is invariant. In order to establish uniqueness, we let $I^\prime\subset\Lambda_{\mathbb{F}}$ be any non-trivial invariant ideal. Suppose that $H_\lambda\in I^\prime$ for some partition $\lambda$ such that $(n+1,m+1)\notin\lambda$. Then, by following the discussion preceding the Theorem, we obtain that $1\in I^\prime$. We must therefore have $I^\prime\subseteq I$. In addition, we can exclude the possibility that $I^\prime\neq I$ by observing that, starting from any $\mathcal{H}_\lambda\in I$, we can obtain any other Hermite symmetric function $\mathcal{H}_\mu$ by applying the recurrence relation~\eqref{FirstHermitePieri}.
\end{proof}

\section{Laguerre symmetric functions}\label{LaguerreSec}
In this section we shall introduce and study Laguerre symmetric functions as eigenfunctions of the dif\/ferential operator
\begin{gather}\label{LaguerreOp}
\mathcal{L}^L=D^1+(a+1)E^0-\nu E^1
\end{gather}
with $a$ an indeterminate, and $\nu\in\mathbf{F}$. We shall follow closely our treatment of the Hermite symmetric functions in Section~\ref{HermiteSec}. To avoid unnecessary repetitions, the discussion will therefore be brief, and statements that can be obtained as straightforward generalisations from the Hermite case will be stated without proofs.

It is clear from Lemma \ref{actionLemma} that
\begin{gather*}
    \mathcal{L}^LP_\lambda = -\nu|\lambda|P_\lambda + \sum_{\mu\subset\lambda}c_{\lambda\mu}P_\mu
\end{gather*}
for some coef\/f\/icients $c_{\lambda\mu}\in\mathbf{K}$; c.f., \eqref{fieldExtensions}. Theorem \ref{TheoExistenceEigenfunction} thus guarantees that we can make the following def\/inition:

\begin{Definition}\label{LaguerreDef}\sloppy
Let $\lambda$ be a partition. We then def\/ine the Laguerre symmetric function $L_\lambda(\alpha,p_0,a,\nu)$ as the unique symmetric function such that
\begin{enumerate}\itemsep=0pt
\item $L_\lambda = P_\lambda + \sum\limits_{\mu\subset\lambda}u_{\lambda\mu}P_\mu$ for some $u_{\lambda\mu}\in\mathbf{K}$,
\item $\mathcal{L}^LL_\lambda=-\nu|\lambda|L_\lambda$.
\end{enumerate}
\end{Definition}

Following the proof of Proposition \ref{PropLassalle}, we can establish a constructive def\/inition.

\begin{Proposition}\label{propexplaguerre}
Let $\lambda$ be a partition and set $L=|\lambda|$.  Then, we have that
\begin{gather}\label{LaguerreExpRep}
    L_\lambda = \exp_L\left(-\frac{1}{\nu}\big(D^1+(a+1)E^0\big)\right)(P_\lambda).
\end{gather}
\end{Proposition}

\subsection{A symmetry property}
Here, we obtain a symmetry property of the Laguerre symmetric functions that does not have a~counterpart in the case of symmetric polynomials. This fact follows easily once we have shown that the second-order CMS operator~\eqref{LaguerreOp} posess the symmetry in question.

We write $D^1(\alpha,p_0)$ and $E^0(p_0)$ for the operators $D^1$ and $E^0$ of Def\/inition~\ref{EDDef}, and let
\begin{gather}\label{eqap0}
    a'=2/\alpha-a-2,\qquad p_0'=p_0-1+\alpha(a+1).
\end{gather}
Then, by expanding the series def\/ining $D^1(\alpha,p_0^\prime)$ and $E^0(\alpha,p_0^\prime)$ and concentrating on the terms depending on $p_0$, one readily verif\/ies that
\begin{gather*}
    D^1(\alpha,p_0') =D^1(\alpha,p_0)+\frac{2}{\alpha}\big(\alpha(a+1)-1\big)E^0(p_0) +\frac{1}{\alpha}\big(\alpha(a+1)-1\big)\big(\alpha(a+1)-2\big)\partial(p_1),\\
    E^0(p_0') =E^0(p_0)+\big(\alpha(a+1)-1\big)\partial(p_1).
\end{gather*}
These equations, together with  the independence of $E^1$ upon $p_0$, imply the following symmetry property:
\begin{gather*}
    D^1(\alpha,p_0')+(a'+1)E^0(p_0')-\nu E^1=D^1(\alpha,p_0)+(a+1)E^0(p_0)-\nu E^1.
\end{gather*}
Finally, a direct comparison of this equation with Def\/inition \ref{LaguerreDef} or Proposition \ref{propexplaguerre} establishes the desired symmetry relation.

\begin{Proposition}\label{propsymlaguerre}
Let $a'$ and $p_0'$ be as in \eqref{eqap0}.  Then, we have that
\begin{gather*}
    L_\lambda(\alpha,p_0',a',\nu) = L_\lambda(\alpha,p_0,a,\nu).
\end{gather*}
\end{Proposition}

\subsection{A duality relation}
We continue by considering the action of the homomorphism $\omega_\alpha$ on the Laguerre symmetric functions; c.f., \eqref{omegaDef}. Starting from the representation~\eqref{LaguerreExpRep}, a direct application of Lemma~\ref{lemmadualityDE} yields
\begin{gather*}
    \omega_\alpha(L_\lambda(\alpha,p_0,a,\nu)) = \exp_L\left(-\frac{1}{\alpha\nu}\big(D^1(1/\alpha,-\alpha p_0)+(1-\alpha a )E^0\big)\right)\big(Q_{\lambda'}(1/\alpha)\big).
\end{gather*}
We thus arrive at the following proposition:

\begin{Proposition}\label{LaguerreDualityProp}
We have that
\begin{gather*}
    \omega_\alpha\big(L_\lambda(\alpha,p_0,a,\nu)\big)= b_{\lambda'}(1/\alpha) L_{\lambda'}(1/\alpha,-\alpha p_0,-\alpha a,-\alpha \nu).
\end{gather*}
\end{Proposition}

For the remainder of this section, we shall assume that $\nu=1$. This parameter can be reintroduced by using the fact that
\begin{gather*}
    \big(D^1+(a+1)E^0-\nu E^1\big)\circ\sigma_\nu=\nu \sigma_\nu\circ \big(D^1+(a+1)E^0-E^1\big);
\end{gather*}
c.f., \eqref{sigmaD}. In particular, this intertwining relation implies that
\begin{gather*}
    L_\lambda(\alpha,p_0,a,\nu)=\nu^{-|\lambda|} \sigma_\nu\big(L_\lambda(\alpha,p_0,a,1)\big).
\end{gather*}

\subsection{A generating function}
Proceeding as in the proof of Proposition \ref{HermiteGenFuncProp}, but using Proposition 4.1 rather than 3.1 in Baker and Forrester \cite{BF97}, it is straightforward to verify the generating function expansion given in the following proposition:

 \begin{Proposition}\label{PropFincGenLaguerre}
Let $q=1+(p_0-1)/\alpha$. Then, we have that
\begin{gather*}
    \sum_\lambda\frac{1}{h_\lambda [a+q]_\lambda}\frac{L_\lambda\otimes P_\lambda}{\epsilon_{p_0}(P_\lambda)}  = {}_0\mathscr{F}_1(a+q) e^{-1\otimes p_1}.
\end{gather*}
 \end{Proposition}

As detailed in Section \ref{SectionGenFunc} in the case of the Hermite symmetric functions (see the paragraph containing \eqref{DBCH}), by exploiting the Baker--Campbell--Hausdorf\/f formula, we can deduce an inf\/inite-dimensional family of eigenoperators for the Laguerre symmetric functions, parameterised by the shifted symmetric functions.

\begin{Proposition}\label{CoroLfL}
Fix $f\in\Lambda_{\mathbb{F},\alpha}$, and let $k$ be the degree of $f$.  Let, furthermore,
\begin{gather}
    \mathcal{L}^L_f = \mathcal{L}_f -\big\lbrack D^1+(a+1)E^0,\mathcal{L}_f \big\rbrack
    + \frac{1}{ 2!}\big\lbrack D^1+(a+1)E^0,\big\lbrack D^1+(a+1)E^0,\mathcal{L}_f\big\rbrack\big\rbrack\nonumber\\
\hphantom{\mathcal{L}^L_f =}{}
    +\cdots+\frac{(-1)^k}{ k!}\big\lbrack D^1+(a+1)E^0),\ldots,\big\lbrack D^1+(a+1)E^0,\mathcal{L}_f\big\rbrack\cdots\big\rbrack.\label{LaguerreEigenOps}
\end{gather}
Then, $\mathcal{L}^L$ is a differential operator on $\Lambda_{\mathbf{K}}$ of order $2k$. Moreover, it is the unique operator on~$\Lambda_{\mathbf{K}}$ such that
\begin{gather*}
    \mathcal{L}^L_f L^{}_\lambda = f(\lambda)L_\lambda
\end{gather*}
for all partitions $\lambda$.
\end{Proposition}

\subsection{Limit transition from the Laguerre to the Hermite symmetric functions}

As is well known for the one-variable case (see for instance \cite[Section 2.22]{KS96}), there is a simple limit transition from the Laguerre polynomials to the Hermite polynomials.  Explicitly,
\begin{gather}
\lim_{a\to\infty} \left( \frac{2}{a}\right)^{\frac{n}{2}}L^a\big(\sqrt{2a}x+a\big)=\frac{(-1)^n}{n!}H_n(x) ,
\end{gather}
where $L_n$ and $H_n$ respectively denote the standard (and consequently, non-monic) Laguere and Hermite polynomials.

Computer simulations suggest that a very similar transition relation holds between the Laguerre and Hermite symmetric functions.  It seems that the corresponding  transition  for any f\/inite number of variables greater than one, has been overlooked by the authors cited in the list of references.

\medskip

\noindent {\bf Conjecture.} \emph{For any partition $\lambda$ and generic values of the parameters $\alpha$ and  $\nu$, we have}
\begin{gather}\label{transitionLH}
    \lim_{a\to\infty} \left( \frac{1}{2a}\right)^{\frac{|\lambda|}{2}} \sigma_{\sqrt{2a}} \circ t_{a/\nu^{}} \big(L_\lambda(\alpha,a,\nu)\big)=H_\lambda\big(\alpha,\nu^2\big).
\end{gather}

One easily checks that the translation and the scale change  induced by $\sigma_{\sqrt{2a}} \circ t_{a/\nu}$ allow to map, as $a\to\infty$, the  Laguerre operator $\mathcal{L}^L$ and the Hermite  operator $\mathcal{L}^H$.  This is not suf\/f\/icient however for proving the conjecture.  The dif\/f\/iculty resides  in the fact that the Laguerre symmetric function $L_\Lambda(\alpha,a,\nu)$ has degree $|\lambda|$ in the parameter $a$, while $ t_{a/\nu}L_\Lambda(\alpha,a,\nu)$ appears to have degree $|\lambda|/2$.   Proving the latter property would require some new identities about the generalised binomial coef\/f\/icients.

\subsection{A limit from the Jacobi symmetric functions}
Also the Laguerre symmetric functions $L_\lambda(\alpha,p_0,a)$ can be viewed as limits of the Jacobi symmetric functions $\mathcal{J}_\lambda(\alpha,p_0,p,q)$. To make this precise, we set $a=-p-q-1/2$, and observe that
\begin{gather*}
    {q}^{-1}\sigma_{2q^{-1}}\circ\mathcal{L}^J = \big(D^1+(a+1)E^0-E^1\big)\circ \sigma_{2q^{-1}}+{q}^{-1}\big(D^2+(a+3/2)E^1\big)\circ \sigma_{2q^{-1}}.
\end{gather*}
From the representation \eqref{FRep} we can thus infer the following proposition:

\begin{Proposition}\label{LaguerreLimitProp}
Let $\lambda$ be a partition. Then, for generic parameter values, we have that
\begin{gather*}
    L_\lambda(\alpha,p_0,a) = \lim_{q\to\infty}\left({q}/{2}\right)^{|\lambda|}\sigma_{2q^{-1}}\big(\mathcal{J}_\lambda(\alpha,p_0,-a-q-1/2,q)\big)
\end{gather*}
in the sense of term-wise convergence.
\end{Proposition}

\subsection{Structure of Pieri formulae and invariant ideals}
The Laguerre case is in many ways easier to handle than the Hermite case. The reason being that, if we apply $\sigma_{2q^{-1}}$ to \eqref{JacobiPieri}, then the limit $q\to\infty$ is well-def\/ined. By Proposition \ref{LaguerreLimitProp}, this directly yields Pieri type recurrence relations for the Laguerre symmetric functions.

We f\/irst note the appropriate limit of the normalisation factors $\epsilon_0(\mathcal{J}_\lambda)$:
\begin{gather*}
    \lim_{q\to\infty}(-q/2)^{|\lambda|}\epsilon_0\big(\mathcal{J}_\lambda(\alpha,p_0,-q-a-1/2,q)\big) =\frac{C^0_\lambda(p_0/\alpha)C^0_\lambda((p_0-1)/\alpha+a+1)}{C^-_\lambda(1/\alpha)}.
\end{gather*}
In particular, this means that the left-hand side of the recurrence relations \eqref{JacobiPieri} for the Jacobi symmetric functions have the limits
\begin{gather*}
    \lim_{q\to\infty}(\sigma_{2q^{-1}})\left(\frac{\mathcal{J}_\lambda(\alpha,p_0,-q-a-1/2,q)}{\epsilon_0(\mathcal{J}_\lambda(\alpha,p_0,-q-a-1/2,q))}\right)=(-1)^{|\lambda|}\frac{\mathcal{L}_\lambda(\alpha,p_0,a)}{C^0_\lambda(p_0/\alpha)C^0_\lambda((p_0-1)/\alpha+a+1)},
\end{gather*}
where we have introduced the re-normalised Laguerre symmetric functions
\begin{gather*}
    \mathcal{L}_\lambda(\alpha,p_0,a) = C^-_\lambda(1/\alpha)L_\lambda(\alpha,p_0,a).
\end{gather*}

It is straightforward to verify that the $q\to\infty$ limit of $R_{\epsilon(L)}(\rho^J(p_0)+\lambda)$ is, just as in the Hermite case, given by \eqref{RLim}, while
\begin{gather*}
    \lim_{q\to\infty}(-q)^{|L|}\prod_{j\in L}\hat{w}^J\big(\epsilon_j(\rho^J(p_0)+\lambda)_j\big)
      =4^{|L|}\!\prod_{j\in L_+}\!\left(\frac{p_0-j}{\alpha}+\lambda_j+a+1\right)\!\prod_{j\in L_-}\! \left(\frac{p_0-j}{\alpha}+\lambda_j\right)
\end{gather*}
for $L=J,K$. In addition, we have that
\begin{gather*}
    \frac{C^0_\lambda(p_0/\alpha)C^0_\lambda((p_0-1)/\alpha+a+1)}{C^0_{\lambda+e_{J_+}-e_{J_-}}(p_0/\alpha)C^0_{\lambda+e_{J_+}-e_{J_-}}((p_0-1)/\alpha+a+1)}\\
    \qquad{}
    = \frac{\prod\limits_{j\in J_-}\big(p_0-j+1+\alpha(\lambda_j-1)\big)\big(p_0-j+\alpha(\lambda_j+a)\big)}{\prod\limits_{j\in J_+}\big(p_0-j+1+\alpha\lambda_j\big)\big(p_0-j+\alpha(\lambda_j+1+a)\big)}.
\end{gather*}
Now, proceeding in analogy with the proof of Proposition~\ref{HermitePieriFormulaeProp}, we arrive at the corresponding recurrence relations for the re-normalised Laguerre symmetric functions.

\begin{Proposition}\label{LaguerrePieriFormulaeProp}
The re-normalised Laguerre symmetric functions $\mathcal{L}_\lambda$ satisfy recurrence relations of the form
\begin{gather}
    (-1/2)^re_r\mathcal{L}_\lambda = \sum_{J_+,J_-}\hat{W}_{I(\ell(\lambda)+r);J_+,J_-}(\lambda)\mathcal{L}_{\lambda+e_{J_+}-e_{J_-}},
\end{gather}
where the sum extends over all subsets $J_+,J_-\subset\mathbb{N}$ such that $J_+\cap J_- = \varnothing$, $|J_+|+|J_-|\leq r$,
and $\lambda+e_{J_+}-e_{J_-}$ is a partition.

Moreover, the coefficients $\hat{W}_{I(\ell(\lambda)+r);J_+,J_-}$ are given explicitly by
\begin{gather*}
    \hat{W}_{I(\ell(\lambda)+r);J_+,J_-}(\lambda)  = (-1)^{r-p}\hat{V}_{I(\ell(\lambda)+r);J_+,J_-}(\lambda)  \sum_{K_+,K_-}\hat{V}_{I(\ell(\lambda)+r)\setminus(J_+\cup J_-);K_+,K_-}(\lambda)
\end{gather*}
with
\begin{gather*}
    \hat{V}_{I(\ell(\lambda)+r);J_+,J_-}(\lambda) = \prod_{j\in J_-}\left(\frac{p_0-j+1}{\alpha}+\lambda_j-1\right)\left(\frac{p_0-j}{\alpha}+\lambda_j+a\right)\\
    \qquad{}\times\prod_{j\in J_+,j^\prime\in J_-}\left(1+\frac{1}{j^\prime-j+\alpha(\lambda_j-\lambda_{j^\prime})}\right)\left(1+\frac{1}{j^\prime-j+\alpha(\lambda_j-\lambda_{j^\prime}+1)}\right)\\
    \qquad{} \times\prod_{j\in J_-}\big(\ell(\lambda)+r-j+\alpha\lambda_j\big)\prod_{i\in I(\ell(\lambda)+r)\setminus J}\left(1+\frac{1}{j-i+\alpha(\lambda_i-\lambda_j)}\right)\\
    \qquad{} \times\prod_{j\in J_+}\frac{1}{\ell(\lambda)+r+1-j+\alpha\lambda_j}\prod_{i\in I(\ell(\lambda)+r)\setminus J}\left(1-\frac{1}{j-i+\alpha(\lambda_i-\lambda_j)}\right),
\end{gather*}
and where the sum runs over all subsets $K_+,K_-\subset I(\ell(\lambda)+r)\setminus(J_+\cup J_-)$ such that $K_+\cap K_- = \varnothing$ and $|K_+|+|K_-|=r-|J_+|-|J_-|$.
\end{Proposition}

As a direct consequence, we have the following corollary:

\begin{Corollary}\label{LaguerrePieriFormulaeCor}
The re-normalised Laguerre symmetric functions $\mathcal{L}_\lambda$ satisfy the recurrence relation
\begin{gather*}
    (-1/2)e_1\mathcal{L}_\lambda = \sum_{j=1}^{|\lambda|+1}\big(\hat{W}_j(\lambda)\mathcal{L}_{\lambda+e_j} + \hat{W}_{-j}(\lambda)\mathcal{L}_{\lambda-e_j}\big)
\end{gather*}
with the coefficients
\begin{gather*}
    \hat{W}_j(\lambda)  = \frac{1}{\ell(\lambda)+2-j+\alpha\lambda_j}\prod_{\substack{1\leq i\leq|\lambda|+1\\ i\neq j}}\left(1-\frac{1}{j-i+\alpha(\lambda_i-\lambda_j)}\right),\\
    \hat{W}_{-j}(\lambda)  = \left(\frac{p_0-j+1}{\alpha}+\lambda_j-1\right) \left(\frac{p_0-j}{\alpha}+\lambda_j+a\right)\big(\ell(\lambda)+1-j+\alpha\lambda_j\big)\\
    \hphantom{\hat{W}_{-j}(\lambda)  =}{}
    \times\prod_{\substack{1\leq i\leq|\lambda|+1\\ i\neq j}}\left(1+\frac{1}{j-i+\alpha(\lambda_i-\lambda_j)}\right).
\end{gather*}
\end{Corollary}

\begin{Remark}
Corollary \ref{LaguerrePieriFormulaeCor} generalises results that were known in the polynomial case, namely  Proposition~4.7 of~\cite{BF97}.  Similarly, by restric\-ting the symmetric functions in Proposition~\ref{LaguerrePieriFormulaeProp} to a~f\/inite number of variables, one recovers Propositions~2.6 and~2.7 of~\cite{vD97}.
\end{Remark}

It is a straightforward exercise to adapt the discussion succeeding Corollary \ref{FirstHermitePieriCor}, as well as the proof of Theorem \ref{IdealThm}, to the Laguerre case. In this way, one arrives at the main result of this section.

\begin{Theorem}\label{IdealThmLaguerre}
$\Lambda_{\mathbb{F}}$ contains a non-trivial ideal invariant under all differential operators $\mathcal{L}^L_f$, $f\in\Lambda_{\mathbb{F},\alpha}$, for and only for non-zero $p_0$ of the form
\begin{gather}\label{p0ForIdealLaguerre}
    p_0 = n - \alpha m \qquad \text{or}\qquad p_0=n+1-\alpha(m+a+1),\qquad n,m\in\mathbb{N}_0.
\end{gather}
If that is the case, then there is a unique such ideal, spanned by the Laguerre symmetric func\-tions~$L_\lambda$ labelled by the partitions~$\lambda$ such that $(n+1,m+1)\in\lambda$.
\end{Theorem}

We note that the presence of the second family of ideals can be explained by the symmetry property in Proposition~\ref{propsymlaguerre}. Indeed, we have that $p_0'=n+1-\alpha(m+a^\prime+1)$.

\section{Deformed CMS operators and super polynomials}
We set the parameter $p_0=n-\alpha m$ for some $n,m\in\mathbb{N}_0$, and let $I_{n,m}\subset\Lambda_{\mathbb{F}}$ be the subspace given by
\begin{gather*}
    I_{n,m} = \mathbb{F}\big\langle P_\lambda:(n+1,m+1)\in\lambda\big\rangle.
\end{gather*}
Although it is not obvious from its def\/inition, $I_{n,m}$ is an ideal in $\Lambda_{\mathbb{F}}$. This fact follows immediately from well-known Pieri formulae for Jack's symmetric functions; see, e.g., Chapter VI in Macdonald \cite{Mac95}.

In this section we shall show that any CMS operator of the form \eqref{GenCMSOps} admits a restriction onto the quotient ring $\Lambda_{\mathbb{F}}/I_{n,m}$. Moreover, these restrictions will be given explicitly by dif\/ferential operators of so-called deformed CMS type.   Deformed CMS operators were f\/irst introduced by Chalykh, Feigen, Veselov \cite{CFV98} and later studied by many authors \cite{Fei08,HL10,KG05, Ser01,SV04,SV05}.  We shall consider the Hermite and Laguerre cases in some detail. In particular, this restriction procedure will lead to the introduction of  super Hermite and super Laguerre polynomials. Furthermore, the results we have obtained in the context of symmetric functions more or less immediately restrict to corresponding results for these `super' polynomials.

At this point, it is interesting to recall that, in the Hermite and Laguerre cases, we know, by Theorems~\ref{IdealThm} and~\ref{IdealThmLaguerre}, that there exists a unique ideal invariant under the action of all eigenoperators. We can thus conclude that the corresponding restrictions are given precisely by the deformed CMS operators mentioned above. In terms of the ideal $I_{n,m}$ this fact is equivalent to the equalities
\begin{gather}\label{idealEqualities}
    I_{n,m} = \mathbb{F}\big\langle H_\lambda:(n+1,m+1)\in\lambda\big\rangle = \mathbb{F}\big\langle L_\lambda:(n+1,m+1)\in\lambda\big\rangle.
\end{gather}
We stress that these equalities are non-trivial. For example, by def\/inition, the Hermite symmetric functions are of the form
\begin{gather*}
    H_\lambda = P_\lambda + \sum_{\mu\subset\lambda}u_{\lambda\mu}P_\mu,\qquad u_{\lambda\mu}\in\mathbb{F},
\end{gather*}
for $p_0=n-\alpha m$. The former of the above equalities thus implies that the coef\/f\/icient $u_{\lambda\mu}=0$ for all partitions $\lambda,\mu$ such that $(n+1,m+1)\in\lambda$ and $(n+1,m+1)\notin\mu$.

\subsection{Super Jack polynomials}
Our discussion below will involve a few results from the theory of super Jack polynomials which we now brief\/ly recall.
Let $x=(x_1,\ldots,x_n)$ and $y=(y_1,\ldots,y_m)$ be two sequences of indeterminates, and consider the subalgebra
\begin{gather*}
    \Lambda_{\mathbb{F},n,m}\subset\mathbb{F}\lbrack x_1,\ldots,x_n,y_1,\ldots,y_m\rbrack
\end{gather*}
consisting of all polynomials $p(x,y)$ that are separately symmetric in the variables $x$ and $y$, and that satisfy the condition
\begin{gather*}
    \left(\frac{\partial p}{\partial x_i} + \frac{1}{\alpha}\frac{\partial p}{\partial y_I}\right)\Bigg\arrowvert_{x_i=y_I} = 0
\end{gather*}
for all $i = 1,\ldots,n$ and $I = 1,\ldots,m$. This algebra generalises that of so-called supersymmetric polynomials,
which were f\/irst introduced as characters in the representation theory of the superalgebra $gl(n|m)$; see, e.g., Examples~23 and~24 in Section~I.3 of Macdonald~\cite{Mac95} or Moens and Van der Jeugt~\cite{MVdJ04}.
As shown by Sergeev and Veselov~\cite{SV04}, the algebra $ \Lambda_{\mathbb{F},n,m}$ is generated by the `deformed' power sums
\begin{gather*}
    p_{r,\alpha}(x,y) = \sum_{i=1}^n x_i^r - \alpha\sum_{I=1}^m y_I^r,\qquad r\in\mathbb{N}.
\end{gather*}
This yields a surjective homomorphism $\varphi_{n,m}= \Lambda_{\mathbb{F}}\rightarrow\Lambda_{\mathbb{F},n,m}$ by
\begin{gather*}
    \varphi_{n,m}(p_r) = p_{r,\alpha}(x,y),\qquad r\in\mathbb{N}.
\end{gather*}
Moreover, the kernel of this homomorphism is known to be spanned by the Jack's symmetric functions~$P_\lambda$ with $(n+1,m+1)\in\lambda$; see Theorem~2 in~\cite{SV05}.

It is clear from the facts listed above that $\varphi_{n,m}$ yields an isomorphism
\begin{gather*}
    \varphi_{n,m}: \Lambda_{\mathbb{F}}/I_{n,m}\overset\sim\longrightarrow\Lambda_{\mathbb{F},n,m},
\end{gather*}
and that the super Jack polynomials, def\/ined, as in Kerov et al.~\cite{KOO98}, by
\begin{gather*}
    SP_\lambda(x,y) = \varphi_{n,m}(P_\lambda),
\end{gather*}
form a basis for $\Lambda_{\mathbb{F},n,m}$ as $\lambda$ runs through all partitions such that $(n+1,m+1)\notin\lambda$. It is also clear that $\varphi_{n,m}$ intertwines between eigenoperators of Jack's symmetric functions and super Jack polynomials: for each $f\in\Lambda_{\mathbb{F},\alpha}$ there exists a unique dif\/ferential operator $\mathcal{L}_{f,n,m}$ on $\Lambda_{\mathbb{F},n,m}$ such that the diagram
\begin{gather}\label{JackDiagram}
    \begin{CD}
        \Lambda_{\mathbb{F}} @>\mathcal{L}_f>>\Lambda_{\mathbb{F}}\\
        @V\varphi_{n,m}VV @VV\varphi_{n,m}V\\
        \Lambda_{\mathbb{F},n,m} @>\mathcal{L}_{f,n,m}>> \Lambda_{\mathbb{F},n,m}
    \end{CD}
\end{gather}
is commutative. Indeed, the dif\/ferential operators $\mathcal{L}_{f,n,m}$ are given by their action on the super Jack polynomials:
\begin{gather*}
    \mathcal{L}_{f,n,m}SP_\lambda(x,y)=f(\lambda)SP_\lambda(x,y).
\end{gather*}
For further details see, in particular, Sergeev and Veselov \cite{SV05}.

We continue by observing that the duality relation \eqref{JackDuality} for Jack's symmetric functions has an interesting analogue for super Jack polynomials. This observation, which is made precise in the proposition below, is essentially due to the fact that
\begin{gather*}
    -\alpha p_{r,1/\alpha}(y,x) = p_{r,\alpha}(x,y).
\end{gather*}

\begin{Proposition}\label{DualitySuperJack}
The super Jack polynomials satisfy the duality relation
\begin{gather}\label{sJackDuality}
    SP_\lambda(\alpha;x,y)=(-1)^{|\lambda|}\,SQ_{\lambda'}(1/\alpha;y,x)
\end{gather}
for all partitions $\lambda$ such that $(n+1,m+1)\notin\lambda$, where
\begin{gather*}
   SQ_\lambda(x,y) = b_\lambda SP_\lambda(x,y);
\end{gather*}
c.f.~\eqref{QDef}.
\end{Proposition}

\begin{proof}
For clarity of exposition, we shall make explicit the dependence on $\alpha$ by writing $\varphi_{n,m}^{(\alpha)}$ for $\varphi_{n,m}$. In particular, this means that
 \begin{align*}
    \varphi_{m,n}^{(1/\alpha)}(p_r) = p_{r,1/\alpha}(y,x),\qquad r\in\mathbb{N}.
\end{align*}
We proceed to consider the action of the homomorphism $\varphi^{(1/\alpha)}_{m,n}\circ\omega_\alpha\circ\sigma_{-1}$ on one of Jack's symmetric function $P_\lambda$. Since $\omega_\alpha (p_r) =(-1)^{r-1}\alpha p_r$ and $\sigma_{-1}(p_r)=(-1)^rp_r$, we have that
\begin{gather*}
    \big(\varphi^{(1/\alpha)}_{m,n}\circ\omega_\alpha\circ\sigma_{-1}\big)(p_r)=p_{r,\alpha}(x,y),\qquad r\in\mathbb{N}.
\end{gather*}
We can thus conclude that
\begin{gather*}
    \varphi^{(1/\alpha)}_{m,n}\circ\omega_\alpha\circ\sigma_{-1}=\varphi_{n,m}^{(\alpha)}.
\end{gather*}
This yields the left-hand side of \eqref{sJackDuality}. On the other hand, it follows from \eqref{JackDuality}, and the fact that $P_\lambda$ is homogeneous of degree $|\lambda|$, that
\begin{gather*}
    \big(\varphi^{(1/\alpha)}_{m,n}\circ\omega_\alpha\circ\sigma_{-1}\big)(P_\lambda) = (-1)^{|\lambda|}\varphi_{m,n}^{1/\alpha}\big(Q_{\lambda^\prime}(1/\alpha)\big),
\end{gather*}
and we arrive at the right-hand side of \eqref{sJackDuality}.
\end{proof}

We conclude this section by showing that the duality relation for the super Jack polynomials implies a similar duality relation for the following `super' version of the hypergeometric series introduced in~\eqref{eqpFqoneset}:
 \begin{gather*}
    {}_pSF_q(a_1,\ldots,a_p;b_1,\ldots,b_q;\alpha;x,y) = \sum_\lambda
 \frac{1}{h_\lambda}\frac{\lbrack a_1\rbrack_\lambda\cdots\lbrack a_p\rbrack_\lambda}{\lbrack b_1\rbrack_\lambda\cdots\lbrack b_q\rbrack_\lambda} SP_\lambda(\alpha;x,y),
\end{gather*}
where the sum extends over all partition $\lambda$ such that $(n+1,m+1)\not\in\lambda$.

\begin{Proposition}
Let
\begin{gather*}
    \alpha'=1/\alpha,\qquad x'=(-\alpha)^{1+q-p}x,\qquad y'=(-\alpha)^{1+q-p}y,
\end{gather*}
and let
\begin{gather*}
    a_i'  = -\alpha a_i,\quad i=1,\ldots,p,\qquad
    b_j'  = -\alpha b_j,\quad j=1,\ldots,q.
\end{gather*}
Then, we have that
\begin{gather*}
    {}_pSF_q(a_1,\ldots,a_p;b_1,\ldots,b_q;\alpha;x,y) = {}_pSF_q\big(a_1',\ldots,a_p';b_1',\ldots,b_q';\alpha';y',x'\big).
\end{gather*}
\end{Proposition}

\begin{proof}
Applying Proposition \ref{DualitySuperJack} to the def\/inition of ${}_pSF_q$ yields
\begin{gather*}
    {}_pSF_q(a_1,\ldots,a_p;b_1,\ldots,b_q;\alpha;x,y) = \sum_\lambda
 \frac{(-1)^{|\lambda|}b_{\lambda'}(\alpha')}{h_\lambda(\alpha)}\frac{\lbrack a_1\rbrack^{(\alpha)}_\lambda\cdots\lbrack a_p\rbrack^{(\alpha)}_\lambda}{\lbrack b_1\rbrack^{(\alpha)}_\lambda
    \cdots\lbrack b_q\rbrack^{(\alpha)}_\lambda} SP_{\lambda'}(\alpha';y,x).
\end{gather*}
From \eqref{QDef} we infer that
\begin{gather*}
    b_\lambda(\alpha)=\frac{h_{\lambda'}(\alpha^\prime)}{\alpha^{|\lambda|}h_\lambda(\alpha)}=\frac{1}{b_{\lambda'}(\alpha')}.
\end{gather*}
Moreover, it is easily verif\/ied that
\begin{gather*}
    [x]_\lambda^{(\alpha)}=(-\alpha)^{-|\lambda|}[-\alpha x]_{\lambda'}^{(\alpha')}.
\end{gather*}
By applying the last two equations to the series above, we f\/ind that
\begin{gather*}
    {}_pSF_q(a_1,\ldots,a_p;b_1,\ldots,b_q;\alpha;x,y) \\
    \qquad{}= \sum_\lambda
 \frac{(-\alpha)^{(1+q-p)|\lambda|}}{h_{\lambda'}(\alpha')}\frac{\lbrack -\alpha a_1\rbrack^{(\alpha')}_\lambda\cdots\lbrack -\alpha a_p\rbrack^{(\alpha')}_\lambda}{\lbrack -\alpha b_1\rbrack^{(\alpha')}_\lambda
    \cdots\lbrack -\alpha b_q\rbrack^{(\alpha')}_\lambda} SP_{\lambda'}(\alpha';y,x),
\end{gather*}
which is clearly equivalent to the statement.
\end{proof}

\subsection{Deformed CMS operators}
We proceed to show that all CMS operators of the form \eqref{GenCMSOps} admit a restriction onto the quotient ring $\Lambda_{\mathbb{F}}/I_{n,m}\simeq\Lambda_{\mathbb{F},n,m}$, and that these restrictions are given by deformed CMS operators
\begin{gather*}
	\mathcal{L}_{n,m} = \sum_{k=0}^\infty a_kD^k_{n,m} + \sum_{\ell=0}^\infty b_\ell E^\ell_{n,m}
\end{gather*}
with only f\/initely many coef\/f\/icients $a_k,b_\ell\in\mathbb{F}$ non-zero, and where
\begin{gather}\label{EdefDef}
    E^\ell_{n,m} = \sum_{i=1}^nx_i^\ell\frac{\partial }{\partial x_i}+\sum_{I=1}^my_I^\ell\frac{\partial}{\partial y_I}
\end{gather}
and
\begin{gather}
    D^k_{n,m} = \sum_{i=1}^nx_i^k\frac{\partial^2 }{\partial x_i^2}+\frac{2}{\alpha}\sum_{i\neq j}\frac{x_i^k}{x_i-x_j}\frac{\partial}{\partial x_i} -\frac{1}{\alpha}\sum_{I=1}^my_I^k\frac{\partial^2 }{\partial y_I^2}-2\sum_{I\neq J}\frac{y_I^k}{y_I-y_J}\frac{\partial}{\partial y_I}\nonumber\\
\hphantom{D^k_{n,m} =}{}
    -2\sum_{i,I}\frac{1}{x_i-y_I} \left(x_i^k\frac{\partial}{\partial x_i}+\frac{1}{\alpha}y_I^k\frac{\partial}{\partial y_I}\right)-k\left(1+\frac{1}{\alpha}\right)\sum_{I=1}^my_I^{k-1}\frac{\partial}{\partial y_I}.\label{DdefDef}
\end{gather}
Clearly, it is suf\/f\/icient to consider the dif\/ferential operators $E^\ell$ and $D^k$; c.f.\ Lemma~\ref{EDLemma} and Proposition~\ref{CMSOpsProp}.

\begin{Proposition}\label{EDDefProp}
Let $p_0=n-\alpha m$ for some $n,m\in\mathbb{N}_0$. Then, for all $\ell,k\in\mathbb{N}_0$, the differential operators $E^\ell_{n,m}$ and $D^k_{n,m}$ preserve the algebra $\Lambda_{\mathbb{F},n,m}$. Moreover, they are the unique operators on $\Lambda_{\mathbb{F},n,m}$ such that the diagrams
\begin{subequations}
\begin{gather}\label{EdefProj}
    \begin{CD}
        \Lambda_{\mathbb{F}} @>E^\ell>>\Lambda_{\mathbb{F}}\\
        @V\varphi_{n,m}VV @VV\varphi_{n,m}V\\
        \Lambda_{\mathbb{F},n,m} @>E^\ell_{n,m}>> \Lambda_{\mathbb{F},n,m}
    \end{CD}
\end{gather}
and
\begin{gather}\label{DdefProj}
    \begin{CD}
        \Lambda_{\mathbb{F}} @>D^k>>\Lambda_{\mathbb{F}}\\
        @V\varphi_{n,m}VV @VV\varphi_{n,m}V\\
        \Lambda_{\mathbb{F},n,m} @>D^k_{n,m}>> \Lambda_{\mathbb{F},n,m}
    \end{CD}
\end{gather}
\end{subequations}
are commutative.
\end{Proposition}

\begin{proof}
In order to simplify the proof somewhat it is convenient to collect $x$ and $y$ into a sequence $z=(z_1,\ldots,z_{n+m})$ by setting
\begin{gather*}
    z_i= \begin{cases}
        x_i, & 1\leq i\leq n,\\
        y_{i-n}, & n+1\leq i\leq m.
    \end{cases}
\end{gather*}
Also, we introduce a map $\rho:\lbrace 1,\ldots,n+m\rbrace\to\mathbb{F}$ by specifying the value of $\rho(i)$ according to
\begin{gather*}
    \rho(i)= \begin{cases}
        1, & 1\leq i\leq n,\\
        -\alpha, & n+1\leq i\leq m.
    \end{cases}
\end{gather*}
It is clear that
\begin{gather*}
    p_{r,\alpha}(x,y) = \sum_{i=1}^{n+m}\rho(i)z_i.
\end{gather*}
Furthermore, it is readily verif\/ied that $E^\ell_{n,m}$ and $D^k_{n,m}$ are given by
\begin{gather*}
    E^\ell_{n,m} = \sum_{i=1}^{n+m}z_i^{\ell-1}\frac{\partial}{\partial z_i},
\end{gather*}
and
\begin{gather}\label{DkRho}
    D^k_{n,m} = \sum_{i=1}^{n+m}\frac{1}{\rho(i)}z_i^k\frac{\partial^2}{\partial z_i^2} + \frac{2}{\alpha}\sum_{i\neq j}\rho(j)\frac{z_i^k}{z_i-z_j}\frac{\partial}{\partial z_i} -k\sum_{i=1}^{n+m}\left(1-\frac{1}{\rho(i)}\right)z_i^{k-1}\frac{\partial}{\partial z_i},
\end{gather}
respectively.

Proceeding in analogy with the proof of Lemma \ref{EDLemma}, we observe that
\begin{gather*}
    \varphi_{n,m}(E^\ell p_r) = rp_{r,\alpha}(x,y) = E^\ell_{n,m}(\varphi_{n,m}p_r)
\end{gather*}
if we employ the convention $p_{0,\alpha}(x,y)\equiv n-\alpha m$. This yields \eqref{EdefProj}. Turning now to $D^k_{n,m}$, we note that
\begin{gather*}
    D^k_{n,m}\ p_{r,\alpha}(z)p_{s,\alpha}(z) = 2rs p_{r+s+k-2,\alpha}(z) +p_{r,\alpha}(z)\big(D^k_{n,m}\ p_{s,\alpha}(z)\big)+p_{s,\alpha}(z)\big(D^k_{n,m}\ p_{r,\alpha}(z)\big),
\end{gather*}
and deduce by direct computations that
\begin{gather*}
    2\sum_{i\neq j}\rho(j)\frac{z_i^k}{z_i-z_j}\frac{\partial}{\partial z_i}p_{r,\alpha}(z) = r\sum_{i\neq j}\sum_{m=0}^{r+k-2}\rho(i)\rho(j)z_i^{r+k-2}z_j^m\\
 \hphantom{2\sum_{i\neq j}\rho(j)\frac{z_i^k}{z_i-z_j}\frac{\partial}{\partial z_i}p_{r,\alpha}(z)}{}
     = r\left(\sum_{m=0}^{r+k-2}p_{r+k-2-m,\alpha}(z)p_{m,\alpha}(z)-\sum_{i=1}^{n+m}\rho(i)^2z_i^{r+k-2}\right).
\end{gather*}
Using the fact that
\begin{gather*}
    r(r-1)-r(r+k-1)\frac{\rho(i)^2}{\alpha}+rk(1-\rho(i)) = r(r-1)\rho(i)-\frac{1}{\alpha}r(r+k-1)\rho(i),
\end{gather*}
it is now straightforward to deduce that
\begin{gather*}
    D^k_{n,m}\ p_{r,\alpha}(z)p_{s,\alpha}(z) = r(r - 1)p_{r+k-2,\alpha}(z)p_{s,\alpha}(z) + 2rsp_{r+q+k-2,\alpha}(z)\\
     \qquad{} + s(s - s)p_{r,\alpha}(z)p_{s+k-2,\alpha}(z)  +\frac{r}{\alpha}\sum_{m=0}^{r+k-2}\big(p_{r+k-2-m,\alpha}(z)p_{m,\alpha}(z) - p_{r+k-2,\alpha}(z)\big).
\end{gather*}
Comparing this result with Def\/inition \ref{EDDef} we obtain \eqref{DdefProj}.

Finally, uniqueness is clear since $\varphi_{n,m}$ is surjective.
\end{proof}

The following lemma will be useful in establishing duality relations for the super Hermite and Laguerre polynomials:

\begin{Lemma} \label{dualityDnm}
We have that
\begin{gather*}
    D^k_{n,m}(\alpha)=-\frac{1}{\alpha}\big(D^{k}_{m,n}(1/\alpha)+k(1+\alpha)E^{k-1}_{m,n}\big).
\end{gather*}
\end{Lemma}

\begin{proof}
We f\/irst def\/ine
\begin{gather*}
    \Delta_{n,m}(\alpha)=\sum_{i,I}\frac{1}{x_i-y_I}
 \left(x_i^k\frac{\partial}{\partial x_i}+\frac{1}{\alpha} y_I^k\frac{\partial}{\partial y_I}\right).
\end{gather*}
Then, we can decompose $D^k_{n,m}$ as follows:
\begin{gather*}
    D^k_{n,m}(\alpha)=D^k_{n,x}(\alpha)-\frac{1}{\alpha}D^k_{m,y}(1/\alpha)-2\Delta_{n,m}(\alpha)-k\left(1+\frac{1}{\alpha}\right)E^{k-1}_{m,y},
\end{gather*}
where the subscript $x$ indicates that the operator in questions acts in the indeterminates $x$, and similarly for the subscript $y$. We observe that $\Delta_{n,m}(\alpha)=-(1/\alpha)\Delta_{m,n}(1/\alpha)$. Consequently, we can rewrite the operator  $D^k_{n,m}(\alpha)$ in the form
\begin{gather*}
    -\frac{1}{\alpha}\big(D^k_{m,x}(1/\alpha)-\alpha D^k_{n,y}(\alpha)-2\Delta_{m,n}(1/\alpha) +k(1+\alpha)E^{k-1}_{m,y}\big).
\end{gather*}
The desired formula now follows from the fact that $E_{m,n}^\ell=E_{n,m}^\ell$ and the fact that $E_{m,y}^{k-1}=E_{n,m}^{k-1}- E_{n,x}^{k-1}$.
\end{proof}

\subsection{Super Hermite polynomials}\label{sHermitePols}
As before, we assume that $p_0=n-\alpha m$ for some $n,m\in\mathbb{N}_0$. It directly follows from Proposi\-tion~\ref{EDDefProp} that
\begin{gather}
    \varphi_{n,m}\circ\mathcal{L}^H = \mathcal{L}^H_{n,m}\circ\varphi_{n,m}
\end{gather}
with $\mathcal{L}^H$ as in \eqref{HermiteOP} (for $\nu=1$) and $\mathcal{L}^H_{n,m}$ obtained by substituting $D^k_{n,m}$ for $D^k$ and $E^\ell_{n,m}$ for~$E^\ell$.\footnote{The operator $ \mathcal{L}^H_{n,m}$ corresponds to the operator $-\tilde H_{n,m}$ def\/ined in \cite{HL10} by equation~(44) and Case~I of Table~2 together with $\alpha=1/\kappa$, $\nu^2=\omega$, $x_i=z_i$, and $y_i=\tilde z_i$.  $\mathcal{L}^H_{n,m}$ can also be obtained from the operator in Theorem~9 in~\cite{Fei08} given by equation~(15) by f\/irst performing a similarity transformation by the function $\exp\big({-}\omega\sum\limits_{i=1}^{n+m}y_i^2\big)$, then setting $y_i=x_i$, $i=1,\ldots,n$, $y_{n+i}=y_i/\sqrt{k}$, $i=1,\ldots,m$, and f\/inally $k=-1/\alpha$ and $\omega=\nu^2$.} This observation suggests the following def\/inition:

\begin{Definition}
Let $\lambda$ be a partition such that $(n+1,m+1)\notin\lambda$. We then def\/ine the super Hermite polynomial $SH_\lambda(\alpha;x,y)$ by
\begin{gather*}
    SH_\lambda(x,y) = \varphi_{n,m}(H_\lambda).
\end{gather*}
\end{Definition}

Clearly, the super Hermite polynomials satisfy the eigenvalue equation
\begin{gather*}
    \mathcal{L}^H_{n,m}SH_\lambda(x,y) = -2|\lambda|SH_\lambda(x,y).
\end{gather*}
We note that the super Hermite polynomials coincide with the (reduced) eigenfunctions constructed for Case I in Section 4 of \cite{HL10}. This follows from the fact that, for f\/ixed $\lambda$, the eigenvalue equation above has a unique solution with a~triangular expansion in super Jack polynomials; c.f.~Theorem 4 in Appendix B and the discussion in Section~5.3 of~\cite{HL10}.

Using the results obtained in Section \ref{HermiteSec} we can also produce higher-order eigenoperators. For $f\in\Lambda_{\mathbb{F},\alpha}$, we let $\mathcal{L}^H_{f,n,m}$ denote the deformed CMS operator obtained from~\eqref{HermiteEigenOps} by substitu\-ting~$\mathcal{L}_{f,n,m}$ for $\mathcal{L}_f$ and $D^0_{n,m}$ for $D^0$.

\begin{Proposition}\label{deformedHermiteProp}
Let $p_0 = n - \alpha m$ for some $n,m\in\mathbb{N}_0$, and let $f\in\Lambda_{\mathbb{F},\alpha}$. Then, the deformed CMS operator $\mathcal{L}^H_{f,n,m}$ is the unique operator on $\Lambda_{\mathbb{F},n,m}$ such that the diagram
\begin{gather*}
    \begin{CD}
        \Lambda_{\mathbb{F}} @>\mathcal{L}^H_f>>\Lambda_{\mathbb{F}}\\
        @V\varphi_{n,m}VV @VV\varphi_{n,m}V\\
        \Lambda_{\mathbb{F},n,m} @>\mathcal{L}^H_{f,n,m}>> \Lambda_{\mathbb{F},n,m}
    \end{CD}
\end{gather*}
is commutative. Moreover, we have that
\begin{gather*}
    \mathcal{L}^H_{f,n,m}SH_\lambda(x,y) = f(\lambda)SH_\lambda(x,y).
\end{gather*}
\end{Proposition}

\begin{proof}\sloppy
The commutativity of the diagram, as well as uniqueness, follows immediately from~\eqref{JackDiagram}, Proposition~\ref{EDDefProp}, and the def\/inition of $\mathcal{L}^H_{f,n,m}$. The fact that the super Hermite polynomials are eigenfunctions of this operator, with the stated eigenvalues, is then a direct consequence of Corollary~\ref{HermiteEigenOpsCorollary}.
\end{proof}

An immediate consequence of the previous proposition is the integrability of the deformed CMS model in the Hermite case (also called rational CMS model with harmonic conf\/inement).  In the special case $n=1$ or $m=1$, the integrability of this model was already known \cite[Theorem~3]{CFV98}, and for $1/\alpha\in\mathbb{Z}$, but general values of~$n$ and~$m$, a dif\/ferent proof can be found in~\cite{Fei08}.

It also follows from Proposition \ref{deformedHermiteProp} that the kernel $\text{ker}(\varphi_{n,m})$ of $\varphi_{n,m}$ is invariant under all dif\/ferential operators $\mathcal{L}^H_f$, $f\in\Lambda_{\mathbb{F},\alpha}$. The uniqueness result in Theorem~\ref{IdealThm} implies that
\begin{gather*}
    \text{ker}(\varphi_{n,m}) = \mathbb{F}\big\langle H_\lambda:(n+1,m+1)\in\lambda\big\rangle,
\end{gather*}
and we thus arrive at the f\/irst equality in \eqref{idealEqualities}.

We proceed to deduce a duality relation for the super Hermite polynomials, analogous to that given in Proposition \ref{DualitySuperJack} for the super Jack polynomials. To this end, we should re-introduce the parameter $\nu^2$ via the homomorphism $\sigma_\nu$:
\begin{gather*}
    SH_\lambda(\alpha,\nu^2;x,y):= \varphi_{n,m}\big(\sigma_\nu(H_\lambda(\alpha))\big);
\end{gather*}
c.f., the paragraph containing \eqref{sigmaInter}. As for the Hermite symmetric functions, this duality relation can be directly inferred from a constructive def\/inition of the super Hermite polynomials. If we apply the homomorphism $\varphi_{n,m}$ to both sides of \eqref{HermiteExpRep}, then Proposition \ref{EDDefProp} implies the following result:

\begin{Proposition}\label{LassalleSuperHermite}
Let $\lambda$ be a partition such that $(n+1,m+1)\notin\lambda$, and set $L=\lfloor |\lambda|/2 \rfloor$. Then, we have that
\begin{gather*}
   SH_\lambda(x,y)=\exp_L\left(-\frac{1}{4\nu^2}D^{0}_{n,m}\right)SP_\lambda(x,y).
\end{gather*}
\end{Proposition}

Since $-\alpha D^0_{n,m}(\alpha)=D^0_{m,n}(1/\alpha)$ (c.f., Lemma \ref{dualityDnm}), this immediately yields the desired duality relation for the super Hermite polynomials.

\begin{Proposition}\label{PropDualitySuperHermite}
The super Hermite polynomials satisfy the duality relation
\begin{gather*}
    SH_\lambda\big(\alpha,\nu^2;x,y\big)=(-1)^{|\lambda|}b_{\lambda'}^{(1/\alpha)} SH_{\lambda'}\big(1/\alpha,-\alpha\nu^2;y,x\big).
\end{gather*}
\end{Proposition}

By applying the homomorphism $\varphi_{n,m}$, we can restrict further results from Section~\ref{HermiteSec} to the super Hermite polynomials. In particular, in this way we obtain the generating function expansion
\begin{gather*}
    \sum_\lambda \frac{1}{h_\lambda\epsilon_{n-\alpha m}(P_\lambda)} SH_\lambda(x,y)\otimes SP_\lambda(z,w) = {}_0\mathscr{\SF}_0(x,y;z,w)e^{-\frac{1}{4}p_{2,\alpha}(z,w)}
\end{gather*}
with
\begin{gather*}
    {}_0\mathscr{\SF}_0(x,y;z,w):= \sum_\lambda \frac{1}{h_\lambda\epsilon_{n-\alpha m}(P_\lambda)} SP_\lambda(x,y)SP_\lambda(z,w),
\end{gather*}
where, in both formulae, the sums extend over all partitions $\lambda$ such that $(n+1,m+1)\notin\lambda$. We leave the straightforward deduction of further such properties of the super Hermite polynomials to the interested reader.

\subsection{Super Laguerre polynomials}\label{sLaguerrePols}
The theory of super Laguerre polynomials can be developed in complete analogy with that of the super Hermite polynomials. First of all, we have that
\begin{gather}
    \varphi_{n,m}\circ\mathcal{L}^L = \mathcal{L}^L_{n,m}\circ\varphi_{n,m}
\end{gather}
with $\mathcal{L}^L$ as in \eqref{LaguerreOp} (for $\nu=1$) and $\mathcal{L}^L_{n,m}$ obtained by substituting $D^k_{n,m}$ for $D^k$ and $E^\ell_{n,m}$ for~$E^\ell$.\footnote{In~\cite{HL10}, take the operator    $-\tilde H_{n,m}/4$ given in equation~(44), then choose the Case~IV of Table~2, set $z_i=x_i$, $\tilde z_i=y_i$, $\omega=\nu$, and $\kappa=1/\alpha$.  This gives the operator $ \mathcal{L}^L_{n,m}$ def\/ined above. One can also obtain the same operator from the operator in Theorem~9 in~\cite{Fei08} given by equation~(18) in a similar manner to the Hermite case after f\/irst performing a change of variables $y_i\to\sqrt{y_i}$.}
This fact suggests the following def\/inition:

\begin{Definition}
Let $\lambda$ be a partition such that $(n+1,m+1)\notin\lambda$. We then def\/ine the super Laguerre polynomial $SL_\lambda(\alpha,a;x,y)$ by
\begin{gather*}
    SL_\lambda(x,y) = \varphi_{n,m}(L_\lambda).
\end{gather*}
\end{Definition}

In a similar manner to the Hermite case one can verify that the super Laguerre polynomials coincide with the (reduced) eigenfunctions constructed for Case IV in Section 4 of \cite{HL10}.

Higher-order eigenoperators can be constructed using Proposition \ref{CoroLfL}.

\begin{Proposition}\label{deformedLaguerreProp}
Let $p_0 = n - \alpha m$ for some $n,m\in\mathbb{N}_0$, and fix $f\in\Lambda_{\mathbb{F},\alpha}$. Then, the deformed CMS operator $\mathcal{L}^L_{f,n,m}$, obtained from \eqref{LaguerreEigenOps} by substituting $\mathcal{L}_{f,n,m}$ for $\mathcal{L}_f$, $D^1_{n,m}$ for $D^1$ and $E^0_{n,m}$ for $E^0$, is the unique operator on $\Lambda_{\mathbb{F},n,m}$ such that the diagram
\begin{gather*}
    \begin{CD}
        \Lambda_{\mathbb{F}} @>\mathcal{L}^L_f>>\Lambda_{\mathbb{F}}\\
        @V\varphi_{n,m}VV @VV\varphi_{n,m}V\\
        \Lambda_{\mathbb{F},n,m} @>\mathcal{L}^L_{f,n,m}>> \Lambda_{\mathbb{F},n,m}
    \end{CD}
\end{gather*}
is commutative. Moreover, we have that
\begin{gather*}
    \mathcal{L}^L_{f,n,m}SL_\lambda(x,y) = f(\lambda)SL_\lambda(x,y).
\end{gather*}
\end{Proposition}

We note that integrability of the deformed CMS model in Laguerre case (i.e., type B), which now directly follows from the above proposition, had been previously proved for the  case with $n=m=1$ \cite[Theorem~5]{CFV98}, and for $1/\alpha\in\mathbb{Z}$, but general values of $n$ and $m$, in \cite{Fei08}.

Reintroducing the parameter $\nu$, we can easily establish a duality relation using the following constructive def\/inition:

\begin{Proposition}
Let $\lambda$ be a partition and set $L=|\lambda|$.  Then, we have that
\begin{gather*}
    SL_\lambda = \exp_L\left(-\frac{1}{\nu}\big(D^1_{n,m}+(a+1)E^0_{n,m}\big)\right)SP_\lambda.
\end{gather*}
\end{Proposition}

A duality relation for the super Laguerre polynomials is now straightforward to infer.

\begin{Proposition}
We have that
\begin{gather*}
    SL_\lambda(\alpha,a,\nu;x,y)=(-1)^{|\lambda|}b_{\lambda'}^{(1/\alpha)}\,SL_{\lambda'}(1/\alpha,-\alpha a, -\alpha\nu;y,x).
\end{gather*}
\end{Proposition}

Further results can be obtained by applying $\varphi_{n,m}$ to results obtained in Section \ref{LaguerreSec}.

\appendix

\section{Dif\/ferential operators on the symmetric functions}\label{AppDiffOps}
To establish a wider context for our results, we proceed to consider the notion of dif\/ferential operators on the ring of symmetric functions. In particular, this will make it clear that the inf\/inite-dimensional CMS operators we have considered indeed are dif\/ferential operators on this ring.

Throughout this section, $\mathbb{F}$ will denote an arbitrary, but f\/ixed, f\/ield of characteristic zero. Correspondingly, we have $\Lambda_\mathbb{F}\equiv\mathbb{F}\otimes_{\mathbb{Z}}\Lambda$.

For a module over a commutative algebra, there is a natural notion of dif\/ferential operators on said module; see, e.g., Chapter~2 of Heyneman and Sweedler~\cite{HS69}. In the present case of $\Lambda_{\mathbb{F}}$, this general notion results in the following def\/inition: we f\/irst set $\mathscr{D}_{-1}(\Lambda_{\mathbb{F}}) = \lbrace 0\rbrace$, and then let
\begin{gather*}
    \mathscr{D}_k(\Lambda_{\mathbb{F}}) = \big\lbrace u\in\text{End}(\Lambda_{\mathbb{F}}): u\circ p - p\circ u\in\mathscr{D}_{k-1}(\Lambda_{\mathbb{F}})~\forall\, p\in\Lambda_{\mathbb{F}}\big\rbrace
\end{gather*}
for $k\geq 0$. Following standard terminology, we shall refer to the elements in the space $\mathscr{D}_k(\Lambda_{\mathbb{F}})\setminus\mathscr{D}_{k-1}(\Lambda_{\mathbb{F}})$ as the (linear and homogeneous) dif\/ferential operators on $\Lambda_{\mathbb{F}}$ of order $k$. Moreover, we shall write $\mathscr{D}(\Lambda_{\mathbb{F}})$ for the space of all such dif\/ferential operators, i.e., $\mathscr{D}(\Lambda_{\mathbb{F}}) = \cup_{k\geq -1} \mathscr{D}_k(\Lambda_{\mathbb{F}})$.

We record the following two elementary facts: $\mathscr{D}_0(\Lambda_{\mathbb{F}}) = \Lambda_{\mathbb{F}}$ (if we identify any $p\in\Lambda_{\mathbb{F}}$ with the corresponding operator $q\mapsto pq$), for any $u\in\mathscr{D}_k(\Lambda_{\mathbb{F}})$ and $v\in\mathscr{D}_{k^\prime}(\Lambda_{\mathbb{F}})$, the composition $u\circ v\in\mathscr{D}_{k+k^\prime}(\Lambda_{\mathbb{F}})$. Both of these facts are easy to verify directly, and can also be inferred from the general theory mentioned above.

For our purposes, we require a more explicit realisation of $\mathscr{D}(\Lambda_{\mathbb{F}})$, given by a set of generators. We thus note that, since the power sums $p_r$ freely generate $\Lambda_{\mathbb{F}}$, we can, for each $r\in\mathbb{N}$, def\/ine $\partial(p_r)\in\mathscr{D}_1(\Lambda_{\mathbb{F}})$ by $\partial(p_r)1 = 0$ and
\begin{gather}\label{partialprDef}
    \partial(p_r)p_s =
    \begin{cases}
        1, & r=s,\\
        0, & r\neq s.
    \end{cases}
\end{gather}
The assignment $p_r\mapsto\partial(p_r)$ extends uniquely to an $\mathbb{F}$-algebra homomorphism $\partial: \Lambda_{\mathbb{F}}\to\mathscr{D}(\Lambda_{\mathbb{F}})$, given by
\begin{gather*}
    \partial\left(\sum_\lambda a_\lambda p_\lambda\right) = \sum_\lambda a_\lambda \partial(p_{\lambda_1})\partial(p_{\lambda_2})\cdots.
\end{gather*}
More generally, we can consider the subring of $\mathscr{D}(\Lambda_{\mathbb{F}})$ generated by the dif\/ferential opera\-tors~$\partial(p_r)$ over $\Lambda_{\mathbb{F}}$:
\begin{gather*}
    A := \mathbb{F}\lbrack p_1,p_2,\ldots,\partial(p_1),\partial(p_2),\ldots\rbrack\subset\mathscr{D}(\Lambda_{\mathbb{F}}).
\end{gather*}
In the case of the f\/initely generated algebra $\Lambda_{\mathbb{F},n}$ the corresponding subring would coincide with~$\mathscr{D}(\Lambda_{\mathbb{F},n})$. However, this is not the case for $A$ and $\mathscr{D}(\Lambda_{\mathbb{F}})$. In fact, none of the dif\/ferential operators we focus on in this paper are contained in $A$. For example, it is easy to verify that the formal series
\begin{gather*}
    E_1 = \sum_{r=1}^\infty p_r\partial(p_r)
\end{gather*}
def\/ines a f\/irst-order dif\/ferential operator on $\Lambda_{\mathbb{F}}$ by
\begin{gather}\label{E1Acts}
    E_1q\equiv \sum_{r=1}^\infty p_r\partial(p_r)q,\qquad q\in\Lambda_{\mathbb{F}}.
\end{gather}
Indeed, since any $q\in\Lambda_{\mathbb{F}}$ can be written uniquely as a f\/inite linear combination of terms of the form $p_\lambda=p_{\lambda_1}p_{\lambda_2}\cdots$, the inf\/inite sum in \eqref{E1Acts} contains only a f\/inite number of non-zero terms. Hence, $E_1q\in\Lambda_{\mathbb{F}}$, and the fact that $E_1\in\mathscr{D}_1(\Lambda_{\mathbb{F}})$ is now clear from the def\/inition of $\mathscr{D}_1(\Lambda_{\mathbb{F}})$.

We proceed to enlarge the ring $A$ such that we obtain all of $\mathscr{D}(\Lambda_{\mathbb{F}})$. The above example indicates that we should include also formal power series in the dif\/ferential operators $\partial(p_r)$ with coef\/f\/icients in $\Lambda_{\mathbb{F}}$. However, we should ensure that we do not include formal series of inf\/inite order. To illustrate this point, let us consider
\begin{gather}\label{infOrder}
    D = \sum_{k=1}^\infty\partial(p_1)^k.
\end{gather}
Interchanging the order of summation and dif\/ferentiation, as in \eqref{E1Acts}, this formal series does def\/ine an element in $\text{End}(\Lambda_{\mathbb{F}})$. But, suppose that $D\in\mathscr{D}_l(\Lambda_{\mathbb{F}})$ for some $l\in\mathbb{N}_0$. By def\/inition, this means that
\begin{gather*}
    \lbrack\cdots\lbrack\lbrack D,q_1\rbrack,q_2\rbrack,\cdots,q_{l+1}\rbrack=0
\end{gather*}
for any $q_1,\ldots,q_{l+1}\in\Lambda_{\mathbb{F}}$. However, if we set $q_1=\cdots=q_{l+1}=p_1$, then we obtain
\begin{gather*}
    \sum_{k=0}^\infty\frac{(l+k)!}{k!}\partial(p_1)^k,
\end{gather*}
which is clearly non-zero. Hence, $D\notin\mathscr{D}(\Lambda_{\mathbb{F}})$.

In order to make these remarks precise, we note that $A$ inherits a natural f\/iltration from $\mathscr{D}(\Lambda_{\mathbb{F}})$, given by the order of the dif\/ferential operators:
\begin{gather*}
    A^0\subset A^1\subset\cdots\subset A^k\subset\cdots
\end{gather*}
with
\begin{gather*}
    A^k=A\cap\mathscr{D}_k(\Lambda_{\mathbb{F}}) = \left\lbrace \sum_{\ell(\lambda)\leq k}q_\lambda\partial(p_\lambda): q_\lambda\in\Lambda_{\mathbb{F}}\right\rbrace.
\end{gather*}
Indeed, we clearly have $A=\cup_{k\geq 0}A^k$ and $A^kA^l\subset A^{k+l}$ for all $k,l\in\mathbb{N}_0$. On each such submodule $A^k\subset A$, we can then introduce a f\/iltration
\begin{gather*}
    A^k = A^k_0\supset A^k_1\supset\cdots\supset A^k_n\supset\cdots,
\end{gather*}
where
\begin{gather*}
    A^k_n=\left\lbrace \sum_{\substack{\ell(\lambda)\leq k\\ |\lambda|\geq n}}q_\lambda\partial(p_\lambda):q_\lambda\in\Lambda_{\mathbb{F}}\right\rbrace.
\end{gather*}
By $\hat{A}^k$ we denote the corresponding completion of $A^k$; see, e.g., Chapter 10 in Atiyah and Macdonald \cite{AM69} for a general discussion of the process of completion in the context of (commutative) algebra. For the discussion below, it will be important to note that each element in $\hat{A}^k$ can be uniquely identif\/ied with a formal power series in the f\/irst-order dif\/ferential operators $\partial(p_r)$ of order not greater than $k$, and vice versa. We proceed to show that the f\/iltered module
\begin{gather*}
    \hat{A} = \bigcup_{k\geq 0}\hat{A}^k,
\end{gather*}
can be naturally identif\/ied with $\mathscr{D}(\Lambda_{\mathbb{F}})$.

\begin{Proposition}\label{diffOpsProp}
Any element $D = \sum_\lambda q_\lambda\partial(p_\lambda)\in\hat{A}$ defines a differential operator on $\Lambda_{\mathbb{F}}$ by
\begin{gather*}
    Dp\equiv \sum_\lambda q_\lambda\partial(p_\lambda) p,\qquad p\in\Lambda_{\mathbb{F}}.
\end{gather*}
In this sense, $\hat{A}^k\setminus\hat{A}^{k-1} = \mathscr{D}_k(\Lambda_{\mathbb{F}})\setminus\mathscr{D}_{k-1}(\Lambda_{\mathbb{F}})$ for all $k\in\mathbb{N}_0$, and consequently $\hat{A} = \mathscr{D}(\Lambda_{\mathbb{F}})$.
\end{Proposition}

\begin{proof}
Given $D\in\hat{A}$ as in the statement and $p\in\Lambda_{\mathbb{F}}$, it is clear that $Dp$ is a f\/inite linear combination of $q_\lambda\in\Lambda_{\mathbb{F}}$. Hence, $D\in\text{End}(\Lambda_{\mathbb{F}})$. By induction on $k$ it is now easy to verify that $\hat{A}^k\subset\mathscr{D}_k(\Lambda_{\mathbb{F}})$.

We proceed to establish injectivity, i.e., that given $D,D^\prime\in\hat{A}$ we have $(D - D^\prime)p = 0$ for all $p\in\Lambda_{\mathbb{F}}$ if and only if $D = D^\prime$. The f\/irst part of this claim clearly holds true. Suppose therefore that $D\neq D^\prime$. This means that
\begin{gather*}
    D - D^\prime = \sum_\lambda q_\lambda\partial(p_\lambda)
\end{gather*}
with at least one coef\/f\/icient $q_\lambda\neq 0$. Among all such non-zero coef\/f\/icients, f\/ix one with the corresponding partition $\lambda$ having minimal length $\ell(\lambda)$. It follows that
\begin{gather*}
    (D - D^\prime)p_\lambda = q_\lambda\neq 0,
\end{gather*}
which proves the second part of the claim.

In order to establish surjectivity, i.e., that any $D\in\mathscr{D}_k(\Lambda_{\mathbb{F}})\setminus\mathscr{D}_{k-1}(\Lambda_{\mathbb{F}})$ is of the form $D = \sum\limits_{\ell(\lambda)=k} q_\lambda\partial(p_\lambda)\in\hat{A}^k\setminus\hat{A}^{k-1}$. where $A_{-1}:=\lbrace 0\rbrace$, we proceed by induction on the order $k$. For $k=0$, this claim is obvious. Suppose that $k>0$. Fix $D\in\mathscr{D}_k(\Lambda_{\mathbb{F}})\setminus\mathscr{D}_{k-1}(\Lambda_{\mathbb{F}})$, and let
\begin{gather*}
    f_r = \lbrack D,p_r\rbrack,\qquad r\in\mathbb{N}.
\end{gather*}
Then, by the induction assumption, we have that
\begin{gather*}
    f_r\in\mathscr{D}_{k-1}(\Lambda_{\mathbb{F}})\setminus\mathscr{D}_{k-2}(\Lambda_{\mathbb{F}}) = \hat{A}^{k-1}\setminus\hat{A}^{k-2}.
\end{gather*}
It follows that
\begin{gather*}
    D^\prime:= \frac{1}{k}\sum_r f_r\partial(p_r)\in\hat{A}_k\setminus\hat{A}_{k-1}.
\end{gather*}
We note that $\mathscr{D}(\Lambda_{\mathbb{F}})$ is a Lie algebra (with Lie bracket given by $\lbrack u,v\rbrack = u\circ v - v\circ u$); c.f., Lemma~2.1.1. in Heyneman and Sweedler~\cite{HS69}. By the Jacobi identity, we thus infer that
\begin{gather*}
    \lbrack f_r,p_s\rbrack\equiv \lbrack\lbrack D,p_r\rbrack,p_s\rbrack = \lbrack\lbrack D,p_s\rbrack,p_r\rbrack\equiv \lbrack f_s,p_r\rbrack,\qquad \forall\, r,s\in\mathbb{N}.
\end{gather*}
We also note that, for any partition $\lambda$,
\begin{gather*}
    \sum_r\lbrack\partial(p_\lambda),p_r\rbrack = \ell(\lambda).
\end{gather*}
Since each $f_r = \sum\limits_{\ell(\lambda)=k-1}q_\lambda\partial(p_\lambda)$ for some $q_\lambda\in\Lambda_{\mathbb{F}}$, these facts combine to yield
\begin{gather*}
    \lbrack D - D^\prime, p_s\rbrack = f_s - \frac{1}{k}f_s - \frac{1}{k}\sum_r\lbrack f_s,p_r\rbrack \partial(p_r) = 0.
\end{gather*}
Hence, $D-D^\prime\in\mathscr{D}_0(\Lambda_{\mathbb{F}}) = \hat{A}^0$, and surjectivity follows.
\end{proof}

Finally, we provide a simple criterion for a dif\/ferential operator in $\mathscr{D}(\Lambda_{\mathbf{F}})$ to be continuous with respect to the topology given by the ideal $U=\oplus_{k\geq 1}\Lambda_{\mathbf{F}}^k$; c.f.\ the paragraph preceding Proposition~\ref{prop2F1twosets}. To this end, we introduce a notion of degree, $\mathrm{deg}:\mathscr{D}(\Lambda_{\mathbb{F}})\to\mathbb{Z}\cup\lbrace-\infty\rbrace$, by requiring that
\begin{gather*}
    \mathrm{deg}(D)\geq m \qquad \mathrm{if} \quad DU^n\subset U^{n+m} \quad \forall\, n\in\mathbb{N}_0,
\end{gather*}
where we make the identif\/ication $U^{n+m}\equiv U^0$ for $m+n<0$, and setting $\mathrm{deg}(D)=-\infty$ if no such integer $m$ exists. It is clear that we have the following lemma:

\begin{Lemma}\label{degreeLemma}
Let $D\in\mathscr{D}(\Lambda_{\mathbb{F}})$. If $\mathrm{deg}(D)>-\infty$, then $D$ is continuous with respect to the $U$-adic topology.
\end{Lemma}

As a simple example of a dif\/ferential operator $D\in\mathscr{D}(\Lambda_{\mathbb{F}})$ with $\mathrm{deg}(D)=-\infty$, we note
\begin{gather*}
    D = \sum_{r=1}^\infty\partial(p_r).
\end{gather*}
Indeed, this follows directly from the fact that $Dp_r=1$ for all $r\in\mathbb{N}_0$.

\section{CMS operators on the symmetric functions}\label{AppCMSOps}
In this appendix we record a few technical details on the dif\/ferential operators $E^\ell$ and $D^k$ that are used throughout the paper. In addition, we shall isolate certain results that hold true not only in the Hermite and Laguerre cases, including the fact that a generic inf\/inite-dimensional CMS operator of second order has a complete set of eigenfunctions in $\Lambda_{\mathbf{F}}$.

As in the case of a f\/inite number of variables, these operators are not algebraically independent of each other. By direct computations, we deduce the following relations:

\begin{Lemma}\label{relationsLemma}
For $\ell,k\in\mathbb{N}_0$,
\begin{subequations}
\begin{gather}
\label{comEE}
    \lbrack E^k,E^{\ell+1}\rbrack  = (\ell+1) E^{k+\ell} ,\\
\label{lowerk}
    \lbrack E^0,D^{k+1}\rbrack  = (k+1) D^{k} ,\\
    \label{equalk}
    \lbrack E^1,D^{k}\rbrack  = (k-2) D^k ,\\
\label{raisek}
    \lbrack E^2,D^k\rbrack  = (k - 4)D^{k+1} + 2\left(\frac{p_0 - 1}{\alpha}-1\right)E^k,\\
    \label{Etop}
    \lbrack E^k,p_{\ell+1}\rbrack  = (\ell+1) p_{k+\ell},\\
    \lbrack D^k,p_{\ell+1} \rbrack  = 2(\ell+1) E^{k+\ell}+\ell(\ell+1)p_{k+\ell-1} + \frac{\ell+1}{\alpha}\sum_{m=0}^{k+\ell-1}(p_{k+\ell-m-1}p_m-p_{k+\ell-1}).\label{DtoE}
\end{gather}
\end{subequations}
\end{Lemma}

When constructing eigenfunctions of the CMS operators
\begin{gather}\label{genCMSOp}
    \mathcal{L} = \sum_{k=0}^2 a_k D^k + \sum_{\ell=0}^1 b_\ell E^\ell,\qquad a_k,b_\ell\in\mathbf{F},
\end{gather}
it is important to know their action on Jack's symmetric functions. (We could of course use another linear basis for $\Lambda_{\mathbf{F}}$, but the action on Jack's symmetric functions is particularly simple.) Clearly, it is suf\/f\/icient to consider the operators $E^\ell$ and $D^k$.

To this end, we recall that
Lassalle \cite{Las90} def\/ined generalised binomial coef\/f\/icients $\binom{\lambda}{\mu}$ by the series expansion
\begin{gather}\label{binomFormula}
    \frac{P_\lambda(x_1+1,\ldots,x_n+1)}{P_\lambda(1^n)} = \sum_{\mu\subseteq\lambda} \binom{\lambda}{\mu}\frac{P_\mu(x_1,\ldots,x_n)}{P_\mu(1^n)}.
\end{gather}
It was later shown that these binomial coef\/f\/icients are independent of $n$; see, e.g., \cite{OO97}. As a~consequence, we can deduce the following generalisation to Jack's symmetric functions:

\begin{Proposition}\label{genBinomProp}
Let $t_\gamma$ be the translation homomorphism defined in~\eqref{tgammaDef}. For any parti\-tion~$\lambda$, we have
\begin{gather}\label{genBinomFormula}
    \frac{t_1(P_\lambda)}{\epsilon_{p_0}(P_\lambda)} = \sum_{\mu\subseteq\lambda}\binom{\lambda}{\mu}\frac{P_\mu}{\epsilon_{p_0}(P_\mu)}.
\end{gather}
\end{Proposition}

\begin{proof}
We f\/ix the partition $\lambda$, and let $k=|\lambda|$. Then, we expand the dif\/ference between the left- and right-hand side of \eqref{genBinomFormula} in terms of Jack's symmetric functions,
\begin{gather*}
    \frac{t_1(P_\lambda)}{\epsilon_{p_0}(P_\lambda)} - \sum_{\mu\subseteq\lambda}\binom{\lambda}{\mu}\frac{P^{ }_\mu}{\epsilon_{p_0}(P_\mu)} = \sum_{|\mu|\leq|\lambda|}a_{\lambda\mu}(p_0)P_\mu.
\end{gather*}
By applying the restriction homomorphism $\rho_n$ we infer from \eqref{binomFormula} that $a_{\lambda\mu}(n)=0$ as long as $n\geq k$. (If $n<k$, then $P_\mu(x_1,\ldots,x_n)$ migh be zero.) Moreover, it is clear from Stanley's specialisation formula \eqref{epsilonXJackEq} and the def\/inition of $t_1$ that $a_{\lambda\mu}(p_0)$ is a rational function of $p_0$. Since it vanishes at inf\/initely many distinct points, we can thus conclude that $a_{\lambda\mu}\equiv 0$.
\end{proof}

Similarly, starting from the formula in Section 3 of Lassalle \cite{Las90}, we can deduce an expression for the lowest degree Pieri formula for Jack's symmetric functions given in terms of binomial coef\/f\/icients.

\begin{Proposition}\label{PieriJackProp}
For any partition $\lambda$, we have
\begin{gather*}
    p_1\frac{P_\lambda}{h_\lambda} = \sum_i \binom{\lambda^{(i)}}{\lambda}\frac{P_{\lambda^{(i)}}}{h_{\lambda^{(i)}}},
\end{gather*}
where the sum extends over all positive integers $i$ such that $\lambda^{(i)}$ is a partition.
\end{Proposition}

Moreover, the coef\/f\/icients in this formula are known explicitly: $h_\lambda$ is given by \eqref{defhook} and
\begin{gather}\label{binomExpr}
    \binom{\lambda^{(i)}}{\lambda}=\left(\lambda_i+1+\frac{\ell(\lambda^{(i)})-i}{\alpha}\right)\prod_{j\neq i}\frac{\alpha(\lambda_i+1-\lambda_j)+j-i-1}{\alpha(\lambda_i+1-\lambda_j)+j-i};
\end{gather}
see Section 14 in Lassalle \cite{Las98} for the latter fact.

Using the two Propositions above it is now straightforward to compute the action of the operators $E^\ell$ and $D^k$ on Jack's symmetric functions. In particular, we have the following:

\begin{Lemma}\label{actionLemma}
For any partition $\lambda$,
\begin{subequations}
\begin{gather}
\label{E2Action}
    E^2 \frac{P_\lambda }{h_\lambda} = \sum_i \binom{\lambda^{(i)}}{\lambda}\left(\lambda_i -\frac{i-1}{\alpha}\right)\frac{P_{\lambda^{(i)}}}{h_{\lambda^{(i)}}},\\
\label{E1Action}
    E^1 P_\lambda  = |\lambda|P^{}_\lambda,\\
\label{E0Action}
    E^0 \frac{P^{ }_\lambda}{\epsilon_{p_0}(P_\lambda)}  = \sum_i \binom{\lambda}{\lambda_{(i)}}\frac{P_{\lambda_{(i)}}}{\epsilon_{p_0}\big(P_{\lambda_{(i)}}\big)},\\
\label{D2Action}
    D^2 P_\lambda  = d_\lambda P_\lambda,\qquad d_\lambda = \sum_i \lambda_i\left(\lambda_i - 1 + \frac{2}{\alpha}(p_0 - i)\right),\\
\label{D1Action}
    D^1 \frac{P_\lambda}{\epsilon_{p_0}(P_\lambda)}  = \sum_i \binom{\lambda}{\lambda_{(i)}}\left(\lambda_i - 1 + \frac{p_0-i}{\alpha}\right)\frac{P_{\lambda_{(i)}}}{\epsilon_{p_0}\big(P_{\lambda_{(i)}}\big)},\\
\label{D0Action}
    D^0 \frac{P_\lambda}{\epsilon_{p_0}(P_\lambda)}  = \sum_{i,j} \binom{\lambda}{\lambda_{(i)}}\binom{\lambda_{(i)}}{\lambda_{(i,j)}}\left(\lambda_i-\lambda_j + \frac{j-i}{\alpha}+\delta_{ij}\right)
    \frac{P_{\lambda_{(i,j)}}}{\epsilon_{p_0}\big(P_{\lambda_{(i,j)}}\big)}.
\end{gather}
\end{subequations}
\end{Lemma}

\begin{proof}
We shall verify the equations in the order in which they are listed. According to \eqref{DtoE},
\begin{gather*}
    E^2=\frac{1}{2}[D^2,p_1]-\frac{1}{\alpha}(p_0-1)p_1.
\end{gather*}
Using this expression for $E^2$, \eqref{E2Action} is readily inferred from Proposition \ref{PieriJackProp} and \eqref{D2Action}, which is established below. Equation \eqref{E1Action} is a direct consequence of the def\/inition of $E^1$ and the fact that $P_\lambda$ is homogeneous of degree $|\lambda|$. We observe that
\begin{gather*}
    t_1(P_\lambda) = P_\lambda + E^0 P_\lambda + \text{l.d.}
\end{gather*}
(where l.d.~stands for terms of lower degree); c.f., \eqref{tgammaDef} for $\gamma=1$. Inserting this expression into the left-hand side of \eqref{genBinomFormula}, and then comparing coef\/f\/icients with the right-hand side, yields~\eqref{E0Action}.

It is clear from the discussion in Section~\ref{JackFuncsSection} and Lemma~\ref{EDLemma} that $D:=D^2-(2/\alpha)(p_0-1)E^1$ is the inverse limit of $D_n$, as def\/ined in \eqref{stableOp}. Consulting Example~3 in Section~VI.4 of Mac\-do\-nald~\cite{Mac95}, we thus conclude that $P_\lambda$ is an eigenfunction of $D$ with eigenvalue $\sum\limits_{i=1}^{\ell(\lambda)} \lambda_i(\lambda_i-1-2(i-1)/\alpha)$. Hence, \eqref{D2Action} follows from~\eqref{E1Action}. The remaining two equations \eqref{D1Action}, \eqref{D0Action} can now be deduced by a direct computation using~\eqref{lowerk}.
\end{proof}

\begin{Remark}
Combining the relations in Lemma \ref{relationsLemma} with Pieri formulae for Jack's symmetric functions, we could, in principle, compute the action of the dif\/ferential operators $E^\ell$ or $D^k$ for any $\ell,k\in\mathbb{N}$. However, the the resulting formulae become more and more complex as the values of $\ell$ and $k$ are increased. Since we shall only make use of the formulae obtained in Lemma \ref{actionLemma}, we therefore refrain from any further such computations.
\end{Remark}

As a direct application of Lemma \ref{actionLemma}, we have the following theorem:

\begin{Theorem}  \label{TheoExistenceEigenfunction}
In \eqref{genCMSOp} fix $a_k,b_\ell\in\mathbf{F}$ such that $a_2$ and $b_1$ are not both zero. Then, for any partition $\lambda$, there exists a unique symmetric function ${F}_\lambda\in\Lambda_{\mathbf{F}}$ such that
\begin{enumerate}\itemsep=0pt
\item $ {F}_\lambda=P_\lambda+\sum\limits_{\mu\subset \lambda}u_{\lambda\mu}  {P}_\lambda$ for some $u_{\lambda\mu}\in\mathbf{F}$ $($triangularity$)$;
\item $\mathcal{L}  {F}_\lambda=\varepsilon_\lambda  {F}_\lambda$ for some $\epsilon_\lambda\in\mathbf{F}$ $($eigenfunction$)$;
\end{enumerate}
Moreover, the eigenvalue $\epsilon_\lambda$ is given explicitly by
\begin{gather*}
    \varepsilon_\lambda=a_2\sum_i\lambda_i\left(\lambda_i-1+\frac{2}{\alpha}(p_0-i)\right)+b_1|\lambda|.
\end{gather*}
\end{Theorem}

\begin{proof}
According to Lemma \ref{actionLemma}, we have that
\begin{gather}\label{eqtheotriangularaction}
    \mathcal{L}( P_\lambda)= \varepsilon_\lambda P_\lambda+\sum_{\mu\subset\lambda}c_{\lambda\mu}P_\mu,\qquad \varepsilon_\lambda =a_2d_\lambda+b_1|\lambda|,
\end{gather}
for some coef\/f\/icients $c_{\lambda\mu}$, and with $d_\lambda$ as specif\/ied in~\eqref{D2Action}. Furthermore, it follows from~\eqref{defhook}, \eqref{epsilonXJackEq} and \eqref{binomExpr} that $c_{\lambda\mu}\in\mathbf{F}$. Let us now make the ansatz $F_\lambda=\sum\limits_{\mu\subseteq\lambda}u_{\lambda\mu}P_\mu$ with $u_{\lambda\lambda}\equiv 1$. By a direct computation, we then f\/ind that the eigenvalue equation $\mathcal{L}F_\lambda=\epsilon_\lambda F_\lambda$ holds true if and only if the coef\/f\/icients $u_{\lambda\mu}$ satisfy the recurrence relation
\begin{gather}\label{recurrenceRel}
    (\epsilon_\lambda-\epsilon_\mu)u_{\lambda\mu} = \sum_{\mu\subset\nu\subseteq\lambda}u_{\lambda\mu}c_{\nu\mu}.
\end{gather}
We observe that
\begin{gather*}
    \epsilon_\lambda-\epsilon_\mu  = a_2\sum_i\big(\lambda_i^2-\mu_i^2\big) + (b_1-a_2)(|\lambda|-|\mu|) + \frac{2a_2}{\alpha}\sum_i(\lambda_i-\mu_i)(p_0-i),
\end{gather*}
which clearly is non-zero for all $\mu\subset\lambda$. Since we have f\/ixed $u_{\lambda\lambda}\equiv 1$, this means that the coef\/f\/icients $u_{\lambda\mu}$ are uniquely determined by the recurrence relation \eqref{recurrenceRel}. Finally, the fact that $c_{\lambda\mu},\epsilon_\lambda\in\mathbf{F}$ implies that also $u_{\lambda\mu}\in\mathbf{F}$.
\end{proof}

Fix $m\in\mathbb{N}$, and let $\text{Par}_m$ denote the set of partitions $\lambda$ of weight $|\lambda|\leq m$. It is clear from~(2) in Theorem~\ref{TheoExistenceEigenfunction} that there exists a transition matrix $(M_{\lambda\mu})$ from the symmetric functions $\lbrace F_\lambda\rbrace_{\lambda\in\text{Par}_m}$ to the set of Jack's symmetric functions $\lbrace P_\lambda\rbrace_{\lambda\in\text{Par}_m}$, given by
\begin{gather*}
    F_\lambda = \sum_\mu M_{\lambda\mu}P_\mu.
\end{gather*}
Fix a total order $<_t$ on $\text{Par}_m$ that is compatible with the order given by inclusion of diagrams, i.e., if $\mu\subset\lambda$, then $\mu<_t\lambda$. If we order the entries of $(M_{\lambda\mu})$ according to this total order, then Property (2) in Theorem \ref{TheoExistenceEigenfunction} implies that we obtain a lower triangular matrix with one's on the diagonal. Hence, $(M_{\lambda\mu})$ is invertible, i.e., Jack's symmetric functions can be expressed as linear combinations of the $F_\lambda$. Since Jack's symmetric functions form a basis for $\Lambda_{\mathbb{F}}$, we thus arrive at the following corollary:

\begin{Corollary}\label{basisCor}
Let $a_k$ and $b_\ell$ be as in Theorem~{\rm \ref{TheoExistenceEigenfunction}}. Then, as $\lambda$ runs through the set of all partitions, the symmetric functions $F_\lambda$ form a basis for $\Lambda_{\mathbf{F}}$.
\end{Corollary}

We note that the eigenfunctions in Theorem~\ref{TheoExistenceEigenfunction} have the following useful representation:
\begin{gather}\label{FRep}
    F_\lambda=\prod_{\mu\subset\lambda}\frac{\mathcal{L}-\varepsilon_\mu}{\varepsilon_\lambda-\varepsilon_\mu}(P_\lambda).
\end{gather}
To establish this formula, we f\/irst note that the triangularity property~(1) can be directly inferred from Lemma~\ref{actionLemma}. In addition, from the proof of Theorem~\ref{TheoExistenceEigenfunction} we recall the triangular action~\eqref{eqtheotriangularaction} of $\mathcal{L}$ on $P_\lambda$. It follows that the dif\/ferential operator $\prod\limits_{\mu\subset\lambda}(\mathcal{L}-\epsilon_\mu)$ annihilates the subspace spanned by Jack's symmetric functions $P_\mu$ with $\mu\subset\lambda$. Clearly, this fact implies the eigenfunction property~(2).

Since we make use of it in Section~\ref{HermiteSec}, we also note the following lemma:

\begin{Lemma}\label{tCommuteLemma}
For $\ell,k\in\mathbb{N}_0$,
\begin{gather*}
    t_\gamma\circ E^\ell  = \left(\sum_{m=0}^\ell \gamma^{\ell-m}\binom{\ell}{m}E^m\right)\circ t_\gamma,\qquad
    t_\gamma\circ D^k = \left(\sum_{m=0}^k \gamma^{k-m}\binom{k}{m}D^m\right)\circ t_\gamma.
\end{gather*}
\end{Lemma}

\begin{proof}
We f\/irst observe that the analogous statement for f\/initely many variables is easily verif\/ied. Fix $n\in\mathbb{N}$, and def\/ine a homomorphism $t_{\gamma,n}:\Lambda_{\mathbb{F},n}\to\Lambda_{\mathbb{F},n}$ by setting
\begin{gather*}
    t_{\gamma,n}(p)(x_1,\ldots,x_n)=p(x_1+\gamma,\ldots,x_n+\gamma),\qquad \forall\, p\in\Lambda_{\mathbb{F},n}.
\end{gather*}
We note the intertwining relation $\varphi_n\circ t_\gamma=t_{\gamma,n}\circ\varphi_n$; c.f., the discussion succeeding \eqref{tgammaDef}. Clearly, we have that
\begin{gather*}
    t_{\gamma,n}\circ E^\ell_n = \left(\sum_{i=1}^n(x_i+\gamma)^\ell\frac{\partial}{\partial x_i}\right)\circ t_{\gamma,n} = \left(\sum_{m=0}^\ell\gamma^{\ell-m}\binom{\ell}{m}E^m_n\right)\circ t_{\gamma,n}.
\end{gather*}
From Lemma \ref{EDLemma} we thus infer that
\begin{gather*}
    \varphi_n\circ t_\gamma\circ E^\ell  = t_{\gamma,n}\circ E^\ell_n\circ\varphi_n = \left(\sum_{m=0}^\ell\gamma^{\ell-m}\binom{\ell}{m}E^m_n\right)\circ t_{\gamma,n}\circ\varphi_n \\
    \hphantom{\varphi_n\circ t_\gamma\circ E^\ell}{}
    = \varphi_n\circ\left(\sum_{m=0}^\ell \gamma^{\ell-m}\binom{\ell}{m}E^m\right)\circ t_\gamma.
\end{gather*}
Finally the fact that this equation holds true for all $n\in\mathbb{N}$ implies the statement for $E^\ell$; c.f., the proof of Lemma~\ref{EDLemma}. The dif\/ferential operator $D^k$ can be treated similarly.
 \end{proof}

\section{Proof of Theorem \ref{JacobiPieriFormulaeThm}}\label{JacobiPieriFormulaeThmProof}
The starting point is a sequence of recurrence relations for the Jacobi symmetric polynomials~$\mathcal{J}_\lambda(x)$, as deduced by van Diejen \cite{vD99} (see Theorem 6.4), which we now recall. In doing so, we shall essentially employ the formulation in Section~2 of Sergeev and Veselov~\cite{SV09}.

\begin{Theorem}[van Diejen \protect{\cite{vD99}}]
For generic parameter values, the generalised Jacobi polynomials $\mathcal{J}_\lambda(x)$ satisfy the recurrence relations
\begin{gather}\label{JacobiPolsRecurrenceRels}
    2^re_r(x)\frac{\mathcal{J}_\lambda(x)}{\mathcal{J}_\lambda(0^n)} = \!\sum_{\epsilon(J),\epsilon(K)} \! (-1)^{|K|}\hat{V}^{(+)}_{I(n),\epsilon(J)}(\rho^J(n)+\lambda)\hat{V}^{(-)}_{J^c,\epsilon(K)}
    (\rho^J(n)+\lambda)\frac{\mathcal{J}_{\lambda+e_{\epsilon(J)}}(x)}{\mathcal{J}_{\lambda+e_{\epsilon(J)}}(0^n)},\!\!\!
\end{gather}
where the sum is over all sequences of signs $\epsilon(J)$ and $\epsilon(K)$ with $J,K\subseteq I(n)\equiv\lbrace 1,\ldots,n\rbrace$ such that $J\cap K=\varnothing$, $|J| + |K| = r$, and $\lambda+e_{\epsilon(J)}$ is a partition.
\end{Theorem}

In order deduce the corresponding recurrence relations for the Jacobi symmetric functions it is important to know how~\eqref{JacobiPolsRecurrenceRels} depends on the number of variables~$n$. To this end, we f\/irst observe that the sum over~$\epsilon(J)$ is limited by the requirement that $\lambda+e_{\epsilon(J)}$ should be a partition, and consequently does not depend in any essential way on~$n$. However, there is no such obvious limitation on the sum over $\epsilon(K)$. In addition, it is not a priori clear how the coef\/f\/icients in~\eqref{JacobiPolsRecurrenceRels}, as well as the specialization $\mathcal{J}_\lambda(0^n)$, depend on~$n$. The f\/irst problem was resolved by Sergeev and Veselov \cite{SV09} (see Lemma~4.1) through the following lemma:

\begin{Lemma}[Sergeev and Veselov~\protect{\cite{SV09}}]
Let $J,K\subseteq I(m)\equiv\lbrace 1,\ldots,m\rbrace$ for some positive integer $m\geq\ell(\lambda)$, and $\epsilon(J)$ a corresponding conf\/iguration of signs. If $\lambda+\epsilon(J)$ is a partition, and $\max (K)\geq\ell(\lambda)+|K|+1$ and $\max (K)\geq\max (J)+|K|+1$, then $\hat{V}^{(-)}_{J^c,\epsilon(K)}(\rho^J+\lambda) = 0$.
\end{Lemma}

Let $r = |J| + |K|$. In order for $\lambda+e_{\epsilon(J)}$ to be a partition, it is clear that we must have $\max(J)\leq\ell(\lambda)+|J|$. It follows that the inequalities in the Lemma are satisf\/ied for all $K\subseteq I$ such that $\text{max}(K)\geq \ell(\lambda)+r+1$. Hence, we can restrict our attention to $J,K\subseteq I(\ell(\lambda)+r)$ irrespective of the specif\/ic value of $n$.

We proceed to consider the second problem: the dependence of the coef\/f\/icients $\hat{V}^{(+)}_{I,\epsilon(J)}(\rho^J+\lambda)$ and $\hat{V}^{(-)}_{J^c,\epsilon(K)}(\rho^J+\lambda)$ on $n$. We shall require a somewhat more detailed resolution of this problem than that stated by Sergeev and Veselov \cite{SV09} in their Lemma 4.2.

\begin{Lemma}
Let $n\in\mathbb{N}$ be such that $n\geq\ell(\lambda)+r+1$. Then,
\begin{subequations}
\begin{gather}\label{VplusExpr}
    \hat{V}^{(+)}_{I(n),\epsilon(J)}\big(\rho^J(n)+\lambda\big) = \hat{V}^{(+)}_{I(\ell(\lambda)+r),\epsilon(J)}\big(\rho^J(n)+\lambda\big)R_{\epsilon(J)}\big(\rho^J(n)+\lambda;\ell(\lambda)+r\big),
\\
\label{VminusExpr}
    \hat{V}^{(-)}_{J^c,\epsilon(K)}\big(\rho^J(n)+\lambda\big) = \hat{V}^{(-)}_{I(\ell(\lambda)+r)\setminus J,\epsilon(K)}\big(\rho^J(n)+\lambda\big)R_{\epsilon(K)}\big(\rho^J(n)+\lambda;\ell(\lambda)+r\big).
\end{gather}
\end{subequations}
\end{Lemma}

\begin{proof}
As previously observed, $\lambda+e_{\epsilon(J)}$ is a partition only if $J\subseteq I(\ell(\lambda)+r)$. It follows that
\begin{gather*}
    I\setminus J = (I(\ell(\lambda)+r)\setminus J)\cup \lbrace \ell(\lambda)+r+1,\ldots,n\rbrace.
\end{gather*}
We observe that, for any $m\in\mathbb{N}$ such that $\ell(\lambda)<m\leq n$,
\begin{gather*}
    \prod_{j\in J}\prod_{i=m}^n \hat{v}^J\big(\epsilon_j(\rho(n)+\lambda)_j + (\rho(n)+\lambda)_i\big)\hat{v}^J\big(\epsilon_j(\rho(n)+\lambda)_j - (\rho(n)+\lambda)_i\big)\\
    \qquad{} = R_{\epsilon(J)}\big(\rho(n)+\lambda;m\big),
\end{gather*}
which thus depends on $n$ only through $\rho(n)$. Using these facts, it is straightforward to infer equation~\eqref{VplusExpr} from the explicit def\/inition of the function $\hat{V}^{(+)}$ in~\eqref{Vplus}. The validity of~\eqref{VminusExpr} follows similarly once it is observed that, since also $K\subseteq I(\ell(\lambda)+r)$,
\begin{gather*}
    J^c\setminus K = (I(\ell(\lambda)+r)\setminus J)\setminus K\cup \lbrace \ell(\lambda)+r+1,\ldots,n\rbrace.\tag*{\qed}
\end{gather*}
\renewcommand{\qed}{}
\end{proof}

We recall that a formula for the specialisation of $\mathcal{J}_\lambda(x)$ at $x = 0^n$ can be obtained from Corollary 5.2 in Opdam \cite{Op89} by specialising to the root system $BC_n$ (and taking into account the relation between the $\mathcal{J}_\lambda$ and the multivariable Jacobi polynomials considered by Opdam; see Beerends and Opdam~\cite{BO93} and Sergeev and Veselov~\cite{SV09}). Although his formula does not directly generalise to the Jacobi symmetric functions, Sergeev and Veselov~\cite{SV09} (see Proposition~2.3) showed that it is given by the right hand side of~\eqref{JacobiSpecialisation} for $p_0 = n$. If we combine this observation with the two Lemmas above, then we can rewrite van Diejen's recurrence relations~\eqref{JacobiPolsRecurrenceRels} in the form~\eqref{JacobiPieri} for $p_0 = n$. Since the sum does not depend on $n$, and the coef\/f\/icients are rational functions of $n$, the validity of Theorem~\ref{JacobiPieriFormulaeThm} follows.

\subsection*{Acknowledgements}

The work of P.D.\ was  supported by FONDECYT grant \#1090034 and by CONICYT through the Anillos de Investigaci\'on RED4 y ACT56.
M.H.\ would like to thank O.~Chalykh, S.N.M.~Ruij\-se\-naars and A.P.~Veselov for helpful discussions.

\LastPageEnding

\end{document}